\documentclass[onefignum,onetabnum]{siamart190516}
\usepackage{bm}
\usepackage{enumitem}

\usepackage{tikz-cd} 

\usepackage{subcaption}
\usepackage{float}
\usepackage[textsize=small,backgroundcolor=white]{todonotes}

\definecolor{highyellow}{RGB}{255,251,194}
\definecolor{armygreen}{rgb}{0.29, 0.33, 0.13}
\definecolor{dgblue}{rgb}{.2, .2, .702}
\definecolor{dgred}{rgb}{1,0.5020,0}

\usetikzlibrary{shapes,decorations}
\tikzstyle{highbox} = [draw=red, fill=highyellow, very thick,
    rectangle, rounded corners, inner sep=12pt, inner ysep=12pt]
\tikzstyle{hightitle} =[fill=red, text=white]

\newcommand{\stext}[1]{{\small{\text{#1}}}}

\DeclareMathOperator{\tr}{tr}
\DeclareMathOperator{\spn}{span}
\DeclareMathOperator{\ran}{ran}

\newcommand{\Inputs}{\textbf{Inputs}}
\newcommand{\Outputs}{\textbf{Outputs}}
\newcommand{\Steps}{\textbf{Steps}}
\newcommand{\Return}{\textbf{Return:} }

\crefname{enumi}{Assumption}{Assumptions}

\usepackage{lipsum}
\usepackage{amsfonts}
\usepackage{graphicx}
\usepackage{epstopdf}
\usepackage[noend]{algpseudocode}
\ifpdf
  \DeclareGraphicsExtensions{.eps,.pdf,.png,.jpg}
\else
  \DeclareGraphicsExtensions{.eps}
\fi

\newsiamremark{remark}{Remark}
\newsiamremark{hypothesis}{Hypothesis}
\crefname{hypothesis}{Hypothesis}{Hypotheses}
\newsiamthm{claim}{Claim}

\headers{Quantum Mechanics for Closure of Dynamical Systems}{David Freeman, Dimitrios Giannakis}

\title{Quantum Mechanics for Closure of Dynamical Systems\thanks{Submitted to the editors DATE.
\funding{D.G.\ acknowledges support from the US National Science Foundation under grants 1842538 and DMS-1854383, the US Office of Naval Research under MURI grant N00014-19-1-242, and the US Department of Defense, Basic Research Office under Vannevar Bush Faculty Fellowship grant N00014-21-1-2946. D.C.F.\ is supported as a PhD student under the last grant.}}}

\author{David C.\ Freeman\thanks{Department of Mathematics, Dartmouth College, Hanover, NH 03755, USA (\email{david.c.freeman.gr@dartmouth.edu}).} \and Dimitrios Giannakis\thanks{Department of Mathematics, Dartmouth College, Hanover, NH 03755, USA; Department of Physics and Astronomy, Dartmouth College, Hanover, NH 03755, USA (\email{dimitrios.giannakis@dartmouth.edu}).} \and Joanna Slawinska\thanks{Department of Mathematics, Dartmouth College, Hanover, NH 03755; USA (\email{joanna.m.slawinska@dartmouth.edu}).}}

\usepackage{amsopn}
\DeclareMathOperator{\diag}{diag}
\DeclareMathOperator{\Id}{Id}

\ifpdf
\hypersetup{
  pdftitle={Quantum Mechanics for Closure of Dynamical Systems},
  pdfauthor={D. C. Freeman, D. Giannakis, J. Slawinska}
}
\fi

\begin{document}

\maketitle

\begin{abstract}
    We propose a scheme for data-driven parameterization of unresolved dimensions of dynamical systems based on the mathematical framework of quantum mechanics and Koopman operator theory. Given a system in which some components of the state are unknown, this method involves defining a surrogate system in a time-dependent quantum state which determines the fluxes from the unresolved degrees of freedom at each timestep. The quantum state is a density operator on a finite-dimensional Hilbert space of classical observables and evolves over time under an action induced by the Koopman operator. The quantum state also updates with new values of the resolved variables according to a quantum Bayes' law, implemented via an operator-valued feature map. Kernel methods are utilized to learn data-driven basis functions and represent quantum states, observables, and evolution operators as matrices. The resulting computational schemes are automatically positivity-preserving, aiding in the physical consistency of the parameterized system. We analyze the results of two different modalities of this methodology applied to the Lorenz 63 and Lorenz 96 multiscale systems, and show how this approach preserves important statistical and qualitative properties of the underlying chaotic systems.   
\end{abstract}



\section{Introduction}

Among the foundational problems in the modeling of complex dynamical systems is the question of how to account for fine-grain degrees of freedom which are too computationally complex to model directly. The Earth's climate system is a classical example of a multiscale, multiphysics system where direct numerical simulation of all relevant degrees of freedom is not computationally feasible (now, and for the foreseeable future \cite{Fox-KemperEtAl14}), necessitating the use of subgrid-scale models to represent unresolved degrees of freedom. For example, cloud formation (a highly influential process on climate scales) is in part determined by microscopic chemical processes and turbulent convective motions in the atmosphere across the entire globe. Representations of convective cloud physics within global climate models (GCMs) has thus relied on surrogate models of small-scale processes  to approximate their aggregate contribution over the larger spatiotemporal climate scales \cite{GrabowskiEtAl19}. This is a methodology known as \emph{closure}, or \emph{parameterization}, and besides climate dynamics it finds applications in many disciplines dealing with complex time-dependent phenomena, e.g., \cite{Sagaut06,EEtAl07}. 

In this paper, we present a new framework for closure of dynamical systems that models the unresolved degrees of freedom as a quantum mechanical system. Our approach extends a recently developed operator-theoretic framework for data assimilation, called quantum mechanical data assimilation (QMDA) \cite{Giannakis19b,FreemanEtAl22}, to the setting of two-way coupling between classical and quantum systems representing the resolved and unresolved dynamics, respectively.       

\subsection{\label{secParameterizationIntro}Parameterization}

Early approaches to parameterization, e.g., \cite{ArakawaSchubert74,GentMcWilliams89}, were based on low-order bulk formulas representing an average flux from unresolved degrees of freedom to the resolved variables. These formulas are typically constructed using physical reasoning, and feature a small number of parameters that can be tuned, e.g., using observational data, so as to best match the behavior of nature, or a high-resolution reference model, on the resolved scales.  More recently, with the advent of the ``big data'' era, data-driven parameterization schemes have received considerable attention \cite{BrenowitzBretherton18,BoltonZanna19,YuvalOGorman20}. These approaches employ training data generated by high-resolution, targeted simulations to learn models of the unresolved flux terms as functions of the resolved variables. The premise of these machine learning approaches is that by fitting closure models in a high-dimensional hypothesis space one can discover functional relationships between the resolved and unresolved variables that would be difficult to do on the basis of physical reasoning and/or asymptotic analysis alone. Indeed, this approach has been shown to outperform conventional parameterizations in terms of capturing the statistical behavior of the resolved variables in a variety of settings \cite{BrenowitzBretherton18,BoltonZanna19,YuvalOGorman20}.

A common feature of parameterization schemes, including the methods outlined above, is dimension reduction. That is, the contribution of the unresolved variables is made to be a function only of larger-scale variables which can be feasibly simulated or consistently measured. Effectively, this means that the equations governing the parameterized small-scale processes lose whatever independence they had from the resolved processes. In the paper \cite{Palmer01}, Palmer uses ideas from dynamical systems theory (in particular, the Poincar\'e-Bendixson theorem) to show that in some circumstances, dimension-reducing parameterizations will necessarily yield a system which is not chaotic. Palmer points out that since the chaotic nature of the climate is a vital aspect of its long-term behavior, a loss of chaotic dynamics could be a confounding factor in long-term statistical fidelity. In other words, the (partial) independence of small-scale processes could be a fundamental component of the qualitative behavior of large-scale processes. Simplifying the small-scale subsystems with functional dependencies exclusively in terms of other variables could thus result in the loss of dynamical complexity in the climate model, and substantially change the qualitative behavior of long-term climate predictions.

Palmer proposed a solution to this issue, in which unresolved scales were parameterized by a stochastic variable rather than a deterministic model. Altering the system to have a stochastic component effectively corresponds to state space augmentation; that is, the state space of the parameterized system now includes the resolved variables and the event space associated with the stochastic variables. As a result, the act of approximating some component of the state space no longer makes the parameterized variable totally determined by other components of the system through a functional relationship, which means there is no longer necessarily an erasure of the chaotic behavior. The probability density function of the stochastic variable in Palmer's example was based on some known information about the system (making its stochastic behavior mimic, in some sense, the deterministic behavior of the fully defined component). 

More broadly, stochastic parameterization schemes provide a means of restoring the independence between resolved and unresolved degrees of freedom, while maintaining computational tractability, by representing the unresolved degrees of freedom as realizations of a stochastic process \cite{MajdaEtAl99,LinNeelin00,Wilks05,CrommelinVanden-Eijnden08,ChorinLu15,BernerEtAl17}. In systems exhibiting timescale separation, a common route to stochastic parameterization is to derive stochastic differential equations (SDEs) for the unresolved variables by taking homogenization limits of the primitive governing equations \cite{PavliotisStuart08,MelbourneStuart11,KellyMelbourne17}. Alternatively, and particularly in systems where the unresolved dynamics are not consistent with SDEs, stochastic closure schemes are constructed by parameter inference from training data \cite{SchneiderEtAl21}. Recent works have used ideas from data assimilation \cite{MajdaHarlim12,LawEtAl15} to sequentially learn the parameters of stochastic closure models from noisy partial observations \cite{ChenLi21,GottwaldReich21,ChenEtAl22,LevineStuart22}.      

In placing these methods and results in context, it should be kept in mind that stochasticity of the parameterized dynamics is not a necessary ingredient for successful subgrid-scale modeling---rather, it is the act of representing unresolved degrees of freedom by surrogate dynamical models (of either deterministic or stochastic nature) that can lead to improved performance over parameterization schemes based on pure functional relationships with the resolved variables. Indeed, in multiscale dynamical systems exhibiting averaging principles \cite{PavliotisStuart08,ArtsteinEtAl07}, the effective slow dynamics is deterministic. The heterogeneous multiscale method \cite{EEtAl07} is a general multiscale modeling framework that leverages averaging principles to build microscale models for unresolved variables---in a broad range of applications, these models are deterministic. In atmospheric dynamics, one of the most successful approaches to parameterization, called super-parameterization, involves coupling a coarse atmospheric model with embedded column models of moist convection in each gridbox \cite{GrabowskiSmolarkiewicz99}, and these models are again oftentimes (but not always \cite{GroomsMajda13}) deterministic.

\subsection{Quantum mechanical data assimilation}  

QMDA is a technique for sequential data assimilation (filtering) of partially observed dynamical systems \cite{MajdaHarlim12,LawEtAl15}, combining elements of Koopman operator theory \cite{EisnerEtAl15} with the Dirac--von Neumann quantum mechanical evolution and measurement axioms \cite{Takhtajan08}. In QMDA, the density operator from quantum mechanics is employed as a generalization of the probability distribution for the system state in Bayesian data assimilation, the assimilated observables are represented by self-adjoint multiplication operators, and dynamical evolution of observables takes place under the unitary action of the Koopman operator. When observations are made, the density is updated via a projective von Neumann measurement, which is analogous to the Bayesian analysis step. Thus, under the assumption that the assimilated signal is governed by a classical dynamical system, QMDA comprises a one-way-coupled classical-quantum system, where the classical system influences the state of the quantum system through the observation map but the quantum system does not influence the state of the classical system.   

QMDA addresses a number of challenges in classical data assimilation, as it avoids the need for ad hoc Gaussian approximations, does not require diffusion regularization, and automatically preservers sign-definite quantities. Being rooted in linear operator theory, the method also has an asymptotically consistent, data-driven formulation, whereby operators are represented by matrices in a data-driven basis learned from time-ordered data using kernel algorithms \cite{BerryEtAl15,Giannakis19,BerryEtAl20,DasEtAl21}.

\subsection{Our contributions}

\begin{figure}
    \sffamily
    \begin{tikzpicture}
        \node[highbox] (box){%
            \begin{tikzcd}
                \stext{Classical} & [-1.4em] x_n \arrow[r, "\text{classical evolution}"] & [8em] x_{n+1} \arrow[r, "\text{Id}"] \arrow[rd, "\text{quantum feature map}", blue] & [8em] x_{n+1} \\
                \stext{Quantum} & \rho_n \arrow[r, "\text{Koopman evolution}", swap] \arrow[ru, "\text{flux terms}", red]& [8em] \tilde\rho_{n+1} \arrow[r, "\text{quantum state update}", swap] & \rho_{n+1}   
            \end{tikzcd}
        };
    \end{tikzpicture}
    \caption{\label{figSchematic}Schematic representation of a QMCl cycle. The diagram illustrates the update procedure for the classical degrees of freedom (resolved variables) $x_n$ and the quantum state (density operator) $\rho_n$ over one timestep. We advance the resolved variables to their value $x_{n+1}$ at time $t_{n+1} = t_n +\Delta t$ using the values $x_n$ at time $t_n$ and the flux terms obtained by evaluation of a quantum mechanical observable on the state $\rho_n$ (red arrow). Moreover, we advance the quantum state $\rho_n$ to a prior state $\tilde\rho_{n+1}$ at time $t_{n+1}$ through a unitary evolution step induced by the Koopman operator. The state $\tilde\rho_{n+1}$ is then updated to the posterior $\rho_{n+1}$ via a projective update (quantum Bayes' rule) induced by a quantum feature map (blue arrow). The cycle is repeated to advance the coupled classical--quantum states $(x_{n+1},\rho_{n+1})$ at later times.}
\end{figure}
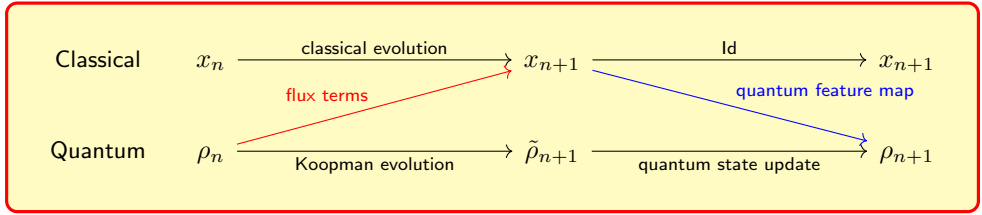

Building on the QMDA framework, we propose a data-driven technique for closure of dynamical systems that employs a quantum mechanical system as a surrogate model of the unresolved degrees of freedom. The surrogate quantum system evolves in tandem with the classical system for the resolved variables via two-way coupling, generalizing the one-way coupling in QMDA. A schematic overview of our approach, which we call Quantum Mechanical Closure (QMCl), is depicted in Fig.~\ref{figSchematic}. 

At any time $t_n$, the state of the quantum mechanical system is represented by a density operator $\rho_n$ acting on a Hilbert space $H_L$ of (finite) dimension $L$. The closure terms needed to advance the resolved variables $x_n$ at time $t_n$ to their state $x_{n+1}$ at time $t_{n+1} = t_n+ \Delta t$ are obtained by evaluation of a quantum mechanical observable $ A$ (i.e., a self-adjoint operator acting on $H_L$) on the state $\rho_n$. Intuitively, the density operator $\rho_n$ can be thought of as a non-commutative analog of a classical probability density, and the evaluation of $ A$ on $ \rho_n$ is analogous to the evaluation of an expectation functional in classical probability theory. Thus, the update $x_n \mapsto x_{n+1}$ given $ \rho_n$ can be thought of as ensemble update based on the statistical information about the unresolved degrees of freedom encoded in $\rho_n$. 

To advance $ \rho_n$ to the quantum state $ \rho_{n+1}$ at time $t_{n+1}$, we first employ a unitary evolution map, $ \rho_{n+1} \mapsto \tilde \rho_{n+1}$, induced by the Koopman operator of the dynamical system \cite{EisnerEtAl15}. Then, we update $\tilde \rho_{n+1} \mapsto \rho_{n+1}$ projectively given the resolved variables $x_{n+1}$ using an operator-valued function that we call a quantum feature map. These steps will be described precisely below, but for now we note that the unitary Koopman evolution can be thought of as being analogous to the forecast step in data assimilation \cite{MajdaHarlim12,LawEtAl15}. In this analogy, the Koopman evolution leads to a ``prior'' state $\tilde\rho_{n+1}$, which is updated projectively to the ``posterior'' $\rho_{n+1}$ given $x_{n+1}$ analogously to the Bayesian analysis step in data assimilation. In quantum information theory, the values of the feature map used here to assimilate the resolved degrees of freedom are known as \emph{quantum effects} \cite{Gudder07}. Quantum effects generalize the notion of events in classical statistics. 

The principal distinguishing aspects of QMCl are as follows.
\begin{enumerate}
    \item We cast the problem of closure of dynamical systems into the framework of operator theory and quantum information theory. Non-abelian quantum systems are known to exhibit new types of behavior (such as entanglement and superposition), and offer greater capacity to encode information than classical systems \cite{GohEtAl18}. In particular, the dimension of the space of density operators on an $L$-dimensional Hilbert space is $O(L^2)$, whereas the dimension of the space of $L$-dimensional probability density vectors is $O(L)$. Our main premise is that the mathematical framework of quantum mechanics provides a flexible arena for building closure schemes for classical dynamical systems enjoying favorable structure preservation and asymptotic consistency properties. 
    \item Positivity preservation: Our procedure for constructing the coupled classical--quantum system depicted in Fig.~\ref{figSchematic} employs, as an intermediate step, an information-preserving embedding of the full dynamical system (i.e., both the resolved and unresolved degrees of freedom) into an infinite-dimensional quantum system. This system is projected to the finite-dimensional system on the Hilbert space $H_L$. The process of embedding the classical system into the infinite-dimensional quantum system \emph{before} projecting into finite dimensions allows us to take advantage of properties of non-abelian operator spaces to ensure that the parameterization scheme is positivity preserving. That is, positive-definite flux terms from the unresolved variables to the resolved variables in the original system are represented by positive-definite operators in the QMCl-parameterized system (and similarly for negative-definite terms). This addresses an important problem in data-driven subgrid-scale modeling, where failure to preserve the sign of sign-definite quantities (e.g., moisture in an atmospheric circulation model) is recognized as a significant source of bias and numerical instability \cite{YuvalOGorman20}. 
    \item Data-driven formulation: QMCl employs a closely-related data-driven formulation to QMDA, which uses eigenfunctions of kernel integral operators as data-driven basis functions for representation of linear operators \cite{BerryEtAl15}. The training data requirements of QMCl are comparable to other data-driven closure schemes, and typically comprise of time series of the resolved variables and the unresolved fluxes. In what follows, we explore scenarios where training data for building the basis are either directly available from integrations of the full model (as would be the case, e.g., when coarsening a high-resolution climate model \cite{BoltonZanna19,BrenowitzBretherton18,YuvalOGorman20}), or they are constructed from the resolved variables (as would be the case, e.g., when learning residuals of imperfect models from observations of nature). In the latter scenario, we use delay-coordinate embedding techniques \cite{SauerEtAl91} to improve the the richness of our approximation models. Using techniques for pointwise and spectral approximation of linear operators \cite{Chatelin11, VonLuxburgEtAl08}, the data-driven QMCl models converge in the large-data limit to data-independent finite-dimensional models on the Hilbert space $H_L$. Increasing the dimension parameter $L$ allows the finite-dimensional models to capture increasingly high-dimensional subgrid-scale dynamics. 
\end{enumerate}

In this paper, we describe the mathematical framework and associated computational algorithms of QMCl. Moreover, we demonstrate the behavior of the scheme in two test systems: the Lorenz 63 (L63) \cite{Lorenz63} and Lorenz 96 (L96) multiscale system \cite{Lorenz96,FatkullinVandenEijnden04}. The L63 system was employed in \cite{Palmer01} as a motivating example for stochastic parameterization. The L96 multiscale system has been employed as a prototype test system in several studies on  parameterization and surrogate modeling, e.g., \cite{CrommelinVanden-Eijnden08,JiangHarlim20b,BurovEtAl21,SchneiderEtAl21}.     

\subsection{Plan of the paper} 

In \cref{secStatement}, we state the parameterization problem under study. \Cref{secQuantumMechanics} gives an overview of the axioms of quantum mechanics. In \cref{secQMDA}, we describe the construction of the quantum mechanical system that we will use to model the unresolved degrees of freedom. This construction largely parallels the QMDA approach from \cite{FreemanEtAl22}. In \cref{secMethodology}, we describe the QMCl parameterization scheme based on the quantum system from \cref{secQMDA}. This is followed by a description of the data-driven implementation of our approach in  \cref{secDataDriven}. In \cref{secL63,secL96}, we present applications of QMCl to the L63 and L96 multiscale systems, respectively. \Cref{secDiscussion} contains a discussion and perspectives on future work. Details on numerical implementation and pseudocode are collected in \cref{appNumerics}.

\section{\label{secStatement}Problem Statement}

Consider a discrete-time dynamical system with state space $\Omega$ and evolution map $\Phi : \Omega \to \Omega$ preserving a probability measure $\mu$. Suppose that the state space is decomposable as the product $ \Omega = \mathcal X \times \mathcal Y$, where $\mathcal X$ and $\mathcal Y$ are state spaces for the resolved and unresolved degrees of freedom, respectively. Throughout the paper, we will assume that $\Omega$ is a ``nice'' (separable, completely metrizable) topological space, and $\mu$ is a Borel measure with compact support. Additional assumptions will be introduced as needed.     

We view the $\Phi$ dynamics on $\Omega$ as being as being inaccessible to direct simulation, e.g., due to the dimension of $\mathcal Y$ being prohibitively large, lack of complete knowledge of the evolution map $ \Phi$, or a combination thereof. The general problem of parameterization is to build a surrogate dynamical system $\tilde \Phi : \tilde \Omega \to \tilde \Omega$ on a state space $ \tilde \Omega = \mathcal X \times \tilde{\mathcal Y} $ which is (i) feasible to integrate; and (ii) compatible with the evolution of the resolved degrees of freedom in $\mathcal X$. Here, $\tilde{\mathcal Y}$ is the state space of a surrogate model of the unresolved degrees of freedom in $\mathcal Y$. 

Let $P_{\mathcal X} : \Omega \to \mathcal X$ and $\tilde P_{\mathcal X} : \tilde \Omega \to \mathcal X $ denote the canonical projection maps onto $\mathcal X$ from $\Omega$ and $\tilde \Omega$, respectively. Let also $\alpha : \Omega \to \tilde \Omega$ be a map that assigns a dynamical state in $\tilde \Omega$ to each dynamical state in $ \Omega$. A parameterized system on $\tilde \Omega$ which is consistent with the original system on $\Omega$ would satisfy the following commutative diagram:
\begin{equation}
    \label{eqComm}
    \begin{tikzcd}
        \Omega \ar[rrr,"\Phi"] \ar[dr,"P_{\mathcal X}", swap] \ar[dd,"\alpha", swap] & & & \Omega \ar[dl,"P_{\mathcal X}"] \ar[dd,"\alpha"] \\
        &X & X\\
        \tilde \Omega \ar[rrr,"\tilde \Phi"]  \ar[ur,"\tilde P_{\mathcal X}"] & & & \tilde \Omega \ar[ul,"\tilde P_{\mathcal X}", swap]  
    \end{tikzcd}
\end{equation}
This diagram represents the fact that for initial conditions $\tilde \omega = \alpha(\omega) \in \tilde \Omega$ that are consistent with $\alpha \in \Omega$ in terms of the resolved variables, i.e., $\tilde P_{\mathcal X}(\tilde \omega) = P_{\mathcal X}(\omega)$, the $\tilde \Phi$ dynamics should produce consistent outcomes with $\Phi$, (again, in terms of the resolved variables), i.e., $\tilde P_{\mathcal X} \circ \tilde \Phi(\tilde \omega) = P_{\mathcal X} \circ \Phi(\omega)$. In practice, such a strong form of consistency is seldom achievable, and one seeks to identify an appropriate state space $\tilde{\mathcal Y}$ and dynamics $\tilde \Phi$ such that~\eqref{eqComm} holds as closely as possible in some sense, e.g., in terms of low-order statistics such as probability distributions and time-autocorrelation functions of the resolved variables with respect to the invariant measure.       

In general, it is advantageous to include any partial knowledge about the $\Phi$ dynamics in the construction of $\tilde \Phi$. To incorporate that possibility in our formalism, we introduce (without loss of generality) a space $\mathcal Z$ and maps $ \phi: \mathcal X \times \mathcal  Z \to \mathcal X$ and $ Z : \mathcal Y \to \mathcal Z$ such that   
\begin{displaymath}
    P_{\mathcal X} \circ \Phi(\omega) = \phi(x, Z(y)), \quad \omega = (x,y). 
\end{displaymath}
Here, $\phi$ captures the functional form  of the evolution map for the resolved variables conditioned on a function $Z$ (a ``flux term'') of the unresolved variables. Assuming that $ \phi $ is known, we can construct $\tilde \Phi$ as a map of the form   
\begin{equation}
    \label{eqPhiDecomp}
    \tilde \Phi( x, \tilde y) = (\phi(x, \tilde Z(\tilde y ), \tilde \psi(x, \tilde y)),
\end{equation}
where $\tilde \psi: \mathcal X \times \tilde {\mathcal Y} \to \tilde{\mathcal Y}$ represents the evolution of the surrogate unresolved variables $\tilde y$ given the resolved variables $x$, and $\tilde Z : \tilde{\mathcal Y} \to \mathcal Z$ is a function acting as a surrogate for $Z$. If $\phi$ is not known, we replace~\eqref{eqPhiDecomp} by  
\begin{equation}
    \label{eqPhiDecomp2}
    \tilde \Phi( x, \tilde y) = (\tilde \phi(x, \tilde Z(\tilde y ), \tilde \psi(x, \tilde y)),
\end{equation}
where $\tilde \phi : \mathcal X \times \mathcal Z \to \mathcal X$ is an approximation of $\phi$.

\subsection{Examples}
We outline how some of the commonly used parameterization schemes fit within the framework described above.    

\subsubsection{Deterministic functional closure}
In closure schemes modeling the unresolved variables through functional relationships with the resolved variables, we formally set $\tilde{\mathcal Y} = \mathcal X$, so that the flux term $\tilde Z : \mathcal X \to \mathcal Z$ becomes a function of the resolved variables. The map $\tilde\psi(x,\tilde y)$ is also formally set to a copy of $\phi(x,\tilde Z(\tilde y))$ so it need not be explicitly retained in calculations. A number of the classical parameterization schemes (e.g., \cite{ArakawaSchubert74,GentMcWilliams89} in geophysical fluid dynamics) derive $\tilde Z$ heuristically using prior physical knowledge of the problem at hand. For certain classes of systems exhibiting timescale separation between the $x$ and $y$ variables, $\tilde Z$ can be systematically obtained from the governing equations of the original system using asymptotic analysis techniques such as averaging \cite{PavliotisStuart08}. In data-driven approaches, the function $\tilde Z$ is constructed via a supervised learning algorithm using $(x,Z(y))$ pairs as training data; e.g., \cite{BrenowitzBretherton18,BoltonZanna19,YuvalOGorman20}.

\subsubsection{\label{secStochasticExamples}Stochastic closure}
To treat stochastic parameterization schemes, we set $\tilde{\mathcal Y}$ to a probability sample space, $\tilde Z : \tilde{\mathcal Y} \to \mathcal Z$ to a random variable, and choose $\tilde \psi$ such that $\{ \tilde \Phi^n : n \in \mathbb N\}$ defines a stochastic process over $\tilde \Omega$. As a concrete example, we may understand Palmer's Gaussian-distributed stochastic closure of the L63 system \cite{Palmer01} in this framework, as follows. 

In a coordinate system obtained via empirical orthogonal function (EOF) analysis, the L63 governing equations with the standard parameter values become
\begin{equation}
    \label{eqL63Palmer}
    \begin{aligned}
        \dot{a_1} & = v_1(a_1,a_2,a_3) := 2.3\,a_1 - 6.2\,a_3 - 0.49\, a_1 a_2 - 0.57\, a_2 a_3, \\ 
        \dot{a_2} & = v_2(a_1,a_2,a_3) := -62 - 2.7\,a_2 + 0.49\, a_1^2 - 0.49\, a_3^2 + 0.14\, a_1 a_3, \\
        \dot{a_3} & = v_3(a_1,a_2,a_3) :=  -0.63\,a_1 - 13\, a_3 + 0.43\, a_1 a_2 + 0.49\, a_2 a_3,
    \end{aligned}
\end{equation}
where the coordinates $a_1$, $a_2$, and $a_3$ have decreasing variance. Let $\Omega= \mathbb R^3$ and $\Phi^t : \Omega \to \Omega$ with $ t \in \mathbb R$ be the flow generated by~\eqref{eqL63Palmer}. In this example, we consider the discrete-time system $\Phi : \Omega \to \Omega $ with $ \Phi = \Phi^{\Delta t}$ induced by the continuous-time flow at a fixed timestep $\Delta t$ as the ``true'' dynamics. Moreover, we let $\Omega = \mathcal X \times \mathcal Y$ with $ x= (a_1,a_2) \in \mathcal X \equiv \mathbb R^2 $ and $y = a_3 \in \mathcal Y \equiv \mathbb R$ represent the resolved and unresolved variables, respectively. The corresponding flux term is $Z : \mathcal Y \to \mathcal Z$ with $ \mathcal Z = \mathbb R  $ and $Z(y) = y$. As approximate resolved dynamics, we employ the map $\tilde \phi : \mathcal X \times \mathcal Z \to \mathcal X$ obtained by a forward Euler discretization of~\eqref{eqL63Palmer} at the timestep $\Delta t$; that is, we have    
\begin{displaymath}
    \tilde\phi(x,z) = x + (v_1(x,z), v_2(x,z)) \, \Delta t.
\end{displaymath}

In Palmer's approach, the evolution of the $a_3$ variable is approximated by an i.i.d.\ Gaussian process. To represent this process within the framework described above, we let $\tilde{\mathcal Y} = \mathbb R^{\mathbb N}$ be a space of $\mathbb R$-valued sequences equipped with a $\sigma$-algebra $\Sigma_{\tilde{\mathcal Y}}$ of cylinder sets and a shift-invariant Gaussian measure $\nu : \Sigma_{\tilde{\mathcal Y}} \to [0,1]$. That is, $\nu$ has the properties that (i) $\nu \circ T^{-1} = \nu $, where $T : \tilde{\mathcal Y} \to \tilde{\mathcal Y}$ is the shift map, $T(\tilde y_0,\tilde y_1, \ldots, ) = (\tilde y_1,\tilde y_2,\ldots)$; and (ii) $\tilde Z_0, \tilde Z_1,\ldots$ with $\tilde Z_n = u \circ T^n $ are i.i.d.\ Gaussian random variables with equal variance. Here, $ u: \tilde{\mathcal Y} \to \mathcal Z $ projects onto the first component, $u(\tilde y_0, \tilde y_1, \ldots) = \tilde y_0$, and the measure $\nu$ is chosen such that $Z_n$ has variance equal to the variance of $a_3$. With these definitions, the parameterized dynamics from~\eqref{eqPhiDecomp2} becomes $\tilde \Phi(x,\tilde y) = (\tilde \phi(x, \tilde Z(\tilde y), \tilde \psi(x,\tilde y)))$,
with $\tilde \psi(x,\tilde y) = u(T(\tilde y))$.   

Palmer observed that with this approach the parameterized $x$ dynamics recovers the qualitative features of the $(a_1,a_2) $ variables of the full L63 system, including the characteristic lobes of the Lorenz attractor projected onto the $(a_1,a_2)$ plane; see \cite[Figure~2]{Palmer01} and \cref{figPalmer2} below. In particular, since the state space $\tilde \Omega$ of the stochastically parameterized system is formally infinite-dimensional, the method is able to overcome theoretical limitations of a deterministic functional closure on $\mathcal X = \mathbb R^2$, which by the Poincar\'e-Bendixson theorem cannot exhibit chaos in continuous time under smooth dynamics.        

In \cref{secL63}, we will use Palmer's L63 closure scheme as a reference example to assess the performance of QMCl. Certain classes of SDE schemes derived via homogenization principles \cite{PavliotisStuart08} can similarly be understood as a state-space augmentation of the resolved variables space $\mathcal X$ to include the sample space of a Brownian motion driving the parameterized system.   

\section{\label{secQuantumMechanics}Overview of quantum mechanics}

As a statistical theory, quantum mechanics can be looked at as a generalization of Bayesian probability theory \cite{Holevo01,Gudder07}. Quantum mechanics deals directly with objects (e.g., quantum particles) for which a measurement of any observable at a given time is inherently probabilistic in nature. The theory of quantum mechanics deals with understanding how the generalized probability distributions associated with these objects (which, in some views, constitute the objects themselves) evolve over time and under experimental observations of the objects' properties. The probability theory associated with quantum systems, however, is generalized beyond classical probability theory. Quantum probabilistic dynamics are understood through non-abelian algebras of operators on Hilbert spaces and spectral theory applied to these operators. The QMDA methodology arises from a parallel between the probabilistic dynamics formulated in this way and the mathematical structure of dynamical systems and Koopman operator theory \cite{Giannakis19b,FreemanEtAl22}.  

A standard framework of quantum probabilistic dynamics arises from the Dirac-von Neumann quantum mechanical axioms. In this section, we provide an overview of the mathematical objects and evolution rules associated with this theory in order to get a general view of how the mathematical systematization of quantum mechanics can be applied to dynamical systems more broadly. More detailed explication of the foundations of this theory can be found, e.g., in \cite{Takhtajan08}. 

\subsection{Hilbert spaces and quantum observables} 

In the Dirac-von Neumann framework, associated to every quantum system is a Hilbert space $\left( \mathcal H, \langle \cdot, \cdot \rangle \right)$ over the complex numbers. The precise choice of $\mathcal H$ will depend on the specific system being studied. Quantum observables are defined as elements of the set of self-adjoint (possibly unbounded) linear operators $A: D(A) \to \mathcal H$, where the domain $D(A)$ of $A$ is a dense subspace of $\mathcal H$.

Every quantum observable is associated with a measurable real-valued quantity of the quantum system (such as momentum, spin, etc.). Importantly, if the quantity associated with a quantum observable $A$ is measured in an experimental setting, the possible values of the measurement are contained in the spectrum of $A$. The self-adjointness condition ensures that such measurements are real-valued. 

Let $B(\mathcal H)$ be the space of bounded operators on $\mathcal H$. From the spectral theorem for self-adjoint operators, associated with every quantum observable $A : D(A) \to \mathcal H$ is an operator-valued  measure $E_A: \mathcal{B}(\mathbb{R}) \to B(\mathcal H)$ on the Borel $\sigma$-algebra $\mathcal B(\mathbb R)$ over $\mathbb R$ such that (i) for every Borel set $S \in \mathcal B(\mathbb R)$, $E_A(S)$ is an orthogonal projection; (ii) $E(\emptyset) = 0$; (iii) $E_A(\mathbb R) = \Id $; and (iv) for every countable collection $\{S_i \}$ of pairwise-disjoint sets $S_i \in \mathcal B(\mathbb R)$, we have $E_A(\bigcup_i S_i) = \sum_i E(S_i)$, where the sum over $i$ converges in the strong topology of $B(\mathcal H)$. Using $E_A$, one can reconstruct the operator $A$ through the spectral integral $ A = \int_\mathbb{R} a \,dE_A(a)$.

The measure $E_A$ plays an important role in calculations, and it can be resolved as a finite discrete sum of projections for QMCl applications (see \cref{secDiscretization}). In the general theory, however, the spectrum of $A$ has a continuous component and $A$ cannot be finitely resolved in this way. Unlike classical observables lying in commutative algebras of functions, quantum observables do not, in general, commute. One of the signature characteristics of quantum mechanics arises from this fact---non-commutative observables cannot be measured simultaneously at arbitrarily high precision, giving rise to the famed uncertainty principle.

\subsection{\label{secQuantumStates}Quantum states} 
Let $B_1(\mathcal H) \subseteq B(\mathcal H)$ be the space of trace class operators on $\mathcal H$. The quantum states of the system are defined as elements of the set $Q(\mathcal H) \subset B_1(\mathcal H)$ consisting of all positive, trace class operators of unit trace, i.e.,
\begin{displaymath}
    Q(\mathcal H) = \{\rho \in B_1(\mathcal H) \mid\text{$\tr\rho=1$, $\rho \geq 0$}\}.
\end{displaymath}
Such operators $\rho \in Q(\mathcal H)$ are known as density operators. Quantum states can be thought of as roughly analogous to classical probability densities in an $L^1$ space, representing the degree of knowledge about the system's current state. With this in mind, a density operator's positivity and unity of trace can be thought of as analogous to the positivity and normalization of a classical probability density. We say that $\rho \in Q(\mathcal H)$ is a pure state if there exists a unit vector $\xi \in \mathcal H$ such that $ \rho = \langle \xi, \cdot \rangle \xi$, i.e., $\rho$ is a rank-1 projection along $\xi$. Otherwise, $\rho$ is said to be mixed.  

Every quantum state $\rho \in Q(\mathcal H)$ induces a continuous, positive linear functional $\mathbb E_\rho : B(\mathcal H) \to \mathbb C$ on bounded operators defined as
\begin{equation}
    \label{eqQExp}
    \mathbb E_\rho A = \tr(\rho A);
\end{equation}
intuitively, we think of this formula as being analogous to an expectation with respect to a classical probability density. If $A$ is a self-adjoint quantum observable (not necessarily bounded), the associated projection-valued measure $E_A$ induces a probability distribution for experimental measurements of $A$. Specifically, the probability that a measurement $a \in \mathbb{R}$ of the quantum observable $A$ lies in a Borel set $S \in \mathcal B(\mathbb R)$ is given by 
\begin{equation}
    \label{eqQProb}
    \mathbb P_{\rho,A}\left( a \in S \right) = \mathbb E_\rho E_A(S).
\end{equation}
This general formula holds for all quantum observables, and thus the quantum state $\rho$ alone determines the probability distribution of all quantum observables of a system. 

\subsection{Quantum evolution and measurement}

In the Dirac-von Neumann axioms, there are two evolution rules for the quantum state. The first is the time-evolution equation. Letting $\mathcal T$ denote either $\mathbb Z$ or $\mathbb R$, respectively for discrete- or continuous-time systems, suppose that $\rho_0 \in Q(\mathcal H)$ is the initial quantum state of a system, and $\rho_t$ the quantum state of the system at a time $ t \in \mathcal T$ such that no intervening measurement took place. The quantum mechanical time evolution axiom posits that there exists a group of unitary operators $\{ U^t \in B(\mathcal H) \}_{t\in\mathcal T}$ such that
\begin{equation} 
    \label{eqSchrodinger}
    \rho_t = \mathcal P^t \rho_0 := U^{t*} \rho_0 U^t. 
\end{equation}
The collection $ \{ \mathcal P^t\}_{t \in \mathcal T}$ extends to an evolution group on $B_1(\mathcal H)$, which can be thought of as an analog of the group of transfer operators of a measure-preserving dynamical system \cite{Baladi00} acting on densities in $L^1$. If the quantum system under study is open, i.e., it is allowed to interact with its environment, \eqref{eqSchrodinger} is replaced by evolution under a semigroup $\{ \mathcal P^t : B_1(\mathcal H) \to B_1(\mathcal H) \}_{t \in \mathcal T_+}$ which is required to be trace non-increasing, i.e., we have  $\tr(\mathcal P ^t \rho)  \leq 1$ rather than the stronger property $ \tr (\mathcal P^t \rho) = 1$ under~\eqref{eqSchrodinger}. In either case, the evolution of quantum states under $\mathcal P^t$ is dual to the evolution of quantum observables given by the operators $\mathcal U^t : B(\mathcal H) \to B(\mathcal H)$ with 
\begin{equation}
    \label{eqHeisenberg}
    \mathcal U^t A = U^t A U^{t*};
\end{equation}
i.e., we have $ \mathbb E_{\mathcal P^t \rho_0} A = \mathbb E_{\rho_0} (\mathcal U^t A)$. In quantum mechanics, the evolution of states and observables in \eqref{eqSchrodinger} and~\eqref{eqHeisenberg} is known as the Schr\"odinger and Heisenberg picture, respectively.

The second evolution rule addresses how measurements of quantum observables affect the quantum state.  Suppose a quantum system is in the state $\rho \in Q(\mathcal H)$ at a time $t$, and a quantum observable $A$ is measured at the time $t$. Then, if the measurement takes a value $a \in \mathbb R$ such that $E_A(\{a\}) \neq 0$ (which necessarily means that $a$ is an eigenvalue of $A$), the quantum state of the system is now given by $\rho|_a \in Q(\mathcal H)$, where
\begin{equation}
    \label{eqQBayes}
    \rho|_a = \frac{E_A(\{a\}) \rho E_A(\{a\})}{\tr(E_A(\{a\}) \rho E_A(\{a\}))}. 
\end{equation}
The update $\rho \mapsto \rho|_a$ is oftentimes referred to as a von Neumann measurement, and can be thought of as an analog of the Bayesian conditioning rule of classical probability theory. 

In~\eqref{eqQBayes}, the projections $E_A(\{a\})$ are examples of quantum effects, defined as the operators in the set
\begin{equation}
    \label{eqEffect}
    \mathcal E(\mathcal H) = \{ e \in B(\mathcal H) \mid 0 \leq e \leq \Id \}.
\end{equation}
Intuitively, one can think of effects as analogous to events in classical probability theory. The projections $e \in \mathcal E(\mathcal H)$ represent events which are ``sharp'', in the sense that $e^2 = e$, whereas other effects such that $e^2 < e$ can be thought of as being ``fuzzy''. 

A generalization of~\eqref{eqQBayes} that represents conditioning of $\rho$ by an arbitrary effect $ e \in \mathcal E(\mathcal H)$ is
\begin{equation}
    \label{eqQBayes2}
    \rho|_e = \frac{\sqrt e \rho \sqrt e}{\tr(\sqrt e \rho \sqrt e)}. 
\end{equation}
This update rule reduces to~\eqref{eqQBayes} when $e = e^2 = E_A(\{a\})$ is a projection. In applications involving physical quantum systems, \eqref{eqQBayes2} is used to model the state change of a quantum system being measured due to interactions with the measuring apparatus, a process sometimes called ``state collapse''. Equation~\eqref{eqQBayes2}, however, also applies more broadly as an inference rule within abstract quantum probability theory \cite{Schack02}, which is how we will use it in our quantum mechanical closure schemes.          

\section{\label{secQMDA}Quantum mechanical representation of classical dynamics}

In this section, we construct the quantum mechanical system that we will use to model the unresolved degrees of freedom in QMCl. Our goal is to build a quantum system on a finite-dimensional Hilbert space, $H_L$, and use the quantum state space of this system as the state space of the surrogate dynamical model for the unresolved degrees of freedom, i.e., $\tilde{\mathcal Y} = Q(H_L)$ in the notation of \cref{secStatement}. Following the QMDA approach \cite{FreemanEtAl22}, the first step in this construction is a formal embedding of the classical dynamics on $\Omega$ into a quantum system associated with an infinite-dimensional Hilbert space, $H$. Under this embedding, classical observables (functions of the state) are represented by quantum mechanical observables, and classical probability densities are represented by quantum density operators. We describe this embedding in \cref{secInfDim}, following a brief overview of relevant definitions from Koopman and transfer operator theory. In \cref{secDiscretization}, we construct finite-dimensional quantum systems by discretization (i.e., finite-dimensional projection) of the infinite-dimensional system, choosing $H_L \subset H$ as an $L$-dimensional space in a nested family of finite-dimensional subspaces associated with an orthonormal basis of $H$. In \cref{secHilb}, we make explicit our choice of basis leading to the construction of $H_L$.  

\subsection{Koopman and transfer operators}

The measure-preserving dynamical system $ \Phi: \Omega \to \Omega $ induces groups of linear operators $ U^n : L^p(\mu) \to L^p(\mu) $ on the $L^p$ spaces associated with the invariant measure, acting by composition with the dynamical flow, $U^n f = f \circ \Phi^n $ with $ n \in \mathbb Z$. These operators are known as Koopman operators \cite{Koopman31,EisnerEtAl15}, and preserve the $L^p(\mu) $ norm for any $ p \in [0, \infty]$ since $\mu$ is an invariant measure under $\Phi$. In the Hilbert space case, $ H := L^2(\mu)$, the Koopman operators are unitary, i.e., $(U^n)^* = U^{-n}$ for any $n \in \mathbb Z$. We define the inner product of $H$ as $ \langle f, g \rangle := \int_\Omega f^* g \, d\mu$, and use the notations $\lVert f \rVert_{L^p(\mu)} := (\int_m \lvert f \rvert^p \, d\mu)$ for $ p \in [1,\infty)$ and $ \lVert f \rVert_{L^\infty(\mu)} = \lim_{p\to\infty} \lVert f \rVert_{L^p(\mu)}$ for the standard $L^p(\mu)$ norms.

Identifying $L^p(\mu)$, $p \in (1,\infty] $, with the dual of $L^q(\mu)$, $ q \in [1,\infty)$, with $\frac{1}{p} + \frac{1}{q} = 1$, we have that $U^n : L^p(\mu) \to L^p(\mu) $ is the adjoint (dual) of the operator $P^n : L^q(\mu) \to L^q(\mu)$ with $P^n f = f \circ \Phi^{-n}$, known as the transfer operator \cite{Baladi00}. If we further identify $L^q(\mu)$ as the space of finite Borel measures on $\Omega$ with densities in $L^q(\mu)$, denoted here as $M_q(\mu)$, the transfer operator can be identified with the pushforward map $\Phi^n_* : M_q(\mu) \to M_q(\mu)$, where $\Phi^n_* \nu = \nu \circ \Phi^{-n}$. That is, if $\nu \in M_q(\mu) $ has density $ \varrho = \frac{d\nu}{ d\mu} \in L^q(\mu)$, then the measure $\Phi^n_* \nu$ has density $ \frac{d \Phi^n_* \nu}{d\mu} = P^n \varrho$. 

Koopman and transfer operators are central analytical tools in ergodic theory \cite{Baladi00,EisnerEtAl15}. In recent years, they have found widespread use in data-driven analysis and forecasting methodologies for dynamical systems; e.g., \cite{DellnitzEtAl00,Mezic05,RowleyEtAl09,WilliamsEtAl15,BruntonEtAl17,KlusEtAl20,BerryEtAl20}. We will be addressing aspects of data-driven approximation in \cref{secDataDriven}, but for the remainder of this section our focus will be on quantum mechanical systems induced by the unitary Koopman operators on $H$ and its subspaces.

\subsection{\label{secInfDim}Embedding into an infinite-dimensional quantum system}

The unitary Koopman evolution group $ \{ U^n \in B(H) \}_{n\in\mathbb Z}$ induces a quantum system on $H$ with states $Q(H) \subset B_1(H)$ as described in \cref{secQuantumMechanics}. Moreover, we can identify the unitary evolution of states in~\eqref{eqSchrodinger}  with an induced action of the Koopman operator; that is,
\begin{displaymath}
    \rho_n = \mathcal P^n \rho_0 := U^{n*} \rho_0 U^n,
\end{displaymath}
where $\rho_0 \in Q(H)$ is the initial state and we set $t \equiv n \in \mathbb Z$ since we work in discrete time. This quantum system can be brought into correspondence with the underlying classical dynamical system, as follows.

\subsubsection{\label{secVonNeumann}Von Neumann algebras}Consider the space of classical observables $L^\infty(\mu)$. In addition to being a Banach space as the other $L^p(\mu)$ spaces, $L^\infty(\mu)$ has the distinguished property of being a von Neumann algebra \cite{Takesaki01} with respect to pointwise function multiplication and complex conjugation. This means:
\begin{enumerate}
    \item For any $f, g \in L^\infty(\mu)$, the pointwise product $fg$ lies in $L^\infty(\mu)$, and we have
        \begin{displaymath}
            \lVert f g \rVert_{L^\infty(\mu)} \leq \lVert f \rVert_{L^\infty(\mu)} \lVert g \rVert_{L^\infty(\mu)}, \quad \lVert f^* f \rVert_{L^\infty(\mu)} = \lVert f \rVert_{L^\infty(\mu)}^2,
        \end{displaymath}
        making $L^\infty(\mu)$ a $C^*$-algebra.
    \item $L^\infty(\mu)$ is the continuous dual of the Banach space $L^1(\mu)$, making it a von Neumann algebra.
\end{enumerate}
By commutativity of function multiplication, $L^\infty(\mu)$ is an abelian algebra. Given a probability density $p \in L^1(\mu)$, i.e., a positive function satisfying the normalization condition $ \int_\Omega p \, d\mu = 1$, we have a continuous, positive linear functional $\mathbb E_p : L^\infty(\mu) \to \mathbb C$, defined as the expectation
\begin{displaymath}
    \mathbb E_p f = \int_\Omega p f \,d\mu.
\end{displaymath}
Letting $ P(\mu) \subset L^1(\mu)$ be the space of probability densities in $L^1(\mu)$, i.e., 
\begin{displaymath}
    P(\mu) = \left\{ p \in L^1(\mu) \mid \text{$p \geq 0$, $\textstyle\int_\Omega p \, d\mu = 1$} \right\},
\end{displaymath}
we can think of elements of $P(\mu)$ as states of the algebra $L^\infty(\mu)$.

The similarities between this construction and the construction of quantum states from \cref{secQuantumStates} is not accidental. Given a Hilbert space $\mathcal H$, the space $B(\mathcal H)$ is a Banach space equipped with the operator norm, $\lVert A \rVert_{B(\mathcal H)} = \sup_{u \in \mathcal H} \frac{\lVert A u \rVert_{\mathcal H}}{\lVert u \rVert_{\mathcal H}}$, and it is a von Neumann algebra with respect to operator composition and adjunction. This means:
\begin{enumerate}
    \item For any $A, B \in B(\mathcal H)$, the operator product $AB$ lies in $B(\mathcal H)$, and we have
        \begin{displaymath}
            \lVert A B \rVert_{B(\mathcal H)} \leq \lVert A \rVert_{B(\mathcal H)} \lVert B \rVert_{B(\mathcal H)}, \quad \lVert A^* A \rVert_{B(\mathcal H)} = \lVert A \rVert_{B(\mathcal H)}^2,
        \end{displaymath}
        making $B(\mathcal H)$ a $C^*$-algebra.
    \item $B(\mathcal H)$ is the continuous dual of the Banach space $B_1(\mathcal H)$, equipped with the trace norm $\lVert A \rVert_{B_1(\mathcal H)} = \tr\sqrt{A^* A}$, making it a von Neumann algebra.
\end{enumerate}
Aside from trivial cases, the algebra $B(\mathcal H)$ is non-abelian, which is a fundamental reason for differences in behavior between classical and quantum systems.

\subsubsection{\label{secInfDimObs}Embedding observables}Setting $\mathcal H = H \equiv L^2(\mu)$, the abelian algebra $L^\infty(\mu)$ embeds naturally and isometrically into $B(H)$ through its regular representation, i.e., the map $\pi : L^\infty(\mu) \to B(H)$ that maps elements of $L^\infty(\mu)$ to multiplication operators in $B(H)$, i.e., 
\begin{displaymath}
    (\pi f) g = f g, \quad \forall g \in H.
\end{displaymath}
This map is a $^*$-homomorphism of $L^\infty(\mu)$, i.e.,
\begin{displaymath}
    \pi 1_\Omega = \Id, \quad \pi (fg) = (\pi f)(\pi g), \quad \pi( f^* ) = (\pi f)^*,
\end{displaymath}
where $f,g$ are arbitrary elements of $L^\infty(\mu)$ and $ 1_\Omega \in L^\infty(\mu)$ is the element equal to 1 $\mu$-a.e. Note that $\pi f$ is self-adjoint if and only if $f$ is real-valued (i.e., $f$ is a self-adjoint element of the $L^\infty(\mu)$ algebra), and $\pi f $ is a positive operator if and only if $ f \geq 0 $ $\mu$-a.e. Moreover, the spectrum of the multiplication operator $\pi f$ is equal to the essential range of $f$, which corresponds to the values of $f$ that have nonzero probability of occurring with respect to $\mu$. 

One can also verify that $\pi$ is dynamically compatible with the evolution of classical observables under the Koopman operators $U^n : L^\infty(\mu) \to L^\infty(\mu)$ and the induced evolution $\mathcal U^n : B(H) \to B(H)$ from~\eqref{eqHeisenberg}; that is, we have $ \mathcal U^n \circ \pi = \pi \circ U^n$. This compatibility relationship is exhibited by the following commutative diagram:   
\begin{displaymath}
    \begin{tikzcd}
        L^\infty(\mu)\arrow[r, "U^n"] \arrow[d,"\pi",swap] & L^\infty(\mu) \arrow[d,"\pi"] \\
        B(H) \arrow[r,"\mathcal U^n"] & B(H)
    \end{tikzcd}
\end{displaymath}

\subsubsection{Embedding states}

The states of $L^\infty(\mu)$ can similarly be consistently embedded into states of $B(H)$ by means of the map $\Gamma : P(\mu) \to Q(H)$ that maps a probability density $ p $ to the pure state that projects along the square root of $p$, i.e.,
\begin{equation}
    \label{eqGammaP}
    \Gamma(p) = \langle \sqrt p, \cdot \rangle \sqrt p.
\end{equation}
Note that $\sqrt p$ is a unit vector in $H=L^2(\mu)$ since $p$ is a probability density in $L^1(\mu)$. Analogously to $\pi$, the map $\Gamma$ is compatible with the evolution of probability densities under the transfer operator $P^n : P(\mu) \to P(\mu)$ and the induced operator $\mathcal P^n : Q(H) \to Q(H)$ from~\eqref{eqSchrodinger}; that is, we have $\mathcal P^n \circ \Gamma = \Gamma \circ P^n$, and the following diagram commutes:
\begin{displaymath}
    \begin{tikzcd}
        P(\mu)\arrow[r, "P^n"] \arrow[d,"\Gamma",swap] & P(\mu) \arrow[d,"\Gamma"] \\
        Q(H) \arrow[r,"\mathcal P^n"] & Q(H)
    \end{tikzcd}
\end{displaymath}

Overall, we have the following compatibility relationship between classical and quantum time evolution:
\begin{equation}
    \label{eqDynCompat}
    \mathbb E_{\mathcal P^n(\Gamma p)}(\pi f) = \mathbb E_{P^n p} f, \quad \forall f \in L^\infty(\mu), \quad \forall p \in P(\mu). 
\end{equation}
Note that, in general, $\Gamma$ is not an onto map even if one restricts attention to pure states in $Q(H)$; this follows from the fact that there exist unit vectors in $ H$ which are not square roots of positive functions in $L^1(\mu)$. In other words, there are pure states in $Q(H)$ which have no underlying classical probability distribution in $P(\mu)$.

\subsubsection{\label{secEmbeddingEffects}Embedding effects}
Together, the maps $\pi$ and $\Gamma$ provide a consistent representation of Bayesian conditioning of probability densities by events in terms of conditioning of quantum states by effects.    

Let $S \subseteq \Omega$ be a measurable subset of $\Omega$ (an event) and let $ p \in P(\mu)$ be a probability density such that $\int_\Omega p \chi_S \, d\mu > 0$, where $\chi_S$ is the characteristic function of $S$. According to Bayes' theorem, the posterior density to $p$ given $S$ is 
\begin{equation}
    \label{eqBayes}
    p\rvert_S = \frac{\chi_S p}{ \int_\Omega p \chi_S \, d\mu}.
\end{equation}
Interpreting $\chi_S$ as an effect of the abelian algebra $L^\infty(\mu)$ (in the sense that $0 \leq \chi_S \leq 1_\Omega$; cf.~\eqref{eqEffect}), and rewriting $\chi_S p = \chi_S p \chi_S$, it is apparent that~\eqref{eqBayes} and~\eqref{eqQBayes} are manifestations of the same state conditioning procedure, the former being abelian and the latter being non-abelian. In particular, $\pi \chi_S $ is an effect in $\mathcal E(H)$, and we have the compatibility relation
\begin{equation}
    \label{eqMeasCompat}
    \Gamma(p\rvert_S) = (\Gamma p) \rvert_{\pi \chi_S}.
\end{equation}

By~\eqref{eqDynCompat} and~\eqref{eqMeasCompat}, we conclude that the infinite-dimensional quantum system on $H$ consistently reproduces the statistical behavior of the underlying classical dynamical system on $L^1(\mu)$.

\subsection{\label{secDiscretization}Discretization}

Let $ \{ \phi_0, \phi_1, \ldots \}$ be an orthonormal basis of $H$, and let $H_L \subset H$ be the $L$-dimensional subspace defined as 
\begin{equation}
    \label{eqHL}
    H_L = \spn\{ \phi_0, \ldots, \phi_{L-1} \}. 
\end{equation}
Deferring the task of choosing $\{\phi_l\}$ to \cref{secHilb}, we construct finite-dimensional quantum systems on $H_L$ by discretization (i.e., projection) of the infinite-dimensional system on $H$, as follows.

\subsubsection{\label{secDiscretizationObs}Discretization of observables}

Let $\Pi_L : H \to H $ be the orthogonal projection with $ \ran \Pi_L = H_L$. When convenient, we will view $\Pi_L$ as a map into $H_L$ without change of notation. The projection $\Pi_L$ induces a projection $\bm \Pi_L : B(H) \to B(H) $ defined as $\bm \Pi_L A = \Pi_L A \Pi_L $. The range of $\bm \Pi_L$ can be canonically identified with $B(H_L)$, i.e., the space of linear maps on $H_L$, so we can consider $\bm \Pi_L : B(H) \to B(H_L)$ as a map from (bounded) quantum observables on $H$ to quantum observables on $H_L$ (which are necessarily bounded since $H_L$ is finite-dimensional). An important property of $\bm \Pi_L $ is that it is a positive map, i.e., it maps positive operators in $B(H)$ to positive operators in $B(H_L)$. This is in contrast to $\Pi_L$, which in general does not map positive functions to positive functions.  

As a von Neumann algebra, $B(H_L)$ is isomorphic to the $L^2$-dimensional algebra of $L\times L$ complex matrices, which we denote as $\mathbb M_L$. Such an isomorphism $\bm\beta_L : B(H_L) \to \mathbb M_L$ is non-canonical since it depends on a choice of basis for $H_L$, but for our purposes a natural choice is to define $\bm\beta_L$ as the map that yields the matrix representation of operators in $B(H_L)$ with respect to the $\{ \phi_l \} $ basis; that is, we set $\bm\beta_L A = \bm A$, where $\bm A = [A_{ij} ]_{i,j=0}^{L-1} $ is the $L \times L $ matrix with elements $A_{ij} = \langle \phi_i, A \phi_j \rangle $.   

Using the quantum representation of classical observables $\pi : L^\infty(\mu) \to B(H) $ from \cref{secInfDimObs}, we define a projected representation $\pi_L : L^\infty(\mu) \to B(H_L)$ and its matrix-valued counterpart $\bm \pi_L : L^\infty(\mu) \to \mathbb M_L$ as $ \pi_L = \bm \Pi_L \circ \pi$ and $ \bm \pi_L = \bm\beta_L \circ \pi_L$, respectively. By positivity of $\pi$, $\bm \Pi_L$, and $\bm\beta_L$, the maps $\pi_L $ and $\bm \pi_L$ are both positive, i.e., they map positive functions to positive operators and positive matrices, respectively. For a real-valued element $ f \in L^\infty(\mu)$, $\pi_L f$ is a finite-rank symmetric operator with a spectrum $\{ a_0, \ldots, a_J \} \subset \mathbb R $ of real eigenvalues that can be loosely thought of a ``discretization'' of the essential range of $f$. The corresponding spectral measure $E_{\pi_L f} : \mathcal B(\mathbb R) \to B(H_L)$ is given by 
\begin{equation}
    \label{eqSpecMeasF}
    E_{\pi_L f}(S) = \sum_{j: a_j \in S} E_{\pi_L f}^{(j)},
\end{equation}
where $E_{\pi_L f}^{(j)} \in B(H_L)$ is the orthogonal projection onto the eigenspace of $\pi_L f$ corresponding to $a_j$. 

Equation~\eqref{eqSpecMeasF} may be used to compute the measurement probability from~\eqref{eqQProb} given a quantum state $\rho \in Q(H_L)$, viz.
\begin{equation}
    \mathbb P_{\rho, \pi_Lf}(S) = \sum_{j:a_j \in S} \tr\left( \rho E_{\pi_L f}^{(j)}\right).
    \label{eqQProbF}
\end{equation}
As $L \to \infty$, the operators $\pi_L f$ exhibit spectral convergence to the multiplication operator $\pi f $ in a suitable sense; see \cite{FreemanEtAl22} for further details. 

It is important to note that for a general element $ f \in L^\infty(\mu)$, the matrix representation $ \bm \pi_L f $ is not diagonal; that is, the range of $ \bm \pi_L$ is not included in an abelian subalgebra of $\mathbb M_L$ (every such subalgebra would only contain diagonal matrices). This is in contrast to $\pi$ which, being an algebra homomorphism, maps $L^\infty(\mu)$ to an abelian subalgebra of $ B(H)$ (consisting of multiplication operators). Thus, as a result of discretization, our representation of classical observables by the finite-dimensional quantum system exhibits non-abelian behavior, which is an intrinsically quantum mechanical property.     

\subsubsection{Discretization of states} 

For a given quantum state $\rho \in Q(H) $, let $\sigma_0, \sigma_1, \ldots$ be the projected operators $\sigma_L = \bm \Pi_L \rho \in B(H)$. While, in general, the $\sigma_L$ are not density operators, as $L\to \infty $ they converge to $\rho$ in the trace norm of $Q(H)$. As a result, there exists $L_* \in \mathbb N$ such that for every $L > L_*$, we have $C_L = \tr \sigma_L > 0 $, and thus 
\begin{equation}
    \label{eqRhoL}
    \rho_L := \frac{\sigma_L}{C_L}
\end{equation}
is a quantum state in $Q(H)$. Every such state can be identified with a state of $Q(H_L)$. As one can directly verify, for every quantum observable $A$, we have
\begin{equation}
    \label{eqRhoConv}
    \lim_{L\to \infty} \mathbb E_{\rho_L} A_L = \mathbb E_\rho A,
\end{equation}
where $A = \bm \Pi_L A$, so, in the limit of infinite dimension $L$, the quantum systems on $H_L$ consistently recover the expectation of any bounded quantum observable on $H$. 

Let $\beta_L : H_L \to \mathbb C^L$ be the linear map that maps elements of $H_L$ to their column vector representation with respect to the $\{ \phi_l \} $ basis, i.e., $ \beta_L f = \bm f = (f_0, \ldots, f_{L-1})^\top $ with $ f_l = \langle \phi_l, f \rangle$. In applications, we represent a quantum state $\rho \in Q(H_L)$ by an $L\times L $ density matrix $\bm \rho = \bm\beta_L \rho$. If $\rho = \langle \xi, \cdot \rangle \xi$ is a pure state induced by a unit vector $\xi \in H_L$, then $\bm \rho = \bm \xi \bm \xi^\dag$ is a rank-1 projection matrix along the unit vector $ \bm \xi = \beta_L \xi \in \mathbb C^L$, where $^\dag$ denotes the complex conjugate transpose. 

Next, if $ \rho = \Gamma (p) \in Q(H)$ is a pure state from~\eqref{eqGammaP} induced by a probability density $ p \in L^1(\mu)$, i.e., $ \rho = \langle \xi, \cdot \rangle \xi $ with $ \xi = \sqrt{p}$, then for large-enough $L$ we have $ \rho_L = \langle \xi_L, \cdot \rangle \xi_L$ with $ \xi_L = \Pi_L \xi / \lVert \Pi_L \xi \rVert_H$. It should be noted that, in general, $ \xi_L $ is not the square root of a probability density in $P(\mu)$; in fact, $\xi_L$ need not be a positive function. Thus, analogously to our discretization of observables from \cref{secDiscretizationObs}, our finite-dimensional representations of classical probability densities are through intrinsically quantum mechanical objects.      

\subsubsection{\label{secDiscretizationOps}Discretization of evolution operators}

Let $U \equiv U^1$ and $ P \equiv P^1$ be the single-step Koopman and transfer operators on $H$. For any $L \in \mathbb N$, we define projected operators $U_L = \bm \Pi_L U$ and $P_L = \bm \Pi_L P$ on $H_L$, and the corresponding induced operators $\mathcal U_L : B(H_L) \to B(H_L)$ and $\mathcal P_L : B_1(H_L) \to B_1(H_L)$ as $\mathcal U_L A = U_L A U_L^*$ and $\mathcal P_L A = U_L^* A U_L \equiv P_L A P_L^* $. In general,  $U_L$ and $P_L$ are not unitary operators. Nevertheless, $\mathcal P_L$ can be shown to be trace non-increasing, $\tr(\mathcal P_L \rho) \leq 1$ for any $ \rho \in Q(H_L)$. As a result, the operator semigroups $ \{ \mathcal U^n_L \}_{n \in \mathbb N}$ and $ \{ \mathcal P^n_L \}_{n \in \mathbb N}$ generated by $\mathcal U_L$ and $\mathcal P_L$, respectively, define an open quantum system on $H_L$. Note that if $H_L$ is a $U$-invariant subspace of $H$, then $U_L : H_L \to H_L$ and $P_L : H_L \to H_L$ \emph{are} unitary, and we can define associated unitary evolution groups $ \{ \mathcal U^n_L \}_{n \in \mathbb N}$ and $ \{ \mathcal P^n_L \}_{n \in \mathbb N}$. However, every such $U$-invariant subspace $H_L$ must necessarily admit a basis of Koopman eigenfunctions, and we cannot in general assume existence of such a basis. Irrespective of whether $U$ and $P$ are unitary, for any number of timesteps $n \in \mathbb N$, observable $A \in B(H)$, and quantum state $ \rho \in Q(H)$, we have
\begin{displaymath}
    \lim_{L \to \infty} \mathbb E_{\mathcal P^n_L \rho_L} A_L = \mathbb E_{\mathcal P^n \rho} A,
\end{displaymath}
where $A_L = \bm \Pi_L A $ and $ \rho_L$ is given by~\eqref{eqRhoL}.

\subsubsection{\label{secDiscretizationEffect}Discretization of effects}

We use the projection maps $\bm \Pi_L$ to discretize effects similarly to our discretization of quantum observables from \cref{secDiscretizationObs}. First, observe that since $\bm \Pi_L A \leq A $ whenever $A \in B(H)$ is positive, the projections $\bm \Pi_L $ map effects into effects; that is, we can view $\bm \Pi_L$ as a map from $\mathcal E(H)$ to $\mathcal E(H_L)$. Using~\eqref{eqRhoConv}, it follows that for every observable $A \in B(H)$, state $\rho \in Q(H)$, and effect $e \in \mathcal E(H)$,
\begin{displaymath}
    \lim_{L\to\infty} \mathbb E_{\rho_L\rvert e_L} A_L = \mathbb E_{\rho\rvert_e} A,
\end{displaymath}
where $ A_L = \bm \Pi_L A$, $\rho_L$ is given by~\eqref{eqRhoL}, and $e_L = \bm \Pi_L e$. Thus, conditioning by the projected effects $e_L$ consistently recovers conditioning by $e$ in the infinite-dimension limit. Note that if $e\in \mathcal E(H)$ is ``classical'', i.e., it is a multiplication operator by characteristic function of a classical event, $ e = \pi \chi_S $ for some measurable set $S \subseteq \Omega$, the projected effect $e_L = \pi_L \chi_S$ is in general not a multiplication operator. Thus, in the context of the finite-dimensional quantum systems on $H_L$, state conditioning takes place by generally non-classical events. 

\subsection{\label{secHilb}Choice of basis}

Recall that the family of Hilbert spaces $H_0 \subset H_1 \subset \cdots$ from \cref{secDiscretization} is determined from an orthonormal basis $ \{ \phi_0, \phi_1, \ldots \}$ of $H$. In choosing this basis, one must keep in mind that in applications the invariant measure $\mu$ that defines the inner product of $H$ is typically supported on an unknown, non-smooth subset of state space (e.g., a fractal attractor). In such cases, defining the basis vectors $\phi_l$ through explicit closed-form expressions is generally not possible, so we rely instead on an implicit basis construction through eigendecomposition of kernel integral operators. This approach was introduced in the diffusion forecast technique \cite{BerryEtAl15}, and was later used in related methods for Koopman operator approximation \cite{GiannakisEtAl15,Giannakis19,DasEtAl21}, as well as in QMDA \cite{Giannakis19b,FreemanEtAl22}. For our purposes, an advantage of using eigenfunctions of integral operators is their approximability from training data sampled from the invariant measure $\mu$; we will take up this task in \cref{secDataDriven}.

Let $W : \Omega \to \mathcal W$ be a measurable function into a data space $\mathcal W$ that we will use to build our basis. In the applications presented in \cref{secL63,secL96} below, $\mathcal W$ will be a Euclidean space, $\mathcal W = \mathbb R^{d_\mathcal W}$ for some dimension $d_{\mathcal W}$, and possibilities for $W$ will include:
\begin{itemize}
    \item The identity map $\Id$ for $\mathcal W = \Omega$, which corresponds to using information from the full system state to build the basis. 
    \item The projection map $\mathcal P_{\mathcal X} : \Omega \to \mathcal X$ for $\mathcal W = \mathcal X$, which corresponds to building the basis using information from only the resolved variables. 
    \item A delay-coordinate map \cite{SauerEtAl91} for $\mathcal W = \mathcal X^Q$, where $Q = \mathbb N$ is the number of delays, which corresponds to using the dynamics to implicitly recover some of the information about the full system state lost through projection onto $\mathcal X$. 
\end{itemize}

For any of these choices, let $\kappa_{\mathcal W} : \mathcal W \times \mathcal W \to \mathbb R $ be a symmetric kernel function and $\kappa : \Omega \times \Omega \to \mathbb R $ its pullback to $\Omega$, i.e., $\kappa(\omega,\omega') = \kappa_{\mathcal W}(W(\omega),W(\omega'))$. We assume that $W$ and $\kappa_{\mathcal W}$ have sufficient regularity so that $\kappa$ lies in $L^2(\mu\times\mu)$. In that case, the integral operator $K : H \to H$ defined as
\begin{equation}
    \label{eqKOp}
    K f = \int_\Omega \kappa(\cdot, \omega) f(\omega) \, d\mu(\omega) 
\end{equation}
is a self-adjoint, real, Hilbert-Schmidt integral operator. As a result, there exists a real orthonormal basis $\{ \phi_0, \phi_1, \ldots \}$ of $H$ consisting of eigenvectors of $K$; that is,
\begin{equation}
    \label{eqKEig}
    K \phi_l = \lambda_l \phi_l,
\end{equation}
where the eigenvalues $\lambda_l$ are real and satisfy $ \lvert \lambda_0 \rvert \leq \lvert \lambda_1 \rvert \leq \cdots \searrow 0$. In the experiments of \cref{secL63,secL96} our nominal choice for $\kappa_{\mathcal W}$ is the radial Gaussian kernel on $\mathcal W = \mathbb R^{d_{\mathcal W}}$, 
\begin{equation}
    \label{eqKW}
    \kappa_{\mathcal W}(w,w') = e^{-\lVert w - w' \rVert^2_2/\epsilon_{\mathcal W}^2},
\end{equation}
where $\epsilon_{\mathcal W} > 0$ is a bandwidth parameter that we tune automatically (see \cref{appAlgTraining}). In some of our experiments, we will use a variable-bandwidth generalization of~\eqref{eqKW} proposed in \cite{BerryHarlim16}, which was also used in refs.~\cite{Giannakis19b,FreemanEtAl22} in QMDA applications. We define the Hilbert spaces $H_L$ using~\eqref{eqHL} with the eigenvectors $\phi_l$ from~\eqref{eqKEig}. 

Next, note that whenever the eigenvalue $\lambda_l$ is nonzero, any corresponding eigenvector $\phi_l$ (which is an equivalence class of functions on $\Omega$, defined $\mu$-a.e.) is represented by the everywhere-defined function $\varphi_l : \Omega \to \mathbb R  $ such that 
\begin{equation}
    \label{eqVarphi}
    \varphi_l(\omega) = \frac{1}{\lambda_l} \int_\Omega \kappa(\omega,\omega') \phi_l(\omega') \, d\mu(\omega').  
\end{equation}
It follows from~\eqref{eqVarphi} that every such function $\varphi_l$ inherits the regularity properties (e.g., boundedness, continuity, differentiability) of the kernel $\kappa$. This motivates using kernels with plentiful nonzero corresponding eigenvalues in order to control the regularity of the basis functions. More specifically, noticing that $\varphi_l$ from~\eqref{eqVarphi} is the pullback of a function on $\mathcal W$, i.e., $ \varphi_l = \varphi_l^{(\mathcal W)} \circ W$ for some $\varphi_l^{(\mathcal W)} : \mathcal W \to \mathbb R$, we obtain a ``maximal'' number of nonzero eigenvalues $\lambda_l$ if the $ \varphi_l^{(\mathcal W)}$ form an orthonormal basis of the Hilbert space $L^2(\mu_{\mathcal W})$ on $\mathcal W$ associated with the pushforward measure $\mu_{\mathcal W}:= W_* \mu$. This will hold, for instance, if $\kappa_{\mathcal W}$ is an integrally strictly-positive-definite kernel on $\mathcal W$ \cite{SriperumbudurEtAl11}, e.g., the radial Gaussian kernel on $ \mathcal W = \mathbb R^{d_{\mathcal W}}$ from~\eqref{eqKW}.

\section{\label{secMethodology}Quantum mechanical closure (QMCl)}

Recall from \cref{secStatement} that the state space of the parameterized system has the decomposition $\tilde \Omega = \mathcal{X} \times \tilde{\mathcal{Y}}$, where $\tilde{\mathcal{Y}}$ is the state space of the surrogate model for the unresolved degrees of freedom, and $\mathcal{X} $ is the space of variables we are seek to accurately predict. In QMCl, we set $\tilde{\mathcal Y}$ to the state space $\tilde{\mathcal Y} = Q(H_L)$ of the finite-dimensional quantum system constructed in \cref{secDiscretization} for some $L \in \mathbb N$. Thus, in order to evolve the parameterized system via~\eqref{eqPhiDecomp}, we need to specify (i) the surrogate flux term $\tilde Z : Q(H_L) \to \mathcal Z$; and (ii) the evolution map $\tilde \psi : \mathcal X \times Q(H_L) \to Q(H_L)$ of the surrogate model. These maps will be specified in \cref{secFlux,secStateUpdate,secFeatureMap}. In \cref{secInit}, we define a map $\alpha : \Omega \to \tilde \Omega$ that assigns initial conditions of the parameterized system on $\tilde \Omega$ from initial conditions of the true system (see the commutative diagram~\eqref{eqComm}). 

\subsection{\label{secFlux}Flux term} We also assume throughout that $\mathcal Z$ is a Euclidean space, $\mathcal Z  = \mathbb R^d$. With that assumption, we represent the true and surrogate flux terms componentwise, i.e., $ Z = (Z^{(1)}, \ldots, Z^{(d)}) $ and $ \tilde Z = (\tilde Z^{(1)}, \ldots, \tilde Z^{(d)})$, where $Z^{(i)} : \mathcal Y \to \mathbb R$ and $\tilde Z^{(i)} : Q(H_L) \to \mathbb R$ are real-valued functions. We also assume that for each $ i \in \{ 1, \ldots, d \}$, the pullback function $\zeta^{(i)} : \Omega \to \mathbb R$ with $\zeta^{(i)} = Z^{(i)} \circ P_{\mathcal Y}$ is essentially bounded with respect to $\mu$, so that $\zeta^{(i)}$ is an element of the $L^\infty(\mu)$ algebra. Under that condition, $\zeta^{(i)}$ has a representation in $B(H)$ by the multiplication operator $\pi \zeta^{(i)}$, and we have a corresponding projected operator $\pi_L \zeta^{(i)} \in B(H_L)$. We set $\tilde Z^{(i)}$ to the quantum mechanical expectation associated with this multiplication operator, i.e.,  
\begin{equation}
    \label{eqQMFlux}
    \tilde Z^{(i)}(\rho) = \mathbb E_\rho (\pi_L \zeta^{(i)}).
\end{equation}

Given that the resolved variables and the quantum state at time $t_n$ are $x_n \in \mathcal X$ and $\rho_n \in B(H_L)$, respectively, we use~\eqref{eqQMFlux} and the resolved dynamics $\phi : \mathcal X \times \mathcal Z \to \mathcal X$ to update $x_n$ to the state $x_{n+1}$ at time $t_{n+1}$ via the formula
\begin{equation}
    \label{eqXUpdate}
    x_{n+1} = \phi(x_n, z_n), \quad z_n = \tilde Z(\rho_n).
\end{equation}
This step is depicted in the top line in the schematic of \cref{figSchematic}.

\subsection{\label{secStateUpdate}Quantum state update}

Our approach for updating the quantum state $\rho_n \in Q(H_L)$ given the resolved variables in $ x_n \in \mathcal X$ is inspired by the prediction--correction scheme employed in sequential data assimilation (filtering) \cite{MajdaHarlim12,LawEtAl15}, which was formulated quantum mechanically in refs.~\cite{Giannakis19b,FreemanEtAl22}. First, we use the projected transfer operator $\mathcal P_L$ to evolve $\rho_n$ to the quantum state
\begin{equation}
    \label{eqQMPrior}
    \tilde \rho_{n_+1} := \frac{\mathcal P_L \rho_n}{\tr(\mathcal P_L \rho_n)};
\end{equation}
see the arrow labeled ``transfer operator'' in the bottom row of the schematic in \cref{figSchematic}. The state $\tilde \rho_{n+1} $ plays an analogous role to the prior density in classical data assimilation. Then, we update the quantum state $\rho_{n+1}$ at time $t_{n+1}$ by conditioning $\tilde \rho_{n+1}$ by a quantum effect induced by the classical state $x_{n+1}$, i.e., 
\begin{equation}
    \label{eqQMPosterior}
    \rho_{n+1} = \tilde\rho_{n+1} \rvert_{e_{n+1}} = \frac{\sqrt{e_{n+1}} \tilde \rho_{n+1} \sqrt{e_{n+1}} }{\tr(\sqrt{e_{n+1}} \tilde \rho_{n+1} \sqrt{e_{n+1}} )}, \quad e_{n+1} = \mathcal F_L(x_{n+1}), 
\end{equation}
where $\mathcal F_L : \mathcal X \to \mathcal E(H_L) $ is an effect-valued map that will be described in \cref{secFeatureMap}. This update is represented by the arrow labeled ``state conditioning'' in the second row of the schematic in \cref{figSchematic}. The state $\rho_{n+1}$ obtained in this manner is loosely analogous to the posterior density in classical data assimilation. 

The map $\tilde \psi: \mathcal X \times Q(H_L) \to Q(H_L) $ that evolves the state of the quantum system is defined as $\tilde \psi(x_n,\rho_n) = \rho_{n+1} $, where $\rho_{n+1} $ is obtained from $(x_n, \rho_n)$ by combining~\eqref{eqXUpdate}, \eqref{eqQMPrior}, and~\eqref{eqQMPosterior}. 

In applications, the frequency of updating the state via conditioning on a quantum effect can be varied. Instead of using the map $\tilde{\psi}$ as defined above, one may use $\tilde{\psi}_r : \mathcal{X} \times Q(H_L) \to Q(H_L)$ where \eqref{eqXUpdate} and \eqref{eqQMPrior} are iterated $r$ times between every application of \eqref{eqQMPosterior}; see \cref{algQMCl}. This corresponds to lengthening the time under which the quantum component of the system undergoes Koopman evolution before it is updated via conditioning on the observation of an effect.

\subsection{\label{secFeatureMap}Effect-valued feature map}

What remains to completely define the posterior state in~\eqref{eqQMPosterior} is to specify the effect-valued map $\mathcal F_L : \mathcal X \to \mathcal E(H_L) $. Following \cite{FreemanEtAl22}, we build this map as an operator-valued feature map induced by a measurable kernel function $ k : \mathcal X \times \mathcal X \to [0,1] $ on the resolved variables space. Every such function induces a feature map $ F : \mathcal X \to L^\infty(\mu)$ with $\ran F \subseteq L^\infty(\mu)$, given by $F(x) = k(x,\cdot)$. Composing this map with the representation $\pi : L^\infty(\mu) \to B(H)$ leads to a feature map $ \mathcal F : \mathcal X \to \mathcal E(H)$, $ \mathcal F = \pi \circ F$, which takes values in the effect space of the infinite-dimensional operator algebra $ B(H)$. We obtain $\mathcal F_L$ by projection of $\mathcal F$ onto $B(H_L)$, i.e., $\mathcal F_L = \bm \Pi_L \circ \mathcal F $. 

In the applications presented in \cref{secL63,secL96}, $\mathcal X = \mathbb R^{d_\mathcal{X}}$ is a Euclidean space, and we use a radial Gaussian kernel
\begin{equation}
    \label{eqKGauss}
    k(x,x') = e^{-\lVert x - x' \rVert^2_2/ \epsilon^2}, \quad \epsilon>0.
\end{equation}
The radial Gaussian kernel is strictly positive-definite \cite{HofmannEtAl08}, which implies that the corresponding feature map $F$ is injective. 

\subsection{\label{secInit}Initialization} 

Assuming that the parameterized dynamics $\tilde \Phi : \tilde \Omega \to \tilde \Omega$ have been appropriately defined, we would ideally like to have an initialization map $\alpha : \Omega \to \tilde \Omega$ such that the commutative diagram~\eqref{eqComm} is satisfied. Practically, we cannot assume to have access to the unresolved variables in $\mathcal Y$, so we must contend with maps that depend on $(x,y) \in \Omega$ only through $x \in \mathcal X$. With that in mind, it is natural to choose $\alpha$ using the effect-valued feature map from \cref{secFeatureMap}, defining 
\begin{equation}
    \alpha(x) = (x, \rho_x), \quad \rho_x = \frac{\mathcal F_L(x)}{\tr(\mathcal F_L(x))},
    \label{eqInit}
\end{equation}
whenever $\tr(\mathcal F_X(x))>0$. Note that $\rho_x$ defined in this way is a state since $\mathcal F_L(x)$ is a positive operator for all $ x \in \mathcal X$. Empirically, we find that following the decay of initial transients, the behavior of QMCl-parameterized systems does not depend significantly on the choice of $\rho_x$. For instance, replacing~\eqref{eqInit} by 
\begin{equation}
    \alpha(x) = (x, \bar\rho), \quad \bar\rho = \langle 1_\Omega, \cdot \rangle 1_\Omega,
    \label{eqInit2}
\end{equation}
which ignores any information from $x$ for the assignment of the initial quantum state $\bar \rho$, was found to impart negligible changes in the long-term statistical behavior of the parameterized system. The numerical experiments in \cref{secL63,secL96} (which focus on one- and two-point statistics under the invariant measure) utilize $\alpha$ from~\eqref{eqInit2}. 

Computationally, an advantage of \eqref{eqInit2} over \eqref{eqInit} is that the former is a pure state (whereas the latter is generally mixed), resulting in a significant reduction of computational cost (see \cref{secCost}). However, it is possible that initializing with~\eqref{eqInit} (or a low-rank approximation thereof) as opposed to~\eqref{eqInit2} would be beneficial in initial-value prediction experiments, which we do not address in this paper.  

\subsection{\label{secStochasticQMCL}Stochastic parameterization}

The QMCl scheme described in \cref{secFlux,secStateUpdate,secFeatureMap} is deterministic, in the sense that a given initial condition $(x_0,\rho_0) \in (\mathcal X, Q(H_L))$ completely determines the sequence $(x_0, \rho_0), (x_1, \rho_1), \ldots$  of resolved variables $x_n$ and quantum states $\rho_n$. In another approach, which we call the ``stochastic'' approach, we set $\tilde{\mathcal Y}$ to the sequence space $\tilde{\mathcal Y} = \mathcal Z^{\mathbb N}$. As in \cref{secStochasticExamples}, we equip $\tilde{\mathcal Y}$ with a $\sigma$-algebra $\Sigma_{\tilde{\mathcal Y}}$ and a probability measure $\nu$ that is invariant under the shift map $T : \tilde{\mathcal Y} \to \tilde{\mathcal Y}$. We also let $ u : \tilde{\mathcal Y} \to \mathbb Z$ be the projection map onto the first coordinate, and define the random variables $\tilde Z_n : \tilde{\mathcal Y} \to \mathcal Z$ such that $\tilde Z_n = u \circ T^n$. We formally choose the probability space $ (\tilde{\mathcal Y}, \Sigma_{\tilde{\mathcal Y}}, \nu)$ such that the components $Z_n^{(1)}, \ldots, Z_n^{(d)} : \tilde{\mathcal Y} \to \mathbb R $ are independent, real-valued valued random variables distributed according to the measure $\mathbb P_{\rho_n, \pi \zeta^{(j)}} $ from~\eqref{eqQProb}. Operationally, this means that the values $z_n = (z_n^{(1)}, \ldots, z_n^{(d)}) \in \mathcal Z$ of the flux terms at time $n$ are obtained by independent random draws from the measurement distributions $\mathbb P_{\rho_n, \pi \zeta^{(1)}}, \ldots, \mathbb P_{\rho_n, \pi \zeta^{(d)}} $, respectively, determined via~\eqref{eqQProbF}. The update formulas~\eqref{eqQMPrior} and~\eqref{eqQMPosterior} for the quantum state $\rho_n$ remain unchanged. See \cref{algQMClStoch} for pseudocode.     

\section{\label{secDataDriven}Data-driven formulation}

To move towards a data-driven algorithm, we must first establish a method to generate a discrete analog of $H = L^2(\mu)$ (and then of $H_L$) which can be manipulated in practice. To that end, we make the following standing assumptions.
\begin{enumerate}[label=A\arabic*., ref=A\arabic*]
    \item We have access to samples $w_0, \ldots, w_{N-1} \in \mathcal W$ and $ z_0, \ldots, z_{N-1} \in \mathcal Z$ from the map $W$ and the unresolved fluxes $Z$, respectively, with $ w_m = W(\omega_m)$ and $ z_m = Z(\omega_m)$. The samples are taken along a dynamical trajectory $\omega_0, \ldots, \omega_{N-1} \in \Omega$ with $\omega_m = \Phi^m(\omega_0)$.
    \item \label{assumpErgodic}The invariant measure $\mu$ is ergodic. 
    \item \label{assumpInv} There is a compact, forward-invariant set $\mathcal M \subseteq \Omega$ (i.e., $\Phi(\mathcal M) \subseteq \mathcal M$) that contains the support of $\mu$ and the starting point $\omega_0 $ (and thus the entire orbit $\omega_0, \omega_1, \ldots$).
    \item \label{assumpCont} The maps $\Phi$, $X$, $Z$, and $W$ are continuous on $\mathcal M$. Moreover, the kernels $\kappa$ and $k$ from \cref{secHilb,secFeatureMap}, respectively, are both continuous on $\mathcal M \times \mathcal M$.
\end{enumerate}

With these assumptions, we use the samples $\{ (w_m,z_m)\}_{m=0}^{N-1}$ as training data to build data-driven QMCl models that are structurally very similar to the data-independent models described in \cref{secMethodology}. Our construction follows closely the data-driven formulation of QMDA \cite{Giannakis19b,FreemanEtAl22}. Note that we do not assume knowledge of the states $\omega_m$ underlying the training data $(w_m,z_m)$. Also, we do not require that the $\omega_m$ lie on the support of $\mu$ (which may be a null set with respect to an ambient measure on $\Omega$, such as an attractor of a dissipative system), but allow instead initial conditions $\omega_0$ drawn from the larger set $\mathcal M$ of potentially positive ambient measure \cite{Blank17}.   

For the rest of the paper, we reserve $m$ indices to denote samples in the training data and $n$ indices to denote samples obtained from free-running QMCl models. For example, in our notation, $x_m$ and $x_n$ are independent points in $\mathcal X$, the former representing the resolved variables in the training phase, and the latter output from a QMCl model.       

\subsection{Finite-dimensional Hilbert space}
Let $\mu_N$ be the discrete sampling measure on $\Omega$ supported on the training trajectory, defined as $\mu_N = \sum_{m=0}^{N-1} \delta_{\omega_m}/ N$ where $\delta_{\omega_m}$ is the Dirac $\delta$-measure supported on $\omega_m \in \Omega$. As a finite-dimensional analog of $H=L^2(\mu)$, we employ the finite-dimensional Hilbert space $\hat H_N := L^2(\mu_N)$, equipped with the inner product 
\begin{displaymath}
    \langle f, g \rangle_N := \int_\Omega f^* g \, d\mu_N = \frac{1}{N} \sum_{m=0}^{N-1} f^*(\omega_m) g(\omega_m). 
\end{displaymath}
Note that the elements of $\hat H_N$ are equivalence classes of functions on $\Omega$ with common values at the sampled states $\omega_m$. Assuming, for simplicity of exposition, that the points $\omega_0, \ldots, \omega_{N-1}$ are all distinct (which holds true aside from special cases such as $\omega_0$ being an eventually periodic point under $\Phi$), $\hat H_N$ has dimension $N$. Thus, every element $ f \in \hat H_N$ can be represented by a column vector $\bm f = (f_0, \ldots, f_{N-1} )^\top \in \mathbb C^N$ with components $f_m = f(\omega_m)$, and every column vector $\bm f \in \mathbb C^N$ represents a unique element $f\in \hat H_N$. Correspondingly, every linear map $ A:\hat H_N \to \hat H_N$ can be represented by a matrix $\bm A \in \mathbb C^{N\times N}$ such that $\bm A \bm f$ is the column vector representation of $Af $. 

Let $C(\mathcal M)$ denote the Banach space of continuous, complex-valued functions on $\mathcal M$, equipped with the uniform norm, $\lVert  f \rVert_{C(\mathcal M)} = \max_{\omega \in \mathcal M} \lvert f(\omega)\rvert$. Under \cref{assumpErgodic,assumpInv}, for $\mu$-a.e.\ initial condition $\omega_0 \in \mathcal M$ and in the limit of large data, $N \to \infty $, the measures $\mu_N$ converge to the invariant measure in the weak-$^*$ topology of finite Borel measures on $\mathcal M$; i.e., 
\begin{equation}
    \lim_{N\to\infty} \int_{\mathcal M} f \, d\mu_N \equiv \lim_{N\to\infty} \frac{1}{N} \sum_{m=0}^{N-1} f(\omega_m) = \int_{\mathcal M} f \, d\mu, \quad \forall f \in C(\mathcal M).
  \label{eqWeakConv}
\end{equation}
In \cref{eqWeakConv}, we may replace $f$ in the integrals with respect to $\mu_N$ by a uniformly convergent sequence of functions $f_N$; that is, we have 
\begin{equation}
  \label{eqWeakConv2}
  \lim_{N\to\infty} \int_{\mathcal M} f_N \, d\mu_N = \int_{\mathcal M} f \, d\mu, 
\end{equation}
for $\mu$-a.e.\ $\omega_0 \in \mathcal M$, where $\lim_{N\to \infty} f_N = f$ in $C(\mathcal M)$ norm. In what follows, \cref{eqWeakConv2} will play a key role in ensuring the consistency of the data-driven formulation of QMCl.

\subsection{\label{secDataDrivenBasis}Data-driven basis}

In \cref{secHilb}, we used the eigenfunctions $\phi_l$ of the kernel integral operator $K : H \to H$ to define $L$-dimensional subspaces $H_L \subseteq H$ (see~\eqref{eqHL}) on which we built finite-dimensional quantum systems. To do this in the data-driven setting, we replace $K$ with the operator $K_N: \hat{H}_N \to \hat{H}_N$ defined analogously to~\eqref{eqKOp} as 
\begin{displaymath}
    K_N f = \int_\Omega \kappa(\cdot,\omega) f(\omega) \, d\mu_N(\omega) = \frac{1}{N}\sum_{m=0}^{N-1} \kappa(\cdot,\omega_m) f(\omega_m).  
\end{displaymath}
Computationally, $K_N$ is represented by the $N\times N$ kernel matrix
\begin{displaymath}
    \bm K_N = [\kappa(\omega_i, \omega_j)]_{i,j=0}^{N-1} = [\kappa_{\mathcal W}(w_i, w_j)]_{i,j=0}^{N-1}. 
\end{displaymath}
Once again, we solve the eigenvalue problem for $K_N$,
\begin{equation}
    \label{eqKNEig}
    K_N \phi_{l,N} = \lambda_{l,N} \phi_{l,N},
\end{equation}
and define $L$-dimensional subspaces $H_{L,N} \subseteq \hat H_N$ as (cf.~\eqref{eqHL}) 
\begin{equation}
    \label{eqHLN}
    H_{L,N} = \spn\{ \phi_{0,N}, \ldots, \phi_{L-1,N} \}.
\end{equation}

In~\eqref{eqHLN}, the eigenvectors $\phi_{l,N}$ are orthonormal on $\hat H_N$, and the corresponding eigenvalues $\lambda_{l,N}$ are real and are ordered in order of decreasing modulus. Analogously to $\Pi_L : H \to H$ and $\bm \Pi_L : B(H) \to B(H)$ from \cref{secDiscretization}, we define projection maps $\Pi_{L,N} : \hat H_N \to \hat H_N$ and $\bm \Pi_{L,N} : B(\hat H_N) \to B(\hat H_N)$ such that $\ran \Pi_L = H_L$ and $\bm \Pi_{L,N}A = \Pi_{L,N} A \Pi_{L,N}$. We also let $\beta_{L,N} : H_{L,N} \to \mathbb C^L$ and $\bm\beta_{L,N} : B(H_{L,N}) \to \mathbb M_L$ be the maps from vectors in $H_{L,N}$ and operators in $B(H_{L,N})$ to their corresponding vector and matrix representations with respect to the $\{ \phi_{l,N} \}_{l=0}^{L-1}$ basis, defined analogously to $\beta_L : H_L \to \mathbb C^L$ and $\bm\beta_L : B(H_L) \to \mathbb M_L$ from \cref{secDiscretization}, respectively.

Similarly to the representatives $\varphi_l $ of $\phi_l$ from~\eqref{eqVarphi}, every eigenvector $\phi_{l,N}$ with nonzero corresponding eigenvalue $\lambda_{l,N}$ has an everywhere-defined representative $\varphi_l : \Omega \to \mathbb R$, given by
\begin{equation}
    \varphi_{l,N}(\omega) = \frac{1}{\lambda_{l,N}} \int_{\mathcal M} \kappa(\omega,\omega') \phi_{l,N}(\omega) \, d\mu_N = \frac{1}{\lambda_{l,N}} \frac{1}{N} \sum_{m=0}^{N-1} \kappa(\omega,\omega_m) \phi_{l,N}(\omega_m).
    \label{eqVarphiN}
\end{equation}
Under \cref{assumpCont}, every such $\varphi_l$ and $ \varphi_{l,N}$ is continuous on $\mathcal M$. By results on spectral approximation of kernel integral operators \cite{VonLuxburgEtAl08}, the following can be shown to hold as $N\to\infty $, for $\mu$-a.e.\ initial condition $\omega_0 \in \mathcal M$:
\begin{enumerate}
    \item For every nonzero eigenvalue $\lambda_l $ of $K$, the sequence of eigenvalues $\lambda_{l,N} $ of $K_N$ converges to $\lambda_l$, including multiplicities.
    \item For every continuous representative $\varphi_l \in C(\mathcal M) $ of an eigenfunction $\phi_l \in H$ of $K$ corresponding to $\lambda_l \neq 0$, there exists a sequence of eigenfunctions $\phi_{l,N} \in \hat H_N$ of $K_N$ whose continuous representatives $\varphi_{l,N} \in C(\mathcal M)$ converge to $ \varphi_l$ in the $C(\mathcal M)$ norm. 
\end{enumerate}
We refer the reader to the paper~\cite{FreemanEtAl22} for further details on these results in the context of QMDA. 

\subsection{Operator approximation}

We use the data-driven basis $\{ \phi_{l,N} \}_{l=0}^{N-1}$ of $H_{L,N}$ to consistently approximate matrix representations of operators on $H_L$ by matrix representations of operators on $H_{L,N}$ for a class of operators that behave consistently with operators on continuous functions. In what follows, $\iota : C(\mathcal M) \to L^p(\mu)$ and $\iota_N : C(\mathcal M) \to L^p(\mu_N)$ will denote the standard maps from continuous functions on $\mathcal M$ to their corresponding $L^p$ equivalence classes with respect to $\mu$ and $\mu_N$, respectively. 

Consider an operator $A\in B(H)$ that satisfies
\begin{equation}
  A \circ \iota = \iota \circ \tilde A
  \label{eqOpCompat}
\end{equation}
for some bounded operator $\tilde A : C(\mathcal M) \to C(\mathcal M)$. Consider also a uniformly bounded sequence of operators $\hat A_1, \hat A_2, \ldots$ in $B(\hat H_1), B(\hat H_2), \ldots $, respectively, such that 
\begin{equation}
  \lim_{N\to\infty} \lVert (\hat A_N \circ \iota_N - \iota_N \circ \tilde A) f \rVert_{\hat H_N} =0, \quad \forall f \in C(\mathcal M).
  \label{eqOpCompat2}
\end{equation}
Given any such operator family $\{ A, \tilde A, \hat A_0, \hat A_1, \ldots \}$, our approach is to approximate matrix elements of $A$ with respect to the kernel eigenbasis $\{ \phi_l \} $ of $H$ from \cref{secHilb} by matrix elements of $ \hat A_N$ with respect to the data-driven eigenbasis $\{ \phi_{l,N} \}$ of $\hat H_N$. Specifically, let $\phi_i$ and $\phi_j$ be two basis vectors of $H$ from~\eqref{eqKEig} corresponding to nonzero eigenvalues. Let $\phi_{i,N}$ and $\phi_{j,N}$ be basis vectors of $\hat H_N$ from~\eqref{eqKNEig} chosen such that, as $N\to \infty$, their continuous representatives $\varphi_{i,N}$ and $\varphi_{j,N}$ converge to those of $\phi_i$ and $\phi_j$ (i.e., $\varphi_i$ and $\varphi_j$), respectively, as described in \cref{secDataDrivenBasis}. Then, it can be shown \cite{FreemanEtAl22} that for $\mu$-a.e.\ $\omega_0 \in \mathcal M$, 
\begin{equation}
    \label{eqAConv}
    \lim_{N\to\infty} \langle \phi_{i,N}, \hat A_N \phi_{j,N} \rangle_N = \langle \phi_i, A \phi_j \rangle.
\end{equation}
This means that we can consistently approximate matrix elements of $A$ by matrix elements of $\hat A_N$. In particular, if $L \in \mathbb N $ is such that $\lambda_{L-1} $ is nonzero, it follows from~\eqref{eqAConv} that as $N\to\infty$ the matrix representations $\bm A_{L,N} = \bm\beta_{L,N} \hat A_{L,N}$ of the projected operators $ A_{L,N} = \bm \Pi_{L,N}\hat A_N \in B(H_{L,N})$ converge to the matrix representation $\bm A_L = \bm\beta_L A_L $ of $A_L = \bm \Pi_L A$ in any matrix norm. As we explain in \cref{secDataDrivenQMCl,appDataDriven}, under \cref{assumpCont}, all operators employed in QMCl satisfy the compatibility conditions \eqref{eqOpCompat} and~\eqref{eqOpCompat2} and thus can be approximated in this manner. 

\subsection{\label{secDataDrivenQMCl}Data-driven QMCl framework}

In the data-driven setting, the state space of the quantum mechanical model of the unresolved variables is $\tilde{\mathcal Y} = Q(H_{L,N})$. Analogously to the data-independent formulation in \cref{secMethodology}, we evolve the parameterized system on $ \mathcal X \times \tilde{\mathcal Y}$ using a surrogate flux $ \tilde Z : Q(H_{L,N}) \to \mathcal Z$ and an evolution map $\tilde \psi : \mathcal X \to Q(H_{L,N}) \to Q(H_{L,N}) $. In this subsection, we give an outline of the construction of these maps, focusing on the differences between the data-driven approach and the formulation of  \cref{secMethodology}. Specific formulas and pseudocode relevant to the data-driven setting are included in \cref{appDataDriven}.

\paragraph*{Flux terms}

As in \cref{secFlux}, we prescribe the surrogate flux $\tilde Z : Q(H_{L,N}) \to \mathcal Z$ using quantum mechanical expectations of discrete multiplication operators. For that, we first note that the space $L^\infty(\mu_N)$ is a finite-dimensional, abelian von Neumann algebra, which we use as a data-driven analog of $L^\infty(\mu)$ (see \cref{secVonNeumann}). This algebra has a regular representation $\hat \pi_N : L^\infty(\mu_N) \to B(\hat H_N)$ that maps each vector $ f \in \hat H_N$ to the discrete multiplication operator that multiplies by $f$. Note that $\hat \pi_N f$ is a diagonal operator in the standard basis of $\hat H_N$. We use $\hat \pi_N$ and the projected representation $\pi_{L,N} := \bm \Pi_{L,N} \circ \hat \pi_N $ as data-driven analogs of $\pi : L^\infty(\mu) \to B(H) $ and $\pi_L : L^\infty(\mu) \to B(H_L) $, respectively. The training samples $z_m \in \mathcal Z \equiv \mathbb R^d$ define, componentwise, a collection of elements $\zeta^{(1)}_N, \ldots, \zeta^{(d)}_N \in L^\infty(\mu_N)$ such that $\zeta^{(i)}_N(\omega_m) = z_m^{(i)}$ where $z_m^{(i)}$ is the $i$-th component of $z_m = ( z_m^{(1)}, \ldots, z_m^{(d)})$. Analogously to~\eqref{eqQMFlux}, we define $\tilde Z = (\tilde Z^{(1)}, \ldots, \tilde Z^{(d)})$ with $\tilde Z^{(i)}(\rho) = \mathbb E_\rho(\pi_{L,N} \zeta^{(i)}_N)$. The evolution of the resolved variables in $\mathcal X$ given $\rho \in Q(H_{L,N})$ is carried out via~\eqref{eqXUpdate}. Computationally, the quantum observables $\tilde Z^{(i)}$ and states $\rho$ are represented by their $L\times L$ matrix representations in the $\{ \phi_{l,N} \}$ basis of $H_{L,N}$, i.e., $\bm Z^{(i)}_{L,N} = \bm \beta_{L,N} \tilde Z^{(i)}$ and $ \bm \rho = \bm \beta_{L,N} \rho$.         

\paragraph*{Quantum state update}

We employ a predictor--corrector scheme similar to that in \cref{secStateUpdate}. Given the quantum state $\rho_n \in Q(H_{L,N})$ and resolved variables $x_n \in \mathcal X$ at time $t_n$, we use a data-driven approximation $\mathcal P_{L,N} : B_1(H_{L,N}) \to B_1(H_{L,N})$ of the transfer operator to obtain the prior state $\tilde \rho_{n+1}$ at time $t_{n+1}$ analogously to~\eqref{eqQMPrior} and an effect-valued feature map $\tilde{\mathcal F}_{L,N}: \mathcal X \to \mathcal E(H_{L,N})$ that updates $\tilde \rho_{n+1}$ by conditioning by $ e_{n+1} = \tilde{\mathcal F}_{L,N}(x_{n+1})$ as in~\eqref{eqQMPosterior}. The transfer operator $\mathcal P_{L,N}$ is based on an approximation of the Koopman operator $U : H \to H$ by a shift operator $\hat U_N : \hat H_N \to \hat H_N$ \cite{BerryEtAl15}. That is, we have $\mathcal P_{L,N} A = U_{L,N}^* A \bm U_{L,N} $ where $ U_{L,N} = \bm \Pi_{L,N} \hat U_N$, and $U_{L,N}$ is represented by the $L \times L$ matrix $\bm U_{L,N} = \bm\beta_{L,N} U_{L,N}$; see \cref{appDataDrivenOps} for further details. The effect-valued feature map $\mathcal F_{L,N}$ is constructed analogously to $\mathcal F_L$ from \cref{secFeatureMap} using the radial Gaussian kernel in~\eqref{eqKGauss}. In the $ \{ \phi_{l,N} \}$ basis of $H_{L,N}$, the map $\tilde{\mathcal F}_{L,N}$ is represented by a matrix-valued map $\tilde{\bm F}_{L,N} : \mathcal X \to \mathbb M_L$ with $\tilde{\bm F}_{L,N} = \bm \beta_{L,N} \circ \mathcal F_{L,N}$ that we use in numerical applications; see \cref{appDataDrivenEffect}.  

\paragraph*{Initialization}

In the experiments presented in \cref{secL63,secL96}, we initialize the parameterized system with an uninformative quantum state analogously to \eqref{eqInit2}. Specifically, given any $x \in \mathcal X$, we set the initial state $\alpha(x) = (x, \bar\rho_{N,L}) \in \mathcal X \times Q(H_{L,N})$, where
\begin{equation}
    \label{eqRhoInitDataDriven}
    \bar \rho_{N,L} = \frac{\bm \Pi_{L,N} \bar \rho_N}{\tr(\bm \Pi_{L,N}\bar\rho_N)}, \quad \bar\rho_N= \langle 1_\Omega, \cdot \rangle_N 1_\Omega.
\end{equation}
An alternative approach would be to set $\alpha(x) = (x,\rho_x)$, where the quantum state $\rho_x$ is obtained via the effect-valued feature map $\mathcal F_{L,N}$ analogously to~\eqref{eqInit}. As mentioned in \cref{secInit}, numerically we found that the uninformative initialization approach has negligible impact on the ability of the parameterized system to reproduce the equilibrium statistical behavior of the original system, but initializing with~\eqref{eqInit} is expected to be important in initial-value prediction experiments.

\paragraph*{Stochastic parameterization} The data-driven QMCl framework has a stochastic variant which is entirely analogous to the scheme described in \cref{secStochasticQMCL}.

\subsection{\label{secDataDrivenConvergence}Convergence of data-driven approximation}

With the addition of another approximation parameter ($N$ in addition to $L$), we seek to examine the convergence properties of the system under the iterated limits of $L \to \infty$ after $N \to \infty$. Previous work \cite{FreemanEtAl22} has shown that for quantum observables $\hat A_N \in B(\hat H_N)$ and $A \in B(H)$ satisfying~\eqref{eqOpCompat} and~\eqref{eqOpCompat2} for some $\tilde A : C(\mathcal M) \to C(\mathcal M)$ and for states $\hat\rho_N \in Q(\hat H_N)$ and $\rho \in Q(H)$ satisfying a related compatibility condition with operators on continuous functions, the following asymptotic consistency relationship holds, 
\begin{equation}
    \label{eqConv1}
    \lim_{L \to \infty} \lim_{N\to\infty} \mathbb E_{\mathcal P_{L,N}\rho_{L,N}} A_{L,N} = \lim_{L\to\infty} \mathbb E_{\mathcal P_L \rho_L} A_L = \mathbb E_\rho A.
\end{equation}
Here, $A_{L,N} = \bm \Pi_{L,N} \hat A_N$ and $A_L = \bm \Pi_L A$ are the projected quantum observables associated with $\hat A_N$ and $A$, respectively, and $\rho_{L,N} \in Q(H_{L,N})$ and $\rho_L \in Q(H_L)$ are the projected states associated with $\hat \rho_N$ and $\hat \rho$, respectively (see~\eqref{eqRhoL}). Under \cref{assumpCont}, the fluxes $\zeta^{(i)}$ are continuous, which implies that \eqref{eqOpCompat} and~\eqref{eqOpCompat2} are satisfied with $A = \pi (\iota \zeta^{(i)})$, $\hat A_N = \hat\pi_N(\iota_N \zeta^{(i)})$, and $\tilde A$ set to the multiplication operator by $\zeta^{(i)}$ on continuous functions. Furthermore, the class of states $\rho$ for which~\eqref{eqConv1} holds includes images $\rho = \Gamma(p)$ from~\eqref{eqGammaP} of probability densities in $L^1(p)$ with continuous representatives in $C(\mathcal M)$, as well as higher-rank generalizations, so the data-driven QMCl formulation is asymptotically consistent as $N\to\infty$ in a broad range of scenarios encountered in applications. In addition, an analogous convergence result holds for conditioning by effects $e \in \mathcal E(H)$ and $\hat e_N \in \mathcal E(\hat H_N)$ which satisfy~\eqref{eqOpCompat} and~\eqref{eqOpCompat2} for some operator $ \tilde e : C(\mathcal M) \to C(\mathcal M)$ on continuous functions. We refer the reader to \cite{FreemanEtAl22} for further details. 

\subsection{\label{secCost}Computational cost}

The training data requirements and computational cost of QMCl are generally comparable with those of kernel methods for supervised machine learning. Given that the data space $\mathcal W$ has dimension $d_{\mathcal W}$, the brute-force computation cost of forming the $N\times N$ kernel matrix $\bm K_N$ representing the integral operator $K_N$ is $O(d_{\mathcal W}N^2)$ for radial kernels. This cost can be reduced to $O(N\log N)$ in data spaces of sufficiently low dimension using randomized methods for approximate nearest neighbors, e.g., \cite{JonesEtAl11}. In our numerical experiments, we compute $\bm K_N$ with brute force and then sparsify it, retaining $k_\text{nn} \ll N$ nearest neighbors per data point. The storage cost and matrix--vector multiplication cost for $\bm K_N$ then become $O(k_\text{nn} N)$. We compute the basis vectors $\phi_{l,N}$ using iterative solvers. The cost of this computation depends on the spectral properties of $\bm K_N$ and the number $L$ of requested eigenvectors, but generally scales linearly $k_\text{nn}$ and $N$. Once the basis $\{ \phi_{l,N} \}_{l=0}^{L-1}$ has been computed, we form the $L \times L$ matrix $\bm U_{L,N} = \bm \pi_{L,N} \hat U_N$ representing the projected shift operator and the $L\times L$ observable matrices $ \bm Z^{(1)}_{L,N}, \ldots, \bm Z^{(d)}_{L,N}$ (see \cref{appAlgTraining}), each with an $O(NL^2)$ computational cost. This completes the training phase of QMCl.            

The computational cost associated with advancing the parameterized model over one timestep, given that the classical and quantum states are $x_n \in \mathcal X $ and $\rho_n \in Q(H_{L,N})$, respectively, is as follows:
\begin{itemize}
    \item We compute the flux terms $\tilde Z^{(i)} (\rho_n) = \tr( \bm \rho_n \bm Z_{L,N}^{(i)})$ for $i \in \{ 1, \ldots, d \}$, where $\bm \rho_n \in \mathbb M_L$ is the matrix representation of $\rho_n$ in the $\{\phi_{l,N}\}$ basis of $H_{L,N}$. For a quantum state of rank $r$, the cost of each of these computations is $O(r L^2)$. This can be as high as $O(L^3)$ for quantum states of full rank, but in our experiments we work with pure states, $r=1$, which results in $O(L^2)$ operations.
    \item We advance the resolved variables $x_n$ via~\eqref{eqXUpdate} using the previously computed fluxes. The cost of this procedure is independent of QMCl so we do not consider it further here.
    \item Using the transfer operator, we advance $\bm \rho_n$ to the density matrix $\tilde{\bm\rho}_{n+1}$ representing the prior state $\tilde \rho_n$. This is again an $O(rL^2)$ operation which can be as high as $O(L^3)$ but reduces to $O(L^2)$ for pure states.
    \item We compute the $L\times L$ matrix $\bm F_{L,N}(x_{n+1})$ representing the quantum effect $\mathcal F_{L,N}(x_{n+1})$. This has an $O(d_{\mathcal X}NL^2)$ cost. Using $\bm F_{L,N}(x_{n+1})$ we condition $\tilde{\bm \rho}_{n+1}$ to obtain the density matrix $\bm \rho_{n+1}$ representing the posterior state from~\eqref{eqQMPosterior}. The cost of this operation is $O(rL^2)$, where $r$ is again the rank of $\tilde{\bm \rho}_{n+1}$.
\end{itemize}
Note that the density matrix update, $\tilde{\bm \rho}_{n+1} \mapsto \bm \rho_{n+1}$, is the only step in the online prediction phase whose cost depends on the size $N$ of the training data. In \cref{secFutureWork}, we discuss possible ways of alleviating this dependence using random feature methods for kernel matrix approximation \cite{RahimiRecht07}.

\section{\label{secL63}Quantum mechanical closure of the L63 system}

As our first set of numerical examples, we apply QMCl to the L63 system \cite{Lorenz63}. The L63 system is classically defined as the dynamical system on $\mathbb{R}^3$ where $(x(t), y(t), z(t)) \in \mathbb{R}^3$ evolves as 
\begin{displaymath}
    \dot x(t) = \sigma(y(t) - x(t)), \quad \dot y(t) = x(t) (\rho - z(t)) - y(t), \quad \dot z(t) = x(t) y(t) - \beta z(t)
\end{displaymath}
for the parameter values $\beta = 8/3$, $\rho = 28$, and $\sigma = 10$. However, as noted in \cite{Palmer01}, the L63 system can be expressed in terms of the system~\eqref{eqL63Palmer}, where the variables $(a_1(t), a_2(t), a_3(t)) \in \Omega \equiv \mathbb R^3 $ are obtained by projection of $(x(t),y(t),z(t))$ onto the EOF basis vectors, and have decreasing variance.  

\subsection{Experimental setup}

We follow closely the setup of Palmer~\cite{Palmer01}, which we outlined in \cref{secStochasticExamples}. We assume knowledge of the equations governing the $a_1$ and $a_2$ components, but not of that governing the $a_3$ component. That is, we have $(a_1, a_2) \in \mathcal X \equiv \mathbb R^2$, $a_3 \in \mathcal Y \equiv \mathbb R$, and the flux term is $ Z : \mathcal Y \to \mathcal Z \equiv  \mathbb R $ with $Z(a_3) = a_3$. The true discrete-time system on $\Omega$ evolves under the time-$\Delta t$ flow generated by~\eqref{eqL63Palmer} for a timestep of $\Delta t = 0.01$; that is, we have $\Phi : \Omega \to \Omega $ with $ \Phi = \Phi^{0.01}$. Of course, $\Phi$ is not available in closed form, so practically we consider as the ``true'' L63 dynamics a numerical approximation of $\Phi$ by a high-fidelity ordinary differential equation solver. Here, we use MATLAB's \texttt{ode45} solver which is based on an adaptive Runge-Kutta scheme of order $(4,5)$. As our approximate resolved dynamics $\tilde \phi : \mathcal X \times \mathcal Z \to \mathcal Z$ we use a standard 4th-order Runge-Kutta (RK4) discretization of the $a_1$ and $a_2$ equations in~\eqref{eqL63Palmer}, treating $a_3$ as fixed; see~\eqref{eqRK4} for an explicit formula. We evolve the quantum states in $ \tilde{\mathcal Y} = Q(H_{L,N})$ using the map $\psi_r : \mathcal X \times Q(H_{L,N}) \to Q(H_{L,N})$ from \cref{secStateUpdate} for various choices of the number of timesteps $r$ between state updates by the quantum Bayes' rule~\eqref{eqQMPosterior}.  

\subsection{\label{secL63Big}Experiments with large, full training data}

We first consider QMCl models trained with a long time series of the full system state. In these experiments, the training observation map $ W : \Omega \to \mathcal W$ is the identity map on $\mathcal W = \Omega = \mathbb R^3$, and we use a training time series $\omega_0, \ldots, \omega_{N-1} \in \mathcal W$ consisting of $N=\text{150,000} $ samples with $\omega_m = (x_m,y_m)= \Phi^m(\omega_0)$. The initial training state $\omega_0$ is obtained by integrating~\eqref{eqL63Palmer} with the initial condition $(2,2,2)$ for 500 model time units (i.e., $500/\Delta t = \text{50,000}$ timesteps), and setting $\omega_0$ to the final point of that trajectory. We also use the values $z_0, \ldots, z_{N-1} \in \mathbb R$ of the flux term, $z_m = Z(y_m) = y_m$, from the training trajectory.  

Setting the dimension parameter $L = \text{1500}$, we compute the kernel eigenfunction basis $ \{ \phi_{l,N} \}_{l=0}^{L-1}$ of $H_{L,N}$ from~\eqref{eqKNEig} using the training data $\omega_m$. The kernel function $\kappa_{\mathcal W}$ is the radial Gaussian kernel~\eqref{eqKW} with the bandwidth parameter $\epsilon_{\mathcal W} = \sqrt{19}$ chosen via automatic tuning (see~\cref{appAlgTraining}). We then compute the corresponding $L \times L$ matrix representations $\bm U_{L,N} $ and $ \bm Z_{L,N} $ of the projected Koopman operator and the projected multiplication operator representing the flux term $Z$ as described in \cref{secDataDrivenQMCl}. We also use the basis functions to build the effect-valued feature map $\tilde{\bm F}_{L,N}$ described in \cref{appDataDrivenEffect}. This map is based on the radial Gaussian kernel $k$ in \eqref{eqKGauss} with the bandwidth parameter $\epsilon=2$. This completes the training phase of QMCl. 

To run the parameterized system, we generate a state $\hat\omega_0 \in \Omega$ near the Lorenz attractor analogously to $\omega_0$; that is, we integrate~\eqref{eqL63Palmer} over 500 model time units with the initial condition $(1.99,2,2)$, and set $\hat\omega_0 =( \hat x_0, \hat y_0)$ to the last point of that trajectory. We then use the initialization map~\eqref{eqInit2} to generate an initial condition $\iota(\hat x_0) = (\hat x_0, \hat \rho_0) \in \mathcal X \times Q(H_{L,N})$ of the parameterized system. We condition the quantum state via the $\tilde{\mathcal F}_{L,N}$ effect map every $r=10$ timesteps (i.e., every $ r \, \Delta t = 0.1$ model time units). 

Starting from $(\hat x_0, \hat \rho_0)$, the QMCl system generates via Algorithm \ref{appAlgClosure} a time-ordered sequence of pairs $(\hat x_0, \hat \rho_0), (\hat x_1, \hat \rho_1), \ldots$ of resolved variables $\hat x_n \in \mathcal X_n$ and quantum states $\hat \rho_n \in Q(H_{L,N})$, as well as a corresponding sequence of flux terms $\hat z_0, \hat z_1, \ldots \in \mathcal Z$ given by $\hat z_n = \tilde Z(\rho_n)$ in accordance with \eqref{eqQMFlux}. Under a ``perfect'' closure in the sense of the commutative diagram~\eqref{eqComm}, the sequence of resolved variables $\hat x_0, \hat x_1, \ldots$ should match the time series $\tilde x_0, \tilde x_1, \ldots$ of the $(a_1,a_2)$ state vector components under the L63 flow starting from the same initial condition, i.e., $\tilde x_n = (a_1(t_n), a_2(t_n))$ with $(a_1(t_n), a_2(t_n), a_3(t_n) ) = \Phi^{n\,\Delta t}(\hat \omega_0)$. If, in addition, the flux term $Z=a_3$ is consistently approximated by $\tilde Z$, then the time series $(\hat x_0, \hat z_0 ), (\hat x_1, \hat z_1), \ldots $ generated by the QMCl system should match the full three-dimensional L63 trajectory $\tilde \omega_0, \tilde \omega_1, \ldots$ with $\tilde \omega_n = (a_1(t_n),a_2(t_n),a_3(t_n))$. Here, we assess the performance of the scheme by examining its ability to reproduce salient qualitative features of the Lorenz attractor and to recover the marginal distributions and time-autocorrelation functions of $a_1$, $a_2$, and $a_3$.

\begin{figure}
    \begin{subfigure}{0.49\linewidth}
  \centering
\hspace*{0cm}\vspace*{0cm}\includegraphics[width=\linewidth]{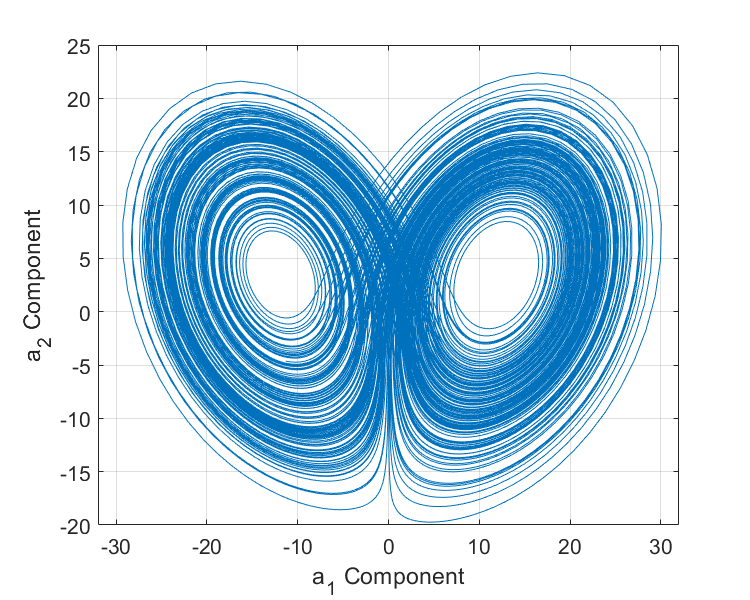}
\caption{True L63 dynamics}
  \label{figRealProj2}
\end{subfigure}
\begin{subfigure}{0.49\linewidth}
\centering
\hspace*{0cm}\vspace*{0cm}\includegraphics[width=\linewidth]{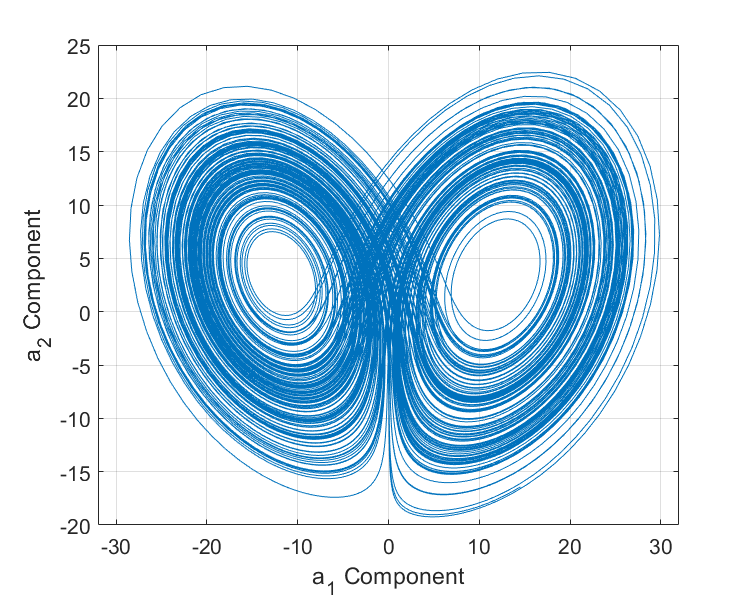}
\caption{QMCl system}
  \label{figQMProj2}
\end{subfigure}

\caption{Trajectory plots of the resolved variables $x=(a_1,a_2)$ under (a) the true L63 dynamics and (b) the QMCl parameterized system. The length of the plotted trajectories is 200 model time units.}
\label{figCoordinateProjections}

\end{figure}

\begin{figure}
\centering
    \begin{subfigure}{0.49\linewidth}
\hspace*{0cm}\vspace*{0cm}\includegraphics[width=\linewidth]{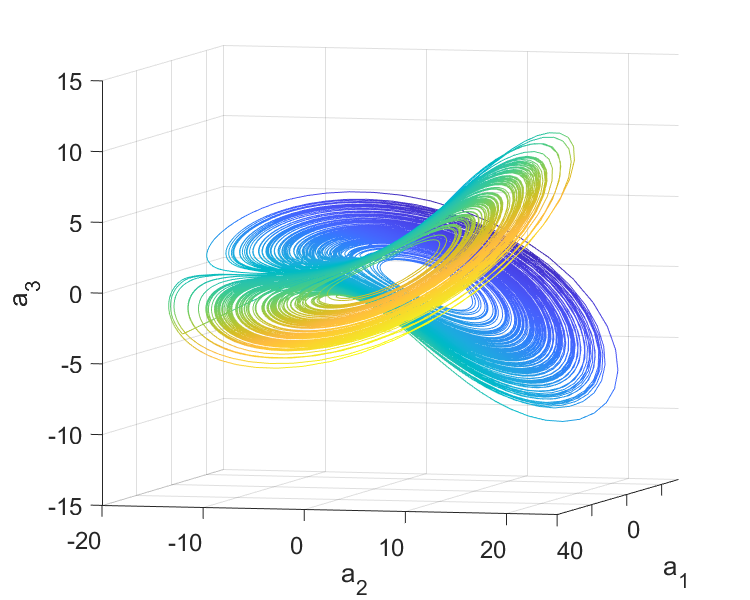}
\caption{True L63 dynamics}
  \label{figRealProj1}
\end{subfigure}
\begin{subfigure}{0.49\linewidth}
  \centering
\hspace*{0cm}\vspace*{0cm}\includegraphics[width=\linewidth]{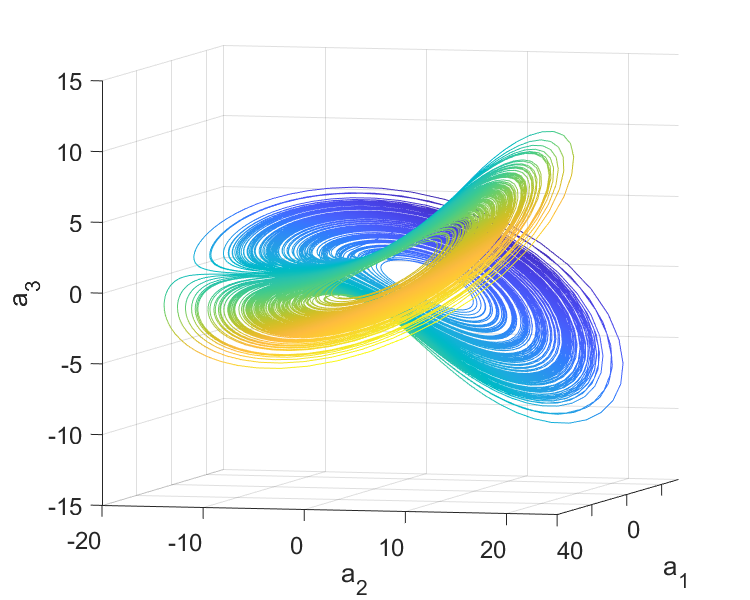}
\caption{QMCl system}
  \label{figQMProj1}
\end{subfigure}
\caption{As in \cref{figCoordinateProjections}, but for trajectories of $(a_1,a_2,a_3)$ in $\mathbb R^3$ showing both the resolved components, $x = (a_1,a_2)$, and the unresolved component, $ y = a_3$, of the L63 state vector. The plotted trajectories spans 200 model time units. Colors correspond to the value of the $a_1$ component, for visual clarity.}

\label{fig3DProjections}

\end{figure}

\begin{figure}
\centering
    \begin{subfigure}{0.49\linewidth}
\hspace*{0cm}\vspace*{0cm}\includegraphics[width=\linewidth]{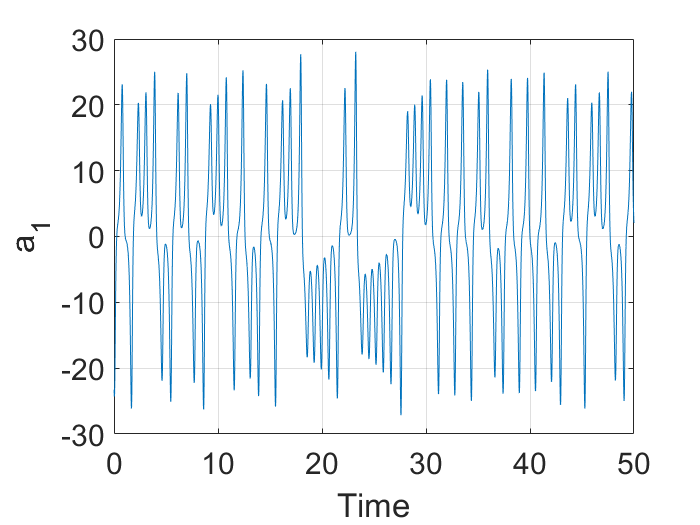}
\caption{True L63 dynamics}
  \label{figRealTimeseries1}
\end{subfigure}
\begin{subfigure}{0.49\linewidth}
  \centering
\hspace*{0cm}\vspace*{0cm}\includegraphics[width=\linewidth]{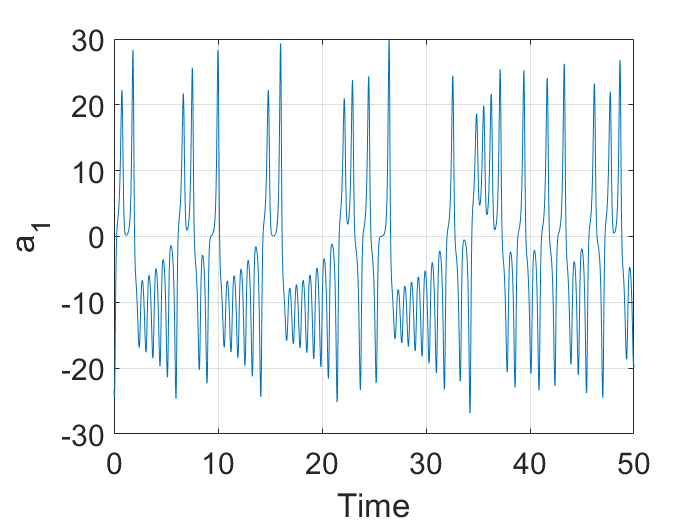}
\caption{QMCl system}
  \label{figQMTimeseries1}
\end{subfigure}
\caption{Time series plots of the $a_1$ component under (a) the true L63 system and (b) the QMCl parameterized system.}

\label{TimeseriesPlots}

\end{figure}

\cref{figCoordinateProjections} compares representative trajectories of the resolved variables $x$ under the true L63 dynamics (\cref{figRealProj2}) and their approximations $\hat x_n$ under the QMCl system (\cref{figQMProj2}). It is readily apparent that the QMCl system generates a structure with a similar geometry to the Lorenz attractor. It is also worth noting that this similarity holds not just in the $a_1, a_2$ projection (which would be sufficient for a parameterization scheme). In fact, as shown by the three-dimensional trajectory plots in \cref{fig3DProjections}, the QMCl algorithm meaningfully reconstructs all three dimensions of the system. Particularly, the trajectories show that the QMCl system evolves in the $a_3$ component in a qualitatively similar manner to the true system. In \cref{TimeseriesPlots}, we show time series plots for the $a_1$ component under the true L63 and QMCl dynamics. The general rate and aperiodicity of the transitions, as well as the tendency for the point to sometimes oscillate in one lobe with increasing amplitude prior to switching lobes, are preserved in the QMCl system. 

\begin{figure}
\begin{subfigure}{0.49\linewidth}
\centering
\hspace*{0cm}\vspace*{0cm}\includegraphics[width=\linewidth]{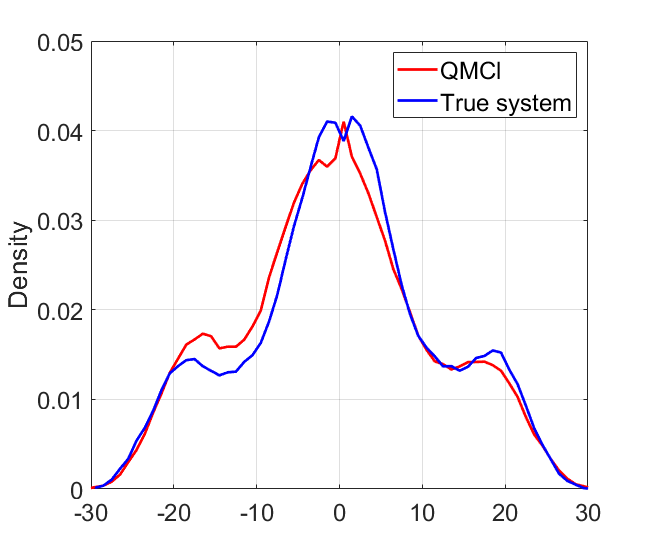}
\caption{$a_1$ histograms}
  \label{figHist1}
\end{subfigure}
\medskip
\begin{subfigure}{0.49\linewidth}
\centering
\hspace*{0cm}\vspace*{0cm}\includegraphics[width=\linewidth]{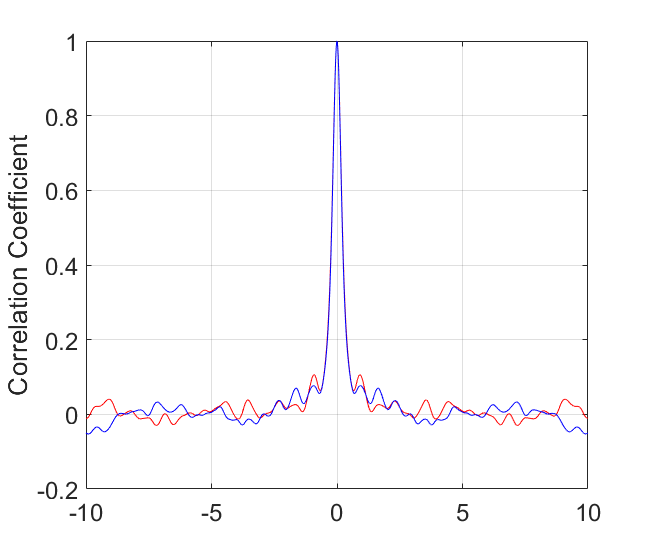}
\caption{$a_1$ autocorrelation functions}
\label{figAutocorr1}
\end{subfigure}

\begin{subfigure}{0.49\linewidth}
\centering
\hspace*{0cm}\vspace*{0cm}\includegraphics[width=\linewidth]{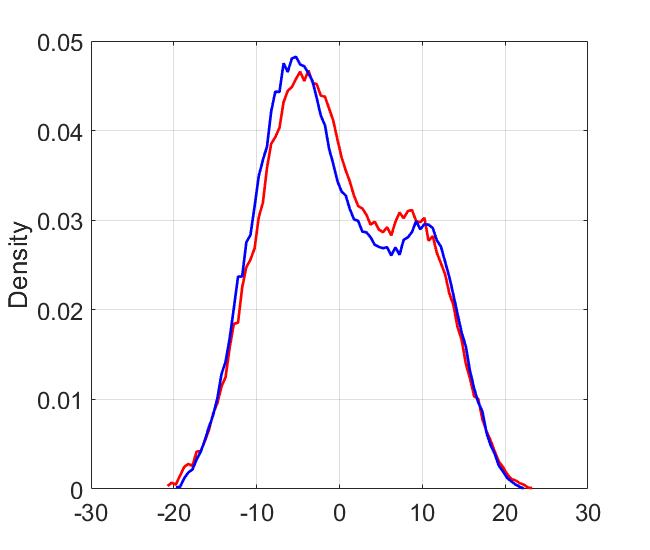}
\caption{$a_2$ histograms}
  \label{figHist2}
\end{subfigure}
\medskip
\begin{subfigure}{0.49\linewidth}
\centering
\hspace*{0cm}\vspace*{0cm}\includegraphics[width=\linewidth]{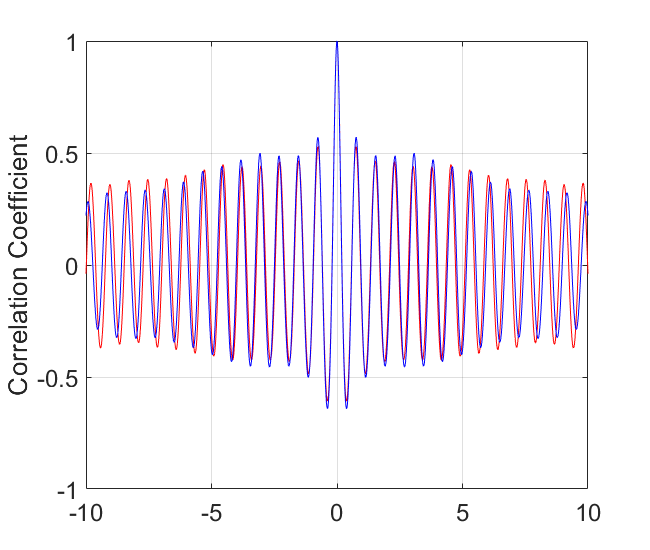}
\caption{$a_2$ autocorrelation functions}
  \label{figAutocorr2}
\end{subfigure}

\begin{subfigure}{0.49\linewidth}
\centering
\hspace*{0cm}\vspace*{0cm}\includegraphics[width=\linewidth]{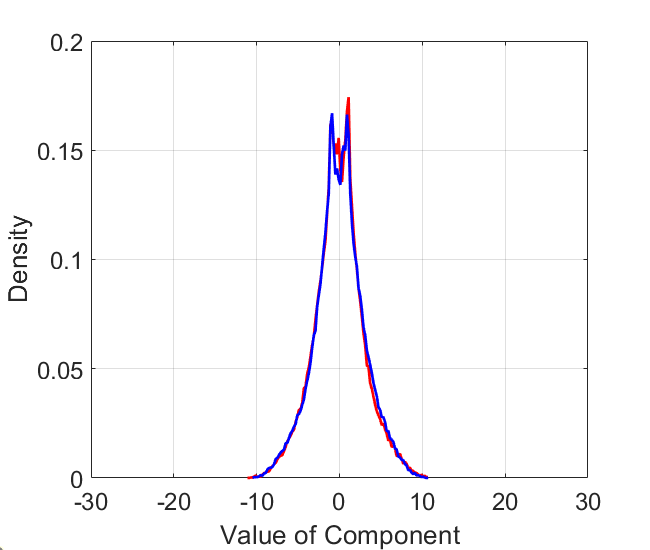}
\caption{$a_3$ histograms}
  \label{Hist3}
\end{subfigure}
\medskip
\begin{subfigure}{0.49\linewidth}
\centering
\hspace*{0cm}\vspace*{0cm}\includegraphics[width=\linewidth]{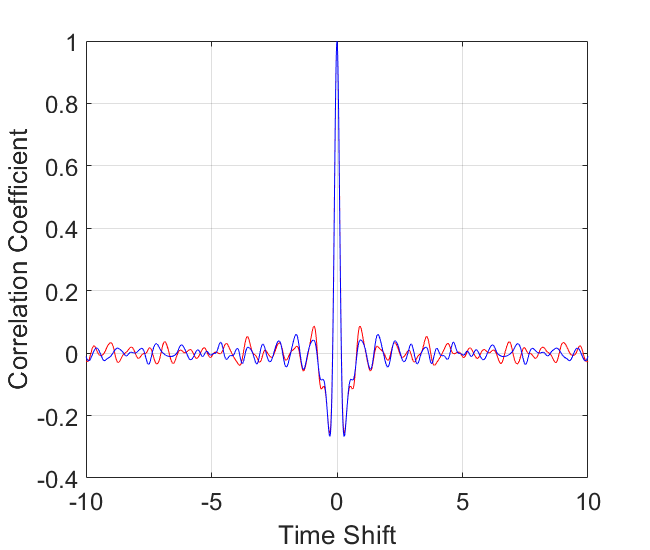}
\caption{$a_3$ autocorrelation functions}
  \label{figAutocorr3}
\end{subfigure}
\caption{Empirical PDFs (a, c, e) and normalized time-autocorrelation functions (b, d, f) of the $a_1$ (a, b), $a_2$ (c, d), and $a_3$ (e, f) components of the L63 system under the true and QMCl dynamics. The PDFs and autocorrelation functions were estimated using time series spanning 1000 model time units.}
\label{figGraphsGrid}

\end{figure}

Next, in \cref{figGraphsGrid}, we examine the ability of QMCl to reproduce some of the one- and two-point statistics under the invariant measure of the L63 system, namely the marginal probability density functions (PDFs) and time-autocorrelation functions of the $a_1$, $a_2$, and $a_3$ coordinates. We estimate these statistics numerically using trajectories generated by the two systems spanning 1000 model time units (i.e., $1000/\Delta t = \text{100,000}$ samples); see Appendix \ref{appHistograms} for further details. It is evident from the results that the QMCl system performs well in terms of reproducing the marginal PDF and time-autocorrelation structure of $a_1$, $a_2$, and $a_3$ with respect to the invariant measure.

Of particular importance in the L63 system is the behavior of the $a_1$ component. Much of the system's chaotic behavior can be understood as chaos in when the system transitions from one lobe of the attractor to another, manifesting as sign changes of the $a_1$ coordinate (see, e.g., \cref{figRealProj2}). The autocorrelation plots in \cref{figAutocorr1} show that the QMCl system accurately captures the initial correlation decay of $a_1$ due to these transitions, which takes place over a timescale comparable to the Lyapunov timescale of the system (approximately 1 model time unit). The corresponding PDF is also reproduced reasonably well (see \cref{figHist1}), though some differences from the true PDF are visible near $a_1 = 0$. These differences are not too surprising given that values of $a_1$ close to 0 correspond to the mixing region between the two lobes of the attractor where the dynamics is particularly sensitive to perturbations. 

Other notable aspects of the results in \cref{figGraphsGrid} are the oscillatory nature of the autocorrelation function for $a_2$ (\cref{figAutocorr2}) and the bimodal corresponding PDF (\cref{figHist2}) due to oscillations about the unstable fixed points in the center of the two ``holes'' of the attractor. The QMCl system is seen to consistently reproduce these features. Meanwhile, the time-autocorrelation of $a_3$ (\cref{figAutocorr3}) exhibits a damped oscillatory behavior that is again reasonably well reproduced by QMCl. 

In calculations not reported here, we observed that replacing the RK4 scheme in the definition of $\tilde \phi$ by a first-order forward Euler scheme results in a noticeable reduction of statistical accuracy of the QMCl model in terms of the PDFs and autocorrelation functions of $a_1$, $a_2$, and $a_3$. This effect is independent of QMCl since integrating the L63 system \eqref{eqL63Palmer} with a forward Euler scheme of fixed timestep $\Delta t = 0.01$ (as opposed to \texttt{ode45}) was found to impart similar changes to the $a_1, a_2, a_3$ statistics.  

In summary, the results presented in this subsection demonstrate that the QMCl-parameterized system provides an accurate surrogate model of the full L63 dynamics.   

\subsection{\label{secL63Small}Small, partial training data} 

The experiment in \cref{secL63Big} assumed access to a large training dataset containing full information about the system state. Here, we show that with the same parameters but much less training data, $N = \text{10,000}$, the QMCl method still performs a reasonable reconstruction of the L63 dynamics. On this experiment, we also impose the constraint that in generating the kernel eigenvectors defining $H_{L,N}$ we do not have access to full state vectors $(a_1,a_2,a_3)$. Instead, only the $a_1$ component is available. 

Let $\tilde W : \Omega \to \mathbb R$ be the map that projects onto the first coordinate, i.e., $\tilde W(\omega) = a_1$ with $\omega = (a_1, a_2, a_3)$. Let also $\tilde w_m = \tilde W(\omega_m)$, where $\omega_0, \omega_1, \ldots, \omega_{N-1} \in \Omega$ is the training trajectory in state space (defined as in \cref{secL63Big}, but assumed here unobserved). To enrich the data $\tilde w_m$ with information lost due to projection by $\tilde W$, delay-coordinate embedding is used on this one-dimensional time series. Setting $ \mathcal W = \mathbb R^Q$, where $Q \in \mathbb N$ is an even parameter corresponding to the number of delays, we define the delay-coordinate map $W : \Omega \to \mathcal W$ such that 
\begin{equation}
    \label{eqWDelay}
    W(\omega) = (\tilde W(\Phi^{-Q/2}(\omega)), \tilde W(\Phi^{-Q/2 +1}(\omega)), \ldots, \tilde W(\Phi^{Q/2}(\omega)).
\end{equation}
From the theory of delay-coordinate maps \cite{SauerEtAl91} it is known that for sufficiently large $Q$, and with ``high probability'' in a suitable sense, $W$ is an injective map on compact subsets of $\Omega$, and thus the support of the invariant measure $\mu$ (which is compact since it is contained in an absorbing ball under the L63 dynamics \cite{LawEtAl14}). Thus, for sufficiently large $Q$, training data obtained through the map $W$ should be theoretically sufficient to build a basis for the entire Hilbert space $H$. Importantly for practical applications, we have
\begin{displaymath}
    W(\omega_m) = ( \tilde w_{m-Q/2}, \tilde w_{m-Q/2+1}, \ldots, \tilde w_{m+Q/2} ),
\end{displaymath}
which means that we can evaluate $W$ on the dynamical states $\omega_m$ underlying the training data without knowledge of these states. In particular, we can compute kernel matrices and build an associated data-driven basis of $H_{L,N} $ as described in \cref{secDataDrivenBasis} using data sampled from $W$. It is worthwhile noting that as the number of delays $Q$ increases, the kernel eigenfunctions obtained from delay-coordinate-mapped data tend to span approximately Koopman-invariant subspaces \cite{DasGiannakis19,Giannakis21a}, which improves the quality of Koopman/transfer operator approximation in these subspaces. This fact motivates using delay-coordinate maps even when full system states are available for training.      

\begin{figure}
    \begin{subfigure}{.49\linewidth}
\centering
\hspace*{0cm}\vspace*{0cm}\includegraphics[width=\linewidth, scale=.75]{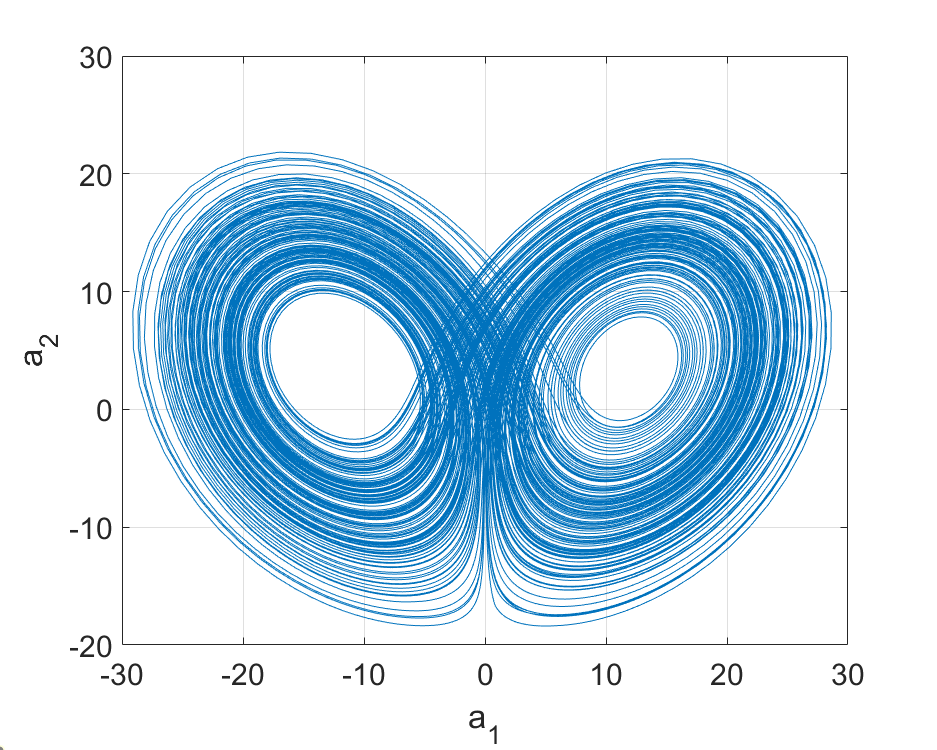}
\caption{$(a_1, a_2)$ projection}
  \label{figLowDataProjA}
\end{subfigure}
\medskip
\begin{subfigure}{.49\linewidth}
  \centering
\hspace*{0cm}\vspace*{0cm}\includegraphics[width=\linewidth, scale=.75]{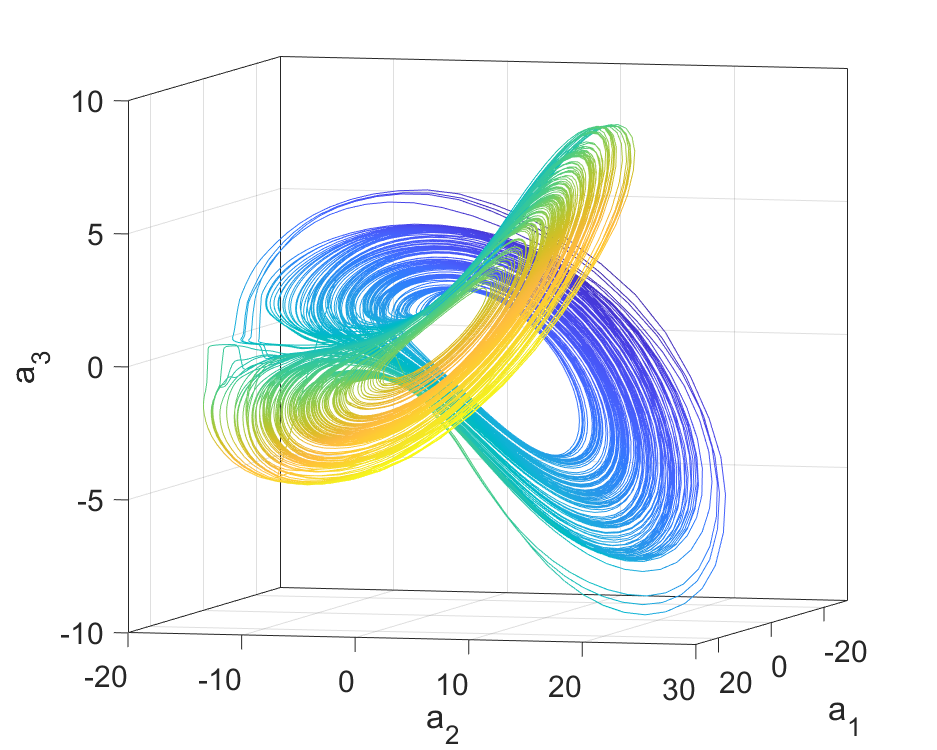}
\caption{$(a_1,a_2, a_3)$ state space}
  \label{figLowDataProjB}
\end{subfigure}

\begin{subfigure}{.49\linewidth}
  \centering
\hspace*{0cm}\vspace*{0cm}\includegraphics[width=\linewidth, scale=.75]{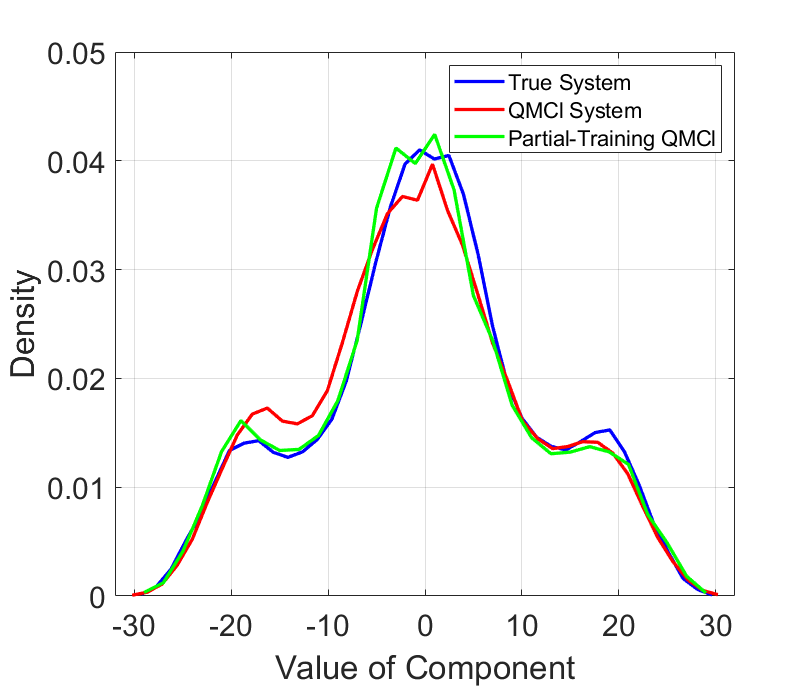}
\caption{$a_1$ histogram}
  \label{figLowDataProjD}
\end{subfigure}
\begin{subfigure}{.49\linewidth}
\centering
\hspace*{0cm}\vspace*{0cm}\includegraphics[width=\linewidth, scale=.75]{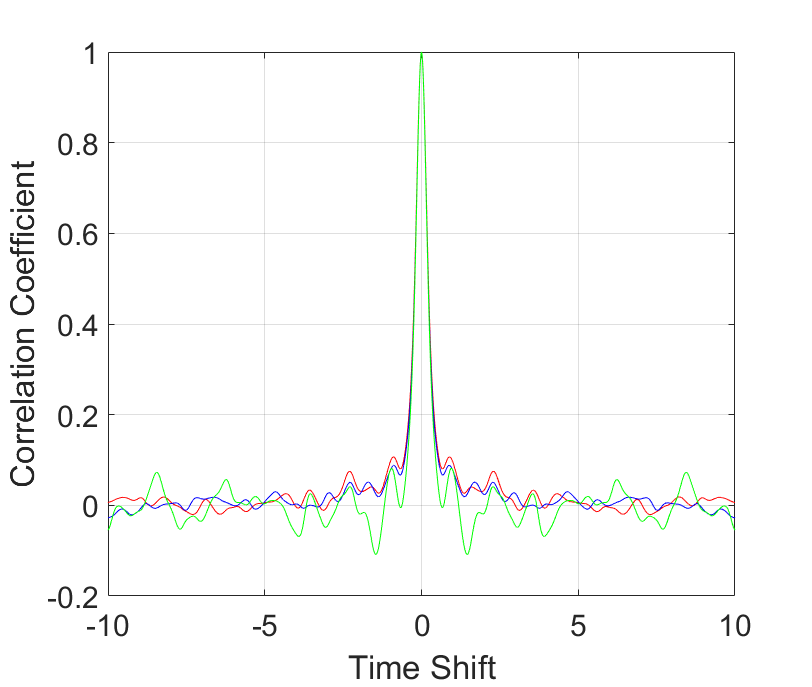}
\caption{$a_1$ autocorrelation}
  \label{figLowDataProjC}
\end{subfigure}
\caption{Trajectories (a, b), marginal PDF of $a_1$ (c), and normalized time-autocorrelation function of $a_1$ (d) for the L63 QMCl system generated from 10,000 training points. The basis functions for this QMCl system were generated only from the $a_1$ L63 component of the training data, using delay-coordinate embedding. The trajectories in (a, b) span 200 model time units. The PDFs and autocorrelation functions for the partial-training system in (c, d) were estimated using time series spanning 200 model time units. PDFs and autocorrelation functions from the true L63 system and QMCl system from \cref{figHist1,figAutocorr1} are shown for reference.}
\label{figLowData}
\end{figure}

\cref{figLowData} displays representative trajectory, marginal PDF, and autocorrelation results obtained using $N=\text{10,000}$ training samples from~\eqref{eqWDelay} and $Q = 10$ delays. For the kernel $\kappa_{\mathcal W}$ used to build the basis we used the variable-bandwidth kernel proposed in \cite{BerryHarlim16}, normalized to a symmetric Markov kernel using the approach of \cite{CoifmanHirn13}. A description of the construction of this kernel can be found in \cite{FreemanEtAl22}. For the feature-map kernel $k$, a bandwidth value $\epsilon = 10$ was chosen. One can notice that this system still closely approximates both the marginal PDFs (\cref{figLowDataProjD}) and autocorrelation functions (\cref{figLowDataProjC}) of the true system under the invariant measure. However, in both cases, it is evident that the larger, more informative training data used in the experiment of \cref{secL63Big} allows for a more accurate approximation of the true system; see, e.g., the autocorrelation plots in \cref{figLowDataProjC}. 

\subsection{\label{secL63Stoch}Stochastic closure experiments}

\begin{figure}
\begin{subfigure}{0.49\linewidth}
\centering
\hspace*{0cm}\vspace*{0cm}\includegraphics[width=\linewidth]{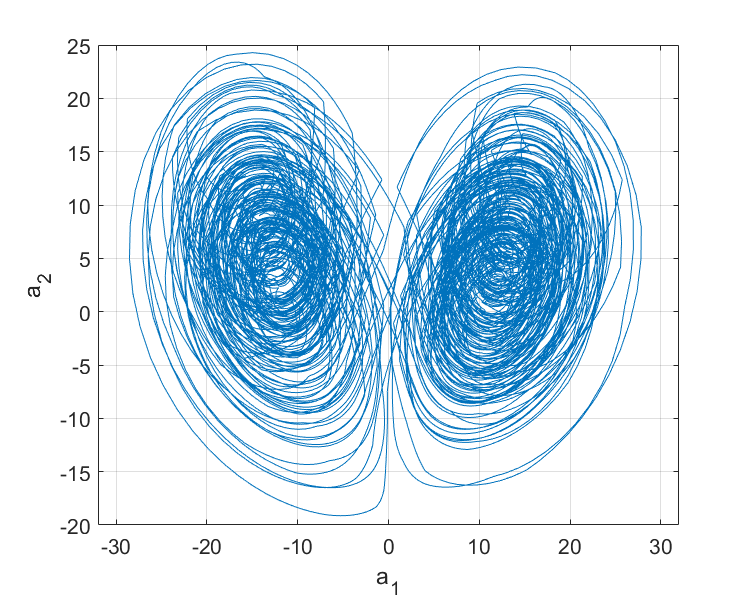}
\caption{$(a_1, a_2)$ projection}
\label{figPalmer2}
\end{subfigure}
\begin{subfigure}{0.49\linewidth}
\centering
\hspace*{0cm}\vspace*{0cm}\includegraphics[width=\linewidth]{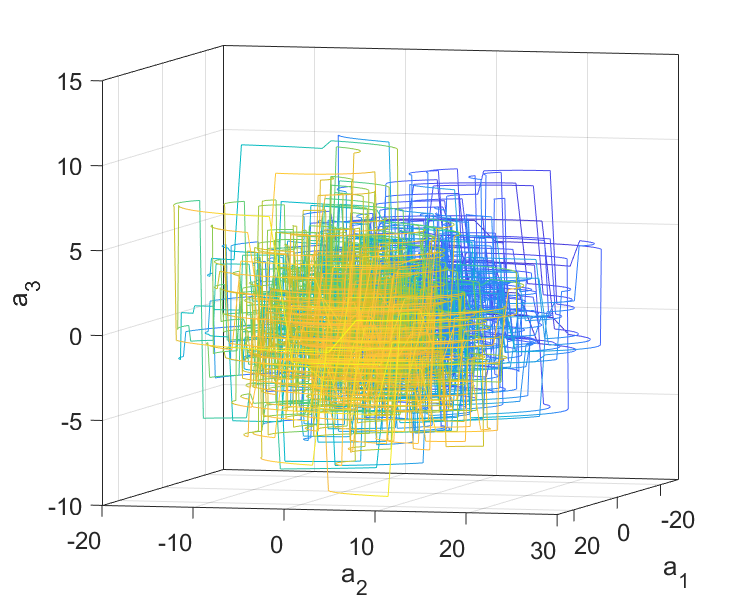}
\caption{$(a_1,a_2,a_3)$ state space}
\label{figPalmer3}
\end{subfigure}
\caption{Projected (a) and three-dimensional (b) trajectories generated by the Gaussian parameterization of the L63 system from \cref{secStochasticExamples}. The plotted trajectories span 200 model time units. In (b), the color corresponds to the value of the $a_1$ component. Notice the lack of structure in the $a_3$ coordinate compared to \cref{figRealProj1} due to modeling of the $a_3$ dynamics as an i.i.d.\ Gaussian process.}
    \label{figPalmer}
\end{figure}

\begin{figure}
\begin{subfigure}{0.49\linewidth}
\centering
\hspace*{0cm}\vspace*{0cm}\includegraphics[width=\linewidth]{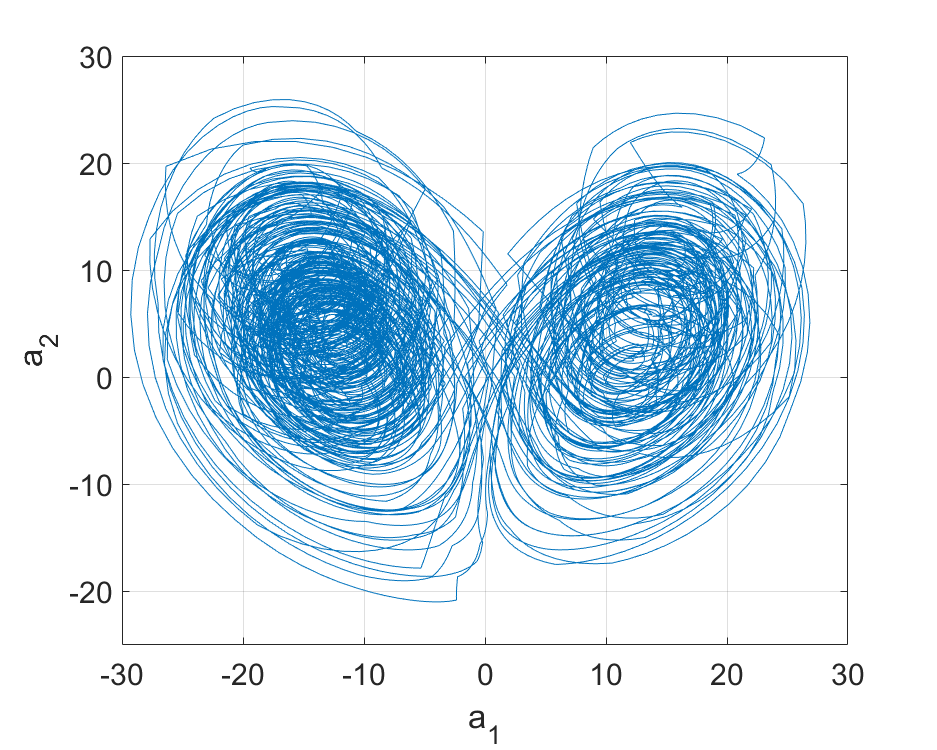}
\caption{$(a_1,a_2)$ projection}
\label{figStochasticTraj}
\end{subfigure}
\medskip
\begin{subfigure}{0.49\linewidth}
  \centering
\hspace*{0cm}\vspace*{0cm}\includegraphics[width=\linewidth]{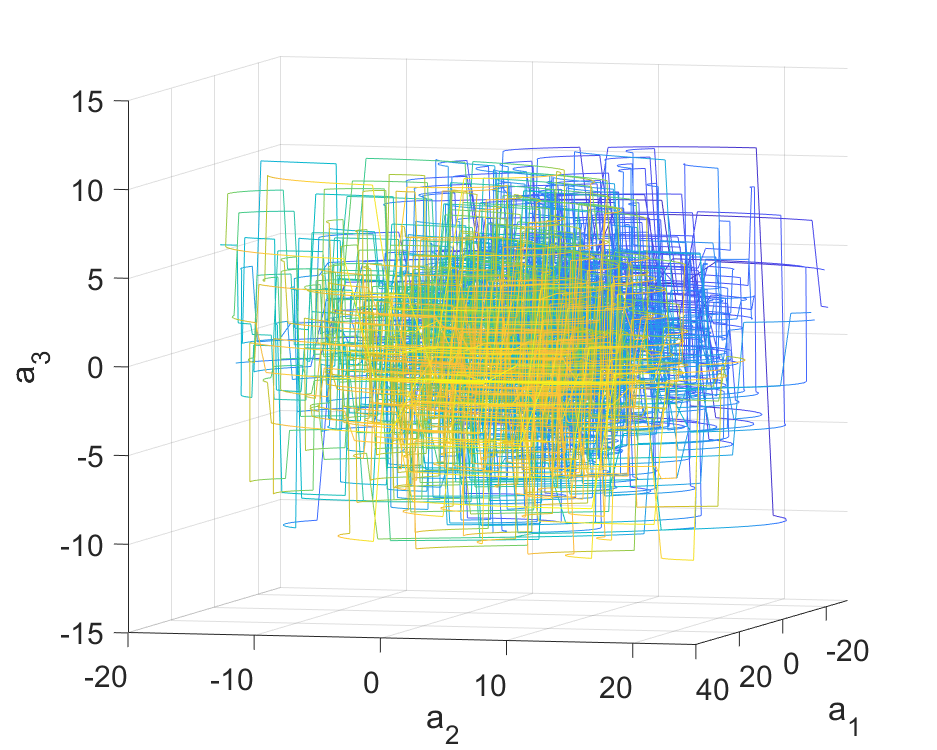}
\caption{$(a_1,a_2,a_3)$ state space}
\label{figStochastic3D}
\end{subfigure}
\begin{subfigure}{0.49\linewidth}
\centering
\hspace*{0cm}\vspace*{0cm}\includegraphics[width=\linewidth]{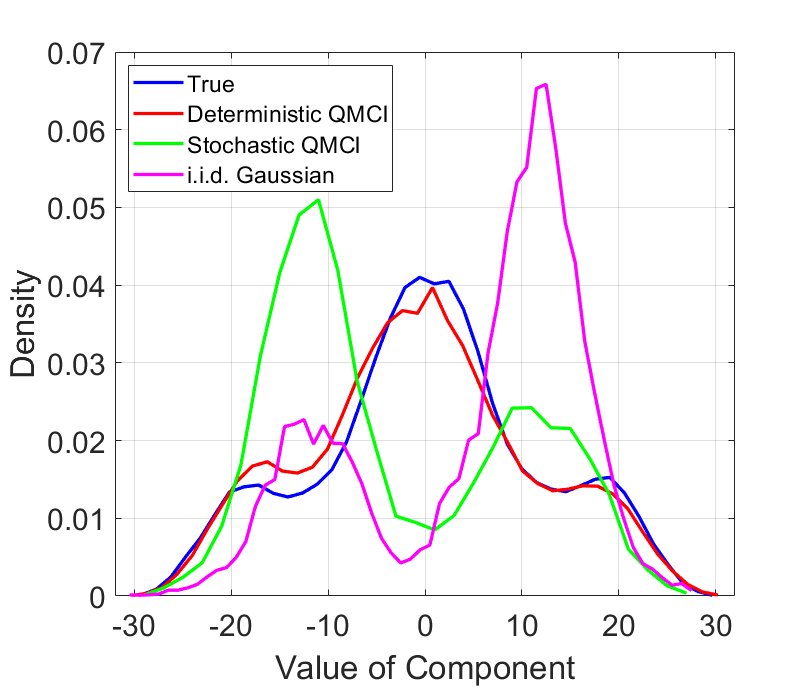}
\caption{$a_1$ histograms}
\label{figStochasticHist}
\end{subfigure}
\begin{subfigure}{0.49\linewidth}
  \centering
\hspace*{0cm}\vspace*{0cm}\includegraphics[width=\linewidth]{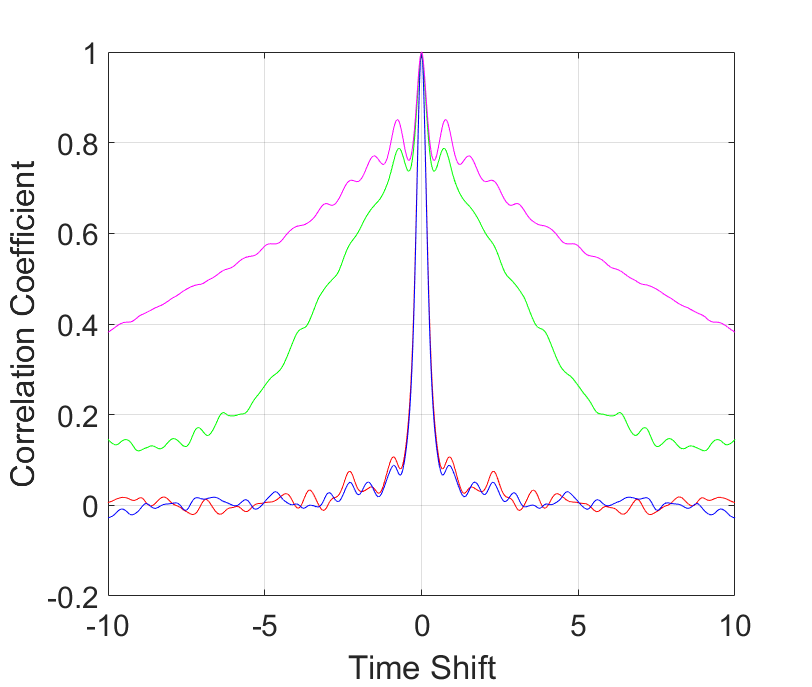}
\caption{$a_1$ autocorrelation functions}
\label{figStochasticAutocorr}
\end{subfigure}

\caption{Trajectories (a, b), marginal PDF of $a_1$ (c), and time-autocorrelation of $a_1$ (d) for the stochastic QMCl closure of the L63 system. The trajectories in (a, b) span 200 model time units. The PDFs and autocorrelation functions for the stochastic QMCl system in (c, d) were estimated using time series spanning 200 time units. PDFs and autocorrelation functions from the true L63 system, the deterministic QMCl system from \cref{figHist1,figAutocorr1}, and the i.i.d.\ Gaussian closure \cite{Palmer01} (see \cref{figPalmer} for trajectories) are also shown.}
\label{figStoch}

\end{figure}

For purposes of establishing a baseline of performance based on a concise and fast stochastic method, we first examine plots generated by the i.i.d.\ Gaussian closure of \cite{Palmer01} (see \cref{secStochasticExamples}). \cref{figPalmer} shows representative trajectories for the $(a_1,a_2)$ and $(a_1,a_2,a_3)$ components generated by this system. Though both attractor lobes are visible, there is noticeably less accuracy in the behavior when it comes to transitions between lobes---this is reflected in the autocorrelation plots shown in \cref{figStochasticAutocorr}. The stochastic closure also fails to reproduce the characteristic ``holes'' around the fixed points in the lobes of the attractor. This is reflected in the $a_1$ PDF from the Gaussian closure in \cref{figStochasticHist} which exhibits high probability density for values of $a_1$ that have lower probability density under the true system due to the holes. 

Viewing \cref{figPalmer3} allows us a glimpse into what is going on ``under the hood'' of the stochastic Gaussian closure. Since the choice of the $a_3$ component is Gaussian i.i.d.\ for each timestep, the $a_3$ component appears simply as random noise, while the useful approximation of the attractor manifests only in the $(a_1, a_2)$ projection (see \cref{figPalmer2}). On the other hand, in \cref{figQMProj1} we can see that the deterministic QMCl approach reconstructs the entire attractor, even including the parameterized dimension. In some applications such as climate dynamics, the increase in computation cost and requisite information in QMCl (see \cref{secCost}) may result in parameterizations based on cheaper parametric stochastic models such as the i.i.d.\ Gaussian closure of the L63 system being favorable. However, the ability to accurately reconstruct entire attractors may allow for the reproduction of more complicated and subtle dynamical properties in the resolved variables, which could in turn have value in areas where additional computation time would be worth investing. 

In comparison with the deterministic QMCl and i.i.d.\ Gaussian closure, we also examine the stochastic QMCl closure described in \cref{secStochasticQMCL}. In this experiment, we use $N = \text{90,000}$ samples of the state vector (i.e., $W=\Id$ as in \cref{secL63Big}) and $L = 1200$ eigenfunctions to build the QMCl Hilbert space $H_{L,N}$. The quantum state was updated through the observation kernel after every ten timesteps (i.e., $r=10$). The bandwidth parameter for the kernel $\kappa_{\mathcal W}$ used to build the basis was algorithmically chosen to be $\epsilon_{\mathcal W} = \sqrt{19}$, and the bandwidth parameter for the feature map kernel $k$ was chosen to be $\epsilon = 35$. 

It is worth noting that a relatively large value of $\epsilon$ was chosen for this experiment (cf.~$\epsilon=2$ and $\epsilon=10$ in \cref{secL63Big,secL63Small}, respectively) due to numerical stability issues. Namely, for values of $\epsilon$ comparable to the deterministic QMCl experiments, the stochastic QMCl system tends to run into numerical errors associated with the effect-update step if the current point lands too far from the training data (meaning that, for some $\hat x_n \in \mathcal X$ under the parameterized dynamics, $k(\hat x_n,x_m)$ is numerically zero for all $x_m$ in the training dataset). Increasing $\epsilon$ makes the system more stable with respect to this particular problem, at the expense of a less informative kernel. The deterministic QMCl system is less noisy, and as it appears inherently less likely to run into this issue, allowing for the value of $\epsilon$ to be chosen significantly smaller. 

In \cref{figStochasticHist,figStochasticAutocorr}, we compare the marginal PDF and time-autocorrelation function of the $a_1$ variable, respectively, of the true system, the deterministic and stochastic QMCl systems, and the i.i.d.\ Gaussian closure. We can notice that the stochastic QMCl closure comes closer to reconstructing the correct PDF and autocorrelation function than the i.i.d.\ Gaussian closure. As is clear from \cref{figStochasticTraj}, the stochastic QMCl algorithm recovers some qualitative aspects of the L63 system, but fails in other important respects. In \cref{figStochastic3D} we can see that the three-dimensional $(a_1,a_2,a_3)$ trajectories are somewhat more coherent than their counterparts from the Gaussian closure in \cref{figPalmer3} (since the stochastic flux terms in QMCl are not i.i.d., and their distribution depends on the resolved variables), but they are still dominated by noise. Moreover, similarly to the Gaussian closure, the ``holes'' of the attractor lobes (\cref{figStochasticTraj}) are not adequately recovered. In \cref{figStochasticHist}, we notice that the $a_1$ histograms have significantly different structures---namely, the trajectories generated from the stochastic QMCl system tend to cluster in the two lobes, while the true system has a high density of points in the central transition region. Evidently, while the ``double wing'' nature of the attractor is present in the stochastic QMCl system, the smoothness of the transition between lobes, the general distribution of points, and the rate and qualitative behavior of the system as it transitions are all less accurate to the true system than the deterministic QMCl system is.

\section{\label{secL96}Lorenz 96 multiscale}

The L96 multiscale model \cite{Lorenz96,FatkullinVandenEijnden04} is a system of equations of $K$ variables $\{ x_k \}_{k=1}^K$ and $JK$ variables $\{y_{j,k}\}_{j,k=1}^{J,K}$ defined by
\begin{equation}
    \label{eqL96}
    \begin{gathered}
        \begin{aligned}
            \dot x_k &= -x_{k-1} \left( x_{k-2} - x_{k+1} \right) - x_k + F - h_x \bar y_k, \\  
            \dot y_{j,k} & =\frac{1}{\varepsilon} \left( - y_{j+1,k} \left( y_{j+2,k} - y_{j-1,k} \right) - y_{j,k} + h_y x_k \right),
        \end{aligned}\\
        x_{k+K} = x_k, \quad y_{j,k+K}=y_{j,k}, \quad y_{j+J,k} = y_{j, k+1},
    \end{gathered}
\end{equation}
where $\bar y_k = \sum_j y_{j,k} / J $, and $F$, $h_x$, $h_y$, and $\varepsilon$ are real parameters. The L96 multiscale system is a more comprehensive model of atmospheric dynamics than the L63 system. For $\varepsilon \ll 1$, it is a multiscale system, in which the variables $x_k$ vary slowly in time and each have an associated set $\{y_{j,k}\}_{j=1}^J$ of variables which vary quickly. Each slow variable $x_k$ is only influenced by the fast variables via the average value $\bar y_k$ of its associated fast variable set. The L96 multiscale system has been extensively used as a testbed for parameterization schemes; see, e.g., \cite{Wilks05,CrommelinVanden-Eijnden08,ArnoldEtAl13,JiangHarlim20b,ChenLi21,GottwaldReich21} among many references. The objective of parameterization in this case is to approximate the behavior of the slow variables. That is, we have $x = (x_1, \ldots, x_K) \in \mathcal X \equiv \mathbb R^K$, $y = ( y_{1,1}, \ldots, y_{J,K} ) \in \mathcal Y \equiv \mathbb R^{JK}$, and $Z: \mathcal Y \to \mathcal Z = \mathbb R^K$ with $Z(y) = (\bar y_1, \ldots, \bar y_K)$ for the resolved variables, unresolved variables, and flux terms, respectively. 

As $\varepsilon \to 0$, the L96 multiscale system is known to exhibit an averaging limit \cite{PavliotisStuart08}, in which the evolution of the $x_k$ variables is Markovian and is governed by the system of equations
\begin{displaymath}
    \dot x_k = - x_{k-1} \left( x_{k-2} - x_{k+1} \right) -x_k + h_x \hat Z_k(x), 
\end{displaymath}
for some functions $\hat Z_1, \ldots, \hat Z_K : \mathcal X \to \mathbb R$. Following refs.~\cite{BurovEtAl21,FreemanEtAl22}, we choose the parameters $K = 9$, $J = 8$, $h_x = - 0.8$, $h_y = 1$, and $\varepsilon = 1/128$. The resulting dynamical regime is chaotic, with approximately Markovian dynamics for the $x_k$ variables.   

\subsection{Quantum mechanical closure experiments}

Let $\Phi^t : \Omega \to \Omega $ with $\Omega = \mathcal X \times \mathcal Y$ and $t \in \mathbb R$ be the flow generated by~\eqref{eqL96}. Similarly to the L63 experiments in \cref{secL63}, we consider a discrete-time system $ \Phi : \Omega \to \Omega $ with $\Phi = \Phi^{\Delta t} $ obtained by temporal subsampling of the flow. Here, the timestep is $\Delta t = 0.01$ model time units. As training data for QMCl, we use time series $w_0, \ldots, w_{N-1} \in \mathcal W \equiv \mathcal X \equiv \mathbb R^K$ and $ z_0, \ldots, z_{N-1} \in \mathcal Z$, where $ w_m = \mathcal P_{\mathcal X}(\omega_m)$, $z_m = Z(\omega_m)$, and $\omega_0, \ldots, \omega_{N-1} \in \Omega$ with $ \omega_m = \Phi^m(\omega_0)$. In particular, we build the basis of $H_{L,N}$ using information from only the slow variables. The numerical trajectory $\omega_m$ is generated using MATLAB's \texttt{ode15s} solver which is appropriate for stiff problems. The number of training samples is $N = \text{40,000}$, and the initial condition $\omega_0$ is taken on the trajectory starting from $(x_1, \ldots, x_K) = ( 1, 0, \ldots, 0)$ and $(y_{1,1}, \ldots, y_{J,K}) = (1.1, 0,\ldots,0)$ after an equilibration time interval of 200 model time units (i.e., similarly to the L63 experiments; see \cref{secL63Big}). 

We build the kernel eigenbasis of the Hilbert space $H_{L,N}$ with dimension $L=1900$ using the radial Gaussian kernel $\kappa_{\mathcal W}$ from~\eqref{eqKW}. As in the L63 experiments from \cref{secL63}, the bandwidth parameter ${\epsilon_{\mathcal W} = \sqrt{10}}$ was tuned automatically. Moreover, the evolution map $\tilde \phi : \mathcal X \times \mathcal Z \to \mathcal X $ is based on the RK4 scheme in~\eqref{eqRK4}. We evolve the quantum states using the map $\psi_r : \mathcal X \times Q(H_{L,N}) \to Q(H_{L,N})$ with $ r = 5 $ Koopman evolution steps (i.e., $ r \, \Delta t = 0.05$ model time units) between each state conditioning via~\eqref{eqQMPosterior}. The effect-valued feature-map $\mathcal F_{L,N}$ was based on a radial Gaussian kernel, here with bandwidth parameter $\epsilon = {2}$.  Similarly to the L63 experiments, we assess the performance of QMCl in terms of its ability to reproduce salient qualitative features of the true dynamics, as well as marginal PDFs and time-autocorrelation functions of the revolved variables. In these tests, the true and QMCl systems are initialized on the trajectory starting from $(x_1, \ldots, x_K) = (1,0,\ldots, 0) $ and $(y_{1,1}, \ldots, y_{J,K}) = (1, 0, \ldots, 0)$ after an equilibration period of 1000 model time units. We use {100,000 samples (i.e., $\text{100,000}\,\Delta t = 1000$ model time units)} on the true and QMCl trajectories starting at that point to compute PDFs and autocorrelation functions. The initial QMCl state is obtained via~\eqref{eqInit2}.

\begin{figure}
    \begin{subfigure}{\linewidth}
        \centering
        \hspace*{0cm}\vspace*{0cm}\includegraphics[width=\linewidth]{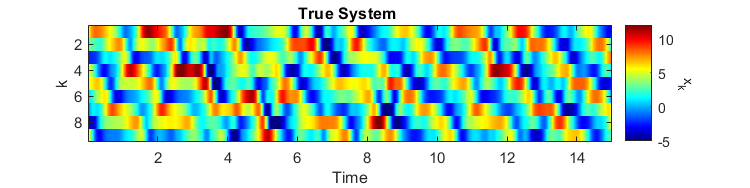}
    \end{subfigure}
    \begin{subfigure}{\linewidth}
        \centering
        \hspace*{0cm}\vspace*{0cm}\includegraphics[width=\linewidth]{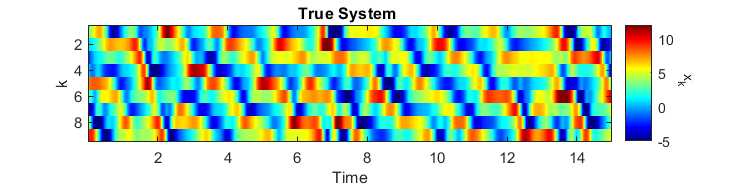}
    \end{subfigure}
    \caption{Hovmoller diagrams of the slow variables $x_k$ for the true L96 multiscale system (top) and the QMCl-parameterized system (bottom).}
    \label{figHeatMaps}
\end{figure}

\cref{figHeatMaps} displays Hovmoller diagrams (space--time heat maps) of the $x_k$ variables under the true L96 and QMCl dynamics. In our chosen dynamical regime, the $x_k$ variables exhibit characteristic propagating patterns which can be thought of as crude representations of eastward-propagating disturbances in the Earth's midlatitude atmosphere. The dynamics of these wave-like structures, including their aperiodic emergence and decay, appear to be qualitatively well-captured by the QMCl system.

\begin{figure}
\begin{subfigure}{0.49\linewidth}
  \centering
\hspace*{0cm}\vspace*{0cm}\includegraphics[width=\linewidth]{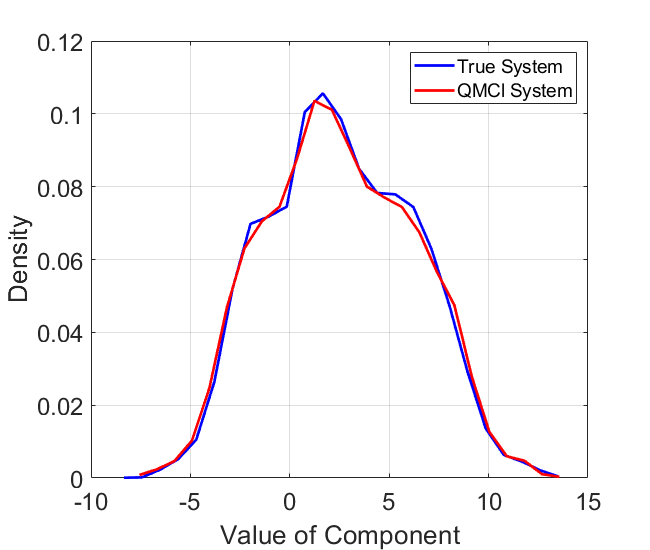}
\caption{Histograms}
\label{figL96Hist}
\end{subfigure}
\begin{subfigure}{0.49\linewidth}
\centering
\hspace*{0cm}\vspace*{0cm}\includegraphics[width=\linewidth]{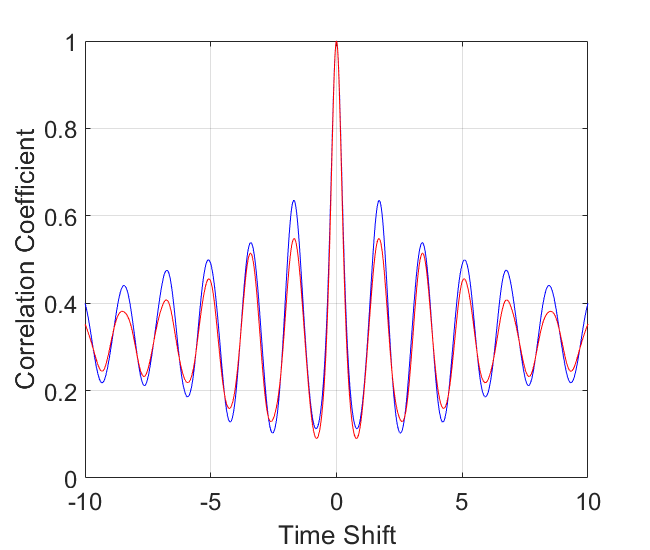}
\caption{Autocorrelation Plots}
\label{figL96Corr}
\end{subfigure}

\caption{Marginal PDFs (a) and time-autocorrelation functions (b) of the $x_1$ component of the L96 multiscale system under the true and QMCl dynamics,  estimated using time series spanning 1000 model time units.}
\label{figStochL96}

\label{figL96}
\end{figure}

As a more quantitative test, in \cref{figL96} we compare the marginal PDFs and time-autocorrelation functions of $x_1$ under the true L96 and QMCl dynamics. We can see in \cref{figL96Corr} that the autocorrelation under the QMCl system matches the periodicity present in that of the true system, though the decay of the correlation amplitude is somewhat faster in the QCMl system. The $x_1$ PDFs (\cref{figL96Hist}) are also in good agreement between the true and QMCl systems. 

\section{Summary and discussion}\label{secDiscussion}

We have developed a data-driven framework for closure of dynamical systems that combines aspects of quantum theory, ergodic theory, and kernel methods for machine learning. Our approach, called Quantum Mechanical Closure (QMCl), models the unresolved degrees of freedom as a finite-dimensional quantum mechanical system that is coupled to a classical model governing the resolved variables. The state of the quantum system is a density operator (a quantum mechanical analog of a probability density) that evolves dynamically under the induced action of the Koopman operator. The fluxes from the unresolved variables to the resolved variables are modeled as expectations of quantum observables (self-adjoint operators) that provide the closure terms needed to advance the resolved variables given the current quantum state. Meanwhile, the state of the resolved variables updates the quantum state by conditioning, using an operator-valued feature map in a step that can be viewed as a quantum mechanical analog of Bayes' rule. The resulting two-way coupling between the classical and quantum systems for the resolved and unresolved degrees of freedom, respectively, resembles the prediction--correction cycle of sequential data assimilation (filtering). 

QMCl has a data-driven implementation in which a dataset of $N$ samples of the resolved variables is used to learn a basis for an $L$-dimensional Hilbert space $H_{L,N}$ consisting of eigenfunctions $\phi_{l,N}$ of a kernel integral operator. The quantum mechanical closure model is then built over $H_{L,N}$, and all operators employed in the scheme are numerically represented as $L\times L$ matrices with respect to the $\{\phi_{l,N}\}$ basis. The data-drive formulation of QMCl has a well-characterized large-data limit, $N\to\infty$, leveraging results from spectral approximation of kernel integral operators. In addition the method has a stochastic variant wherein the closure terms are obtained by random draws at each timestep from the probability distribution of quantum mechanical measurements induced by the quantum state. 

We have applied QMCl to deterministic and stochastic closure experiments involving the L63 (\cref{secL63}) and L96 multiscale (\cref{secL96}) systems; the latter, in a chaotic regime with timescale separation. In these experiments, QMCl was able to reconstruct the qualitative nature of the resolved variables, as well as their marginal distribution and time-autocorrelation functions with respect to the invariant measure. In the L63 examples, the unresolved variable was also recovered, resulting in a three-dimensional reconstruction of the Lorenz attractor.   

We end the paper with a discussion of some of the salient features of QMCl and possible avenues for future work.  

\subsection{Positivity preservation}

A novel aspect of QMCl is that it casts the problem of closure of dynamical systems in the setting of operator algebras and the associated quantum probability theory. What we have argued in this paper is that the structure of these spaces naturally leads to computational algorithms with improved structure-preservation properties compared to formulations in abelian algebras of classical observables. In particular, as a result of basic properties of projections on operator algebras, QMCl represents positive classical observables by positive operators, whereas orthogonal projections onto finite-dimensional function spaces, are not, in general, positivity preserving. By virtue of the same properties, QMCl is able to represent classical probability densities by ``honest-to-god'' density operators in finite dimensions, whereas finite orthogonal basis expansions of probability densities are not, in general, positive functions (and thus not normalizable to probability densities). 

We believe that the positivity preservation enjoyed by QMCl is a useful property in applications. In geophysical fluid dynamics, for instance, many relevant physical quantities such as density, moisture, and pressure, by definition take non-negative values. Since negative values of such variables have no physical meaning, a parameterization scheme which generates negative-valued outputs of these observables could have particularly problematic effects on the system as a whole. Core conservation laws, regarding mass, for example, could be violated by such a negative quantity, leading to poor, or even meaningless, predictions as time advances. In the context of parameterizations of climate models, some of the biases and numerical instabilities exhibited by data-driven parameterization schemes have indeed been attributed to failure to respect physical constraints such as energy conservation and positivity of precipitation \cite{YuvalOGorman20}, and it is possible that QMCl may provide a route to addressing such issues.     

\subsection{Stochastic versus deterministic approaches}

In areas such as climate modeling, stochastic parameterization as a general approach has been increasingly studied and used to successful ends \cite{BernerEtAl17}. The stochastic mode of the QMCl framework (\cref{secStochasticQMCL}) can be thought of as a new stochastic parameterization algorithm. In this view, at each time, the probability distribution associated with the stochastically-drawn flux term $z$ is determined by a generalized probability density function which is the quantum density operator $\rho$. Though the framework of QMCl fundamentally alters the mathematical structure underlying the probabilistic dynamics of the distribution, as a whole this can nonetheless can be thought of in the same general terms as stochastic parameterizations based on classical probability theory. 

The deterministic approach (\cref{secFlux,secStateUpdate,secFeatureMap,secInit}), instead of drawing the value of $z$ randomly from an induced probability measure, directly chooses $z$ to be the expectation value of a quantum observable given the quantum state $\rho$. This is a natural choice as an estimator of the true flux given the unresolved variables, and is similar in spirit to statistical closure methods that estimate fluxes through classical expectations \cite{JiangHarlim20b}, including methods based on ensemble data assimilation \cite{ChenLi21,GottwaldReich21,ChenEtAl22,LevineStuart22}. We argue that the deterministic mode of QMCl is better thought about in a fundamentally different way than the stochastic mode (that is, differently than as some extension or limiting case of the stochastic approach). 

Upon first glance, one might imagine that dropping stochasticity counteracts the primary benefit of this method as an approach compared to more standard approaches to parameterization. Indeed, moving away from deterministic functions of the state as parameterizations, due to their inherent dynamics-flattening properties, is one reason why we were interested in such techniques to begin with. However, while the choice of $z$ at a given time is deterministic in this mode, it is not a function of the classical state at a given time---it is a function of the quantum state. That is, instead of our deterministic parameterization function being a function on the space of unresolved variables $\mathcal Y$, it is now a function on the space of $L \times L$ density matrices. The latter, is a convex subset of the set of $L \times L $ Hermitian matrices of (real) dimension $L^2 - 1$ \cite{BengtssonZyckowski06}. Assuming that the initial state is pure (i.e., a rank-1 projection, as is the case in the experiments in \cref{secL63,secL96}), all subsequent states produced by the QMCl algorithm are pure, and thus lie in a subset of the set of all quantum states of dimension $2(L-1)$.  Since $L$ can very large compared to the dimension of $\mathcal Y$, our parameterization is no longer a dimension reduction, but in fact a massive dimension expansion. These additional dimensions allow for much more information to be carried through time, and allow us to avoid the dynamics-flattening that comes from a reduction approach. 

A number of benefits arise from the deterministic variant of QMCl. For one, it is significantly computationally cheaper than the stochastic version, since the step of generating a probability distribution and drawing points from it is replaced by the comparatively cheaper calculation of computing $\tr(\rho (\pi_L \zeta^{(i)} ))$ from~\eqref{eqQMFlux} (particularly when $\rho$ is a pure state, in which case the computational cost is $O(L)$). Furthermore, in systems such as L96 multiscale (\cref{secL96}), where multiple draws need to be taken at each timestep, we have a simple solution to the problem of choosing how to correlate the approximations of the various unresolved components. In the stochastic setup, all parameters were drawn separately (from the same quantum state, but nonetheless independently). Potentially important relations between what the values of the various parameters could be at a given time were thus lost---regardless of the accuracy of the PDFs, individual variation could result in sets of parameters which could not actually exist in reality. In the deterministic approach, this problem is implicitly resolved. All flux terms are taken to be the expectation of their observable $\pi_L \zeta^{(i)}$ associated with a given quantum state, so there is an implicit correlation between the parameters which naturally arises from the quantum state itself. Finally, it is not unreasonable to think that analyzing convergence properties and error bounds of the deterministic approach may be easier in future research than a stochastic mode. 

In the L63 results of \cref{secL63Stoch}, we can already see clear benefits of using the expectation value over a random draw. What is manifest in the trajectory plots in \cref{fig3DProjections}, along with the marginal PDF and time-autocorrelation plots in \cref{figGraphsGrid}, is that the quantum state appears to be carrying sufficient information to resolve the $a_3$ component at each time in a way that meaningfully corresponds to the underlying dynamics. To achieve this, however, we had to use the expectation value. 
 
\subsection{\label{secFutureWork}Future work}

This work motivates future research in a number of directions. First, it would be fruitful to explore ways of improving the computational scalability of QMCl, both with regards to training and out-of-sample evaluation. To that end, methods for kernel approximation based on random features \cite{RahimiRecht07} appear well-suited to reduce the computational cost of the current brute-force kernel computations, which is $O(N^2)$ for training (basis function computation) and $O(N)$ for out-of-sample evaluation (operator-valued feature map). Random feature methods can reduce these costs to $O(N)$ and $O(1)$, respectively, while allowing streaming data processing in the training phase \cite{GiannakisEtAl21}. Another direction would be to employ methods for kernel learning \cite{OwhadiYoo19} so as to optimize the kernels used in the computation of the basis and operator-valued feature map with respect to an objective. A longer-term goal would be to develop implementations of QMCl on quantum computers. Recent approaches for simulation of dynamical systems on quantum computers based on closely-related mathematical techniques to QMCl have shown promising results for simple ergodic dynamical systems \cite{GiannakisEtAl22}, and it would be interesting to explore whether these methods can be extended in a parameterization context. Here, a challenge would be how to handle sequential two-way interactions between the classical and quantum computational systems representing the resolved and unresolved dynamics, respectively. We believe that addressing this and other related problems would be fruitful areas for future work.  

\appendix

\section{\label{appNumerics}Numerical methods}

This appendix contains details on the numerical implementation of the data-driven formulation of QMCl described in \cref{secDataDriven}. \cref{appDataDriven} provides details on the linear operators employed in the data-driven scheme. \cref{appClassical} gives the formula for classical evolution of the resolved variables used in the experiments of \cref{secL63,secL96}. \cref{appPseudocode} contains a high-level algorithmic description of the QMCl pipeline, along with associated pseudocode. \cref{appHistograms} describes the approach used to compute histograms and time-autocorrelation functions.

\subsection{\label{appDataDriven}Linear operators} 

In \cref{appDataDrivenKernel,appDataDrivenObs,appDataDrivenSpec,appDataDrivenOps,appDataDrivenEffect}, we describe the construction of the kernel integral operators, multiplication operators, spectral measures, evolution operators, and effect-valued feature maps, respectively, used in the data-driven QMCl formulation. We assume throughout availability of the training data described in \cref{secDataDriven}. For simplicity, we suppress $N$ and $L$ subscripts (number of training samples and eigenfunctions, respectively) from our notation of column vectors and matrices representing vectors and linear operators, respectively.  

\subsubsection{\label{appDataDrivenKernel}Kernel integral operators} The first step in the data-driven QMCl algorithm is to use the training data $w_0, \ldots, w_{N-1} \in \mathcal W$ to build the kernel integral operator $K_N : \hat H_N \to \hat H_N$ and compute the associated eigenfunctions $\phi_{l,N}$. The eigenvalue problem for $K_N$ in~\eqref{eqKNEig} is equivalent to the matrix eigenvalue problem
\begin{equation}
    \bm K \bm\phi_l = \lambda_{l,N} \bm\phi_l
    \label{eqKNMat}
\end{equation}
for the $N\times N$ kernel matrix $\bm K = [K_{ij}]_{i,j=0}^{N-1}$ with $K_{ij} = \kappa_{\mathcal W}(w_i,w_j)$, where $ \bm \phi_l = (\phi_{0l},\ldots,\phi_{N-1,l})^\top \in \mathbb R^N $ are column vectors whose elements give the eigenfunction values, $\phi_{ml} = \phi_{l,N}(\omega_m)$. We normalize the eigenvectors such that $\bm \phi_i^\top \bm \phi_j = N \delta_{ij}$, which is equivalent to the orthonormality condition $ \langle \phi_{i,N}, \phi_{j,N} \rangle_N = \delta_{ij}$ on $\hat H_N$. For computation, it is useful to arrange the leading $L$ eigenvectors $\bm \phi_0, \ldots, \bm\phi_{L-1} $ that form the basis of $H_{L,N}$ in an $N\times L$ matrix $\bm \Phi$ whose $l$-th column is equal to $\bm \phi_l$. 

\subsubsection{\label{appDataDrivenObs}Multiplication operators}

Every classical observable $ f : \Omega \to \mathbb C$ induces an element $\hat f_N \in L^\infty(\mu_N)$ by restriction to the finite dynamical trajectory $ \{ \omega_0, \ldots, \omega_{N-1} \} \subset \Omega$ underlying the training data. Since $L^\infty(\mu_N)$ is isomorphic to $\mathbb C^N$ whenever the states $\omega_m$ are distinct (which we assume here), we can represent $\hat f_N$ by the column vector $\bm f = (f_0, \ldots, f_{N-1})^\top \in \mathbb C^N $ where $f_m = f(\omega_m)$. We stress that $\bm f$ is empirically accessible so long as the values $f_m$ of $f$ on the states $\omega_m$ are known, without requiring knowledge of $\omega_m$. This is the case, for instance, for the values $z_m^{(i)} = Z^{(i)}(\omega_m)$ of the fluxes used as our training data. 

In QMCl, we map every such classical observable $f$ to a projected multiplication operator $A_{L,N} := \pi_{L,N} \hat f_N \in B(H_{L,N})$. Computationally, this operator is represented by the $L\times L$ matrix representation $\bm A := \bm \beta_{L,N} A_{L,N} = [ A_{ij} ]_{i,j=0}^{L-1}$ with elements $A_{ij} = \langle \phi_{i,N}, A_{L,N} \phi_{j,N} \rangle_N $. Using $ \odot$ to denote elementwise multiplication of column vectors, we have 
\begin{displaymath}
    A_{ij} = \frac{1}{N} \sum_{m=0}^{N-1} \phi_{i,N}(\omega_m) f(\omega_m) \phi_{j,N}(\omega_m) = \frac{1}{N} \bm \phi_i^\top ( \bm f \odot \bm \phi_j). 
\end{displaymath}
In matrix notation, the above expression becomes
\begin{equation}
    \label{eqAMat}
    \bm A = \frac{1}{N} \bm \Phi^\top (\diag \bm f) \bm \Phi. 
\end{equation}
Note that $\bm A $ is a self-adjoint matrix whenever $f$ is real-valued. Henceforth, we will assume that this the case. 

Given a quantum state $\rho \in Q(H_{L,N})$ with matrix representation 
\begin{equation}
    \label{eqRhoMat}
    \bm \rho = [\langle \phi_{i,N}, \rho \phi_{j,N}\rangle_N ]_{i,j=0}^{L-1},
\end{equation}
the quantum mechanical expectation $\mathbb E_{\rho} A $ from~\eqref{eqQExp}  can be computed as
\begin{equation}
    \label{eqQExpMat}
    \mathbb E_{\rho} A = \tr(\bm \rho \bm A).
\end{equation}
If $\bm \rho = \bm \xi \bm \xi^\dag$ is a pure state associated with a unit vector $\bm \xi \in \mathbb C^L$ (recall that $^\dag$ denotes the complex-conjugate transpose), then~\eqref{eqQExpMat} simplifies to
\begin{displaymath}
    \mathbb E_{\rho} A = \bm \xi^\dag \bm A \bm \xi.
\end{displaymath}

\subsubsection{\label{appDataDrivenSpec}Spectral decomposition}

In the stochastic variant of QMCl (see \cref{secStochasticQMCL}), we use the spectral measure $E_{A_{L,N}}$ of $A_{L,N} = \pi_{L,N} \hat f_N$ in order to compute the discrete probability distributions $\mathbb P_{\rho,A_{L,N}}$ in~\eqref{eqQProbF} for a given quantum state $\rho \in Q(H_{L,N})$. We stochastically generate flux terms by sampling from these distributions, as follows. 

Let $\bm u_0, \ldots, \bm u_{L-1} \in \mathbb R^L$ be a set of orthonormal eigenvectors of $\bm A$ from~\eqref{eqAMat} with corresponding eigenvalues $a_0, \ldots, a_{L-1} \in \mathbb R$ (potentially with multiplicities). In the $ \{ \phi_{l,N} \} $ basis of $H_{L,N}$, the spectral measure $E_{A_{L,N}}$ from~\eqref{eqSpecMeasF} is represented by a matrix-valued measure $\bm E : \mathcal B(\mathbb R) \to \mathbb M_L$ with $\bm E(S) = [\langle \phi_{i,N}, E_{A_{L,N}}(S) \phi_{j,N} \rangle_N]_{i,j=0}^{L-1}$, such that
\begin{displaymath}
    \bm E(S) = \sum_{l: a_l \in S} \bm u_l \bm u_l^\top.
\end{displaymath}
Note that for any eigenvalue $a_j$, the sum $\bm E^{(j)} = \sum_{l: a_l =a_j} \bm u_l \bm u_l^\top$ is the matrix representation of the projection $E^{(j)}_{A_{L,N}}$ onto the eigenspace of $A_{L,N}$ corresponding to $a_j$. 

Given a quantum state $\rho \in Q(H_{L,N})$ with matrix representation $\bm \rho$ from~\eqref{eqRhoMat}, the probability distribution $\mathbb P_{\rho,A_{L,N}}$ can be evaluated as   
\begin{equation}
    \label{eqP}
    \mathbb P_{\rho,A_{L,N}}(S) = \sum_{l:a_l \in S}  \tr\left( \bm \rho (\bm u_l \bm u^\top_l) \right) = \sum_{l:a_l \in S} \bm u_l^\top \bm \rho \bm u_l.
\end{equation}
If $\bm \rho = \bm \xi \bm \xi^\dag$ is pure, then~\eqref{eqP} simplifies to
\begin{displaymath}
    \mathbb P_{\rho,A_{L,N}}(S) = \sum_{l:a_l \in S} \lvert \bm \xi^\dag \bm u_l \rvert^2.
\end{displaymath}
Practically, we draw samples from $\mathbb P_{\rho,A_{L,N}}$ by computing the probability vector $\bm p\in \mathbb R^L$ with 
\begin{equation}
    \label{eqPVector}
    \bm p= (p_0, \ldots, p_{L-1}), \quad p_l = P_{\rho,A_{L,N}}(\{a_l\}), 
\end{equation}
and drawing samples from the spectrum $ \{ a_0, \ldots, a_{L-1} \}$ with distribution $\bm p$ using a sampling algorithm (e.g., \texttt{randsample} in MATLAB).

\subsubsection{\label{appDataDrivenOps}Evolution operators}

Following \cite{BerryEtAl15,Giannakis19b,FreemanEtAl22}, we approximate the Koopman operator by the left shift operator $\hat U_N : \hat H_N \to \hat H_N$ defined as
\begin{equation}
    \label{eqShift}
    \hat U_N f(\omega_m) = 
    \begin{cases}
        f(\omega_{m+1}), & 0 \leq m \leq N-2, \\
        0, & m = N - 1. 
    \end{cases}
\end{equation}
For completeness, we note that an alternative approach with equivalent asymptotic behavior as $N \to \infty$ (which we do not use in experiments of \cref{secL63,secL96}) is to employ the unitary shift
\begin{displaymath}
    \hat U_N f(\omega_m) = 
    \begin{cases}
        f(\omega_{m+1}), & 0 \leq m \leq N-2, \\
        f(\omega_0), & m = N - 1. 
    \end{cases}
\end{displaymath}
With either approach, we obtain the evolution operator $U_{L,N} : H_{L,N} \to H_{L,N}$ by projection onto $H_{L,N}$,  $U_{L,N} = \Pi_{L,N}\hat U_N \Pi_{L,N}$. In the $ \{ \phi_{l,N} \}$ basis of $H_{L,N}$, $U_{L,N}$ is represented by the $L\times L$ matrix 
\begin{displaymath}
    \bm U = [U_{ij}] = [ \langle \phi_{i,N}, \hat U_N \phi_j \rangle_N ]_{i,j=0}^{L-1}. 
\end{displaymath}
In particular, for the left shift operator in~\eqref{eqShift}, we have
\begin{displaymath}
    U_{ij} = \frac{1}{N} \sum_{m=0}^{N-2} \phi_i(\omega_m) \phi_j(\omega_{m+1}) = \frac{1}{N} \bm \phi_i^\top \hat{\bm U} \bm \phi_j,
\end{displaymath}
where $\hat{\bm U}$ is the $N\times N $ left shift matrix
\begin{displaymath}
    \hat{\bm U} = 
    \begin{pmatrix}
        0 & 1 \\
        & \ddots & \ddots \\
        & & 0 & 1 \\
        & & & 0
    \end{pmatrix}.
\end{displaymath}
Equivalently, in matrix notation we have 
\begin{equation}
    \label{eqUMat}
    \bm U = \frac{1}{N} \bm \Phi^\top \hat{\bm U} \bm \Phi. 
\end{equation}
Analogously to the data-independent case described in \cref{secDiscretizationOps}, $U_{L,N}$ has an induced action $\mathcal P_{L,N} : B_1(H_{L,N}) \to B_1(H_{L,N})$ defined by the conjugation formula $\mathcal P_{L,N} A = U_{L,N}^* A U_{L,N}$. (Note that, as a vector space, $B_1(H_{L,N})$ is isomorphic to $B(H_{L,N})$ by finite dimensionality of $H_{L,N}$, but $B_1(H_{L,N})$ is equipped with the trace norm, whereas $B(H_{L,N})$ is equipped with the operator norm.) In the $\{ \phi_{l,N} \}$ basis of $H_{L,N}$, $\mathcal P_{L,N}$ is represented by the linear operator $\mathsf P : \mathbb M_L \to \mathbb M_L $ on $L \times L $ matrices defined as $\mathsf P \bm A = \bm U^\top \bm A \bm U$. That is, if $\bm A = \bm \beta_{L,N} A \in \mathbb M_L $ is the matrix representation of an operator $A \in B_1(H_{L,N})$, then $\mathsf P \bm A $ is the matrix representation of $\mathcal P_{L,N}A$, i.e., $ \mathsf P \bm A = \bm \beta_{L,N} (\mathcal P_{L,N} A)$.       

In general, $U_{L,N}$ is not a unitary operator, so $\mathcal P_{L,N}$ does not necessarily map quantum states in $Q(H_{L,N}) \subset B_1(H_{L,N})$ to quantum states (see \cref{secDiscretizationOps}). However, $ \mathcal P_{L,N}$ is trace non-increasing, so it generates an open quantum system. In our computations, we enforce state preservation by replacing $\mathcal P_{L,N}$ by the nonlinear map
\begin{displaymath}
    \tilde{\mathcal P}_{L,N}(\rho) := \frac{\mathcal P_{L,N}\rho}{\tr(\mathcal P_{L,N}\rho)},
\end{displaymath}
defined on the set of quantum states $ \rho \in Q(H_{L,N})$ for which $\tr(\mathcal P_{L,N} \rho)$ is nonzero. This map is represented by the map  $\tilde{\mathsf P}$ on $L\times L$ density matrices defined as
\begin{displaymath}
    \tilde{\mathsf P}(\bm \rho) = \frac{\mathsf P \bm \rho}{\tr(\mathsf P \bm \rho)}.
\end{displaymath}
If $\bm \rho = \bm \xi \bm \xi^\dag$ is pure, then $\tilde{\mathsf P}(\bm \rho) = \tilde{\bm \xi}^\dag \tilde{\bm \xi} $ is a pure state associated with the unit vector 
\begin{displaymath}
    \tilde{\bm \xi} = \tilde{\bm P}(\bm \xi) := \frac{\bm U^\top \bm \xi}{\lVert \bm U^\top \bm \xi \rVert_2}.
\end{displaymath}
Thus, to compute the evolution of pure states it is sufficient to work with the map $\tilde{\bm P}$ on state vectors rather than explicitly with $\tilde{\mathsf P}$ on density matrices---this results in a reduction of computational cost from $O(L^3)$ to $O(L^2)$; see \cref{secCost}.

\subsubsection{\label{appDataDrivenEffect}Effect-valued feature maps} 

Recall from \cref{secFeatureMap,secDataDrivenQMCl} that the effect-valued feature map $\mathcal F_{L,N} : \mathcal X \to \mathcal E(H_{L,N})$ maps each state $x \in \mathcal X$ of the resolved variables to a projected multiplication operator $ \pi_{L,N}(\hat F_N(x))$, where $\hat F_N(x) = k(x, X(\cdot)) \in L^\infty(\mu_N)$ is the feature vector associated with the kernel $k : \mathcal X \times \mathcal X \to [0,1]$. In the $\{ \phi_{l,N} \}$ basis of $H_{L,N}$, $ \mathcal F_{L,N} $ is represented by the matrix-valued map $\bm F : \mathcal X \to \mathbb M_L$ with $\bm F = \bm\beta_{L,N} \circ \mathcal F_{L,N}$. That is, we have $\bm F(x) = \bm A$, where $\bm A $ is given by~\eqref{eqAMat} with $\bm f = ( k(x,x_0), \ldots, k(x,x_{N-1}))^\top$. State conditioning by the effect $\mathcal F_{L,N}(x)$ (see~\eqref{eqQBayes2} and~\eqref{eqQMPosterior}) can then be computed via the matrix formula
\begin{equation}
    \bm \rho|_{\bm F(x)} = \frac{\sqrt{\bm A} \bm \rho \sqrt{\bm A}}{\tr(\sqrt{\bm A} \bm \rho \sqrt{\bm A})},
    \label{eqQMPosteriorMat}
\end{equation}
where $\bm \rho|_{\bm F(x)}$ is the matrix representation of $\rho|_{\mathcal F_{L,N}(x)}$. If $\bm \rho = \bm \xi \bm \xi^\dag$ is pure, then the conditioned state $\bm \rho |_{\bm F(x)} $ is also pure, and the associated state vector is given by
\begin{equation}
    \label{eqQMPosteriorPure}
    \bm \xi|_{\bm F(x)} = \frac{\sqrt{\bm A}\bm \xi}{\lVert \sqrt{\bm A}\bm \xi \rVert_2}.
\end{equation}

A drawback of~\eqref{eqQMPosteriorMat} and~\eqref{eqQMPosteriorPure} is that they require the computation of the matrix square root $\sqrt{\bm A}$. To avoid the cost of this step, we can modify $\mathcal F_{L,N}$ to the effect-valued map $ \tilde{\mathcal F}_{L,N} (x) = (\pi_{L,N} \hat F^{1/2}_N(x))^2$. This map is represented by the matrix-valued function $\tilde{\bm F} : \mathcal X \to \mathbb M_L $ such that $\tilde{\bm F}(x) = \bm A^2$, where $\bm A$ is given by~\eqref{eqAMat}, now with $\bm f = (\sqrt{k(x,x_0)}, \ldots, \sqrt{k(x,x_{N-1})})^\top$. In this case, the conditioning formula by $\tilde{\bm F}(x)$ becomes 
\begin{equation}
    \bm \rho|_{\tilde{\bm F}(x)} = \frac{\bm A \bm \rho \bm A}{\tr(\bm A \bm \rho \bm A)},
    \label{eqQMPosteriorMat2}
\end{equation}
which avoids the matrix square root. Analogously to~\eqref{eqQMPosteriorPure}, the state vector update under~\eqref{eqQMPosteriorMat2} when $\bm \rho = \bm \xi \bm \xi^{\dag}$ is pure becomes
\begin{equation}
    \label{eqQMPosteriorPure2}
    \bm \xi|_{\tilde{\bm F}(x)} = \frac{\bm A\bm \xi}{\lVert \bm A\bm \xi \rVert_2},
\end{equation}
which again does not require computation of a matrix square root. Our experiments in \cref{secL63,secL96} utilize~\eqref{eqQMPosteriorPure2} for state conditioning. 

Even though, in general, $\mathcal F_{L,N}(x)$ and $\tilde{\mathcal F}_{L,N}(x)$ are not equal, the two maps have the same asymptotic limit as $N\to\infty$ and $L\to\infty$ (see \cref{secDiscretizationEffect,secDataDrivenConvergence}). In that limit, conditioning by any of $\mathcal F_{L,N}(x)$ and $\tilde{\mathcal F}_{L,N}(x)$ recovers conditioning by $\mathcal F(x)$ in the infinite-dimensional quantum system on $H$. The latter, is in turn consistent with classical Bayesian conditioning by the feature vectors $F(x)$ (see \cref{secEmbeddingEffects}). We note that for the radial basis function kernel in \eqref{eqKGauss} used in our experiments the square root kernel function $\sqrt k$ can be obtained simply by scaling the bandwidth parameter $\epsilon$ by a factor of $\sqrt 2$.  

\subsection{\label{appClassical}Classical evolution}

In the experiments of \cref{secL63,secL96}, the classical evolution map $\tilde \phi : \mathcal X \times \mathcal Z \to \mathcal X$ is based on a standard RK4 discretization of the resolved component of the dynamics on $\mathcal X = \mathbb R^{d_{\mathcal X}}$, keeping the flux terms in $\mathcal Z = \mathbb R^d$ fixed. Specifically, given that the resolved components of the dynamics satisfy $\dot x(t) = v(x(t), z(t))$ for a vector field $v: \mathcal X \times \mathcal Z \to \mathcal X$, and the timestep of the parameterized system is $\Delta t$, we set  
\begin{equation}
    \label{eqRK4}
    \begin{gathered}
        \tilde\phi(x,z) = x + \frac{1}{6}(k_1 + 2 k_2 + 2 k_3 + k_4) \, \Delta t,\\
        \begin{aligned}
            k_1 = v( x, z), \quad k_2 = v( x + \Delta t\, k_1 / 2, z), \\
            k_3 = v( x + \Delta t \, k_2 / 2, z), \quad k_4 = v(x + \Delta t\, k_3, z).
        \end{aligned}
    \end{gathered}
\end{equation}

\subsection{\label{appPseudocode}Algorithm structure}

The data-driven QMCl pipeline consists of training, initialization, and simulation/prediction phases, summarized in \cref{appAlgTraining,appAlgInitialization,appAlgClosure} respectively.

\subsubsection{\label{appAlgTraining}Training}

Recall that our training dataset consists of time-ordered samples $w_0, \ldots, w_{N-1} \in \mathcal W$ for computing basis functions, and samples $z_0^{(i)}, \ldots, z_{N-1}^{(i)} \in \mathbb R$ with $ i \in \{ 1, \ldots, d \} $ of the components of the flux $Z : \mathcal Y \to \mathcal Z$. In this paper, we assume that the samples are taken on a single dynamical trajectory $\omega_0, \ldots, \omega_{N-1} \in \Omega$, i.e., $w_m = W(\omega_m) $, $z_m^{(i)} = Z^{(i)}(\omega_m)$, and $ \omega_m = \Phi^m(\omega_0)$ for some initial condition $\omega_0 \in \Omega$. The methods described below can be readily generalized to training with samples from ensembles of shorter trajectories, so long as the sampling measure $\mu_N$ of the data converges to the invariant measure in the sense of~\eqref{eqWeakConv}. The steps of the training phase are as follows:     

\begin{enumerate}
    \item Tune the bandwidth parameter $\epsilon_{\mathcal W}$ of the kernel $\kappa_{\mathcal W}$. We perform this step using the automatic tuning procedure developed in refs.~\cite{CoifmanEtAl08,BerryHarlim16}. The procedure is based on a grid search in a collection of candidate $\epsilon_{\mathcal W}$ values, selecting the value that maximizes a kernel-dependent measure of dimension of the dataset. For pseudocode, see, e.g., Algorithm~B.5 in \cite{FreemanEtAl22}.
    \item With the bandwidth parameter from Step~1, form the kernel matrix $\bm K$ and compute the kernel eigenvectors $\{ \bm \phi_0, \ldots, \bm \phi_{L-1} \} $ from~\eqref{eqKNMat}. Arrange the basis vectors in a matrix $\bm\Phi$ as described in~\cref{appDataDrivenKernel}. In applications, we typically approximate $\bm K$ by a sparse matrix using nearest-neighbor truncation, and solve the resulting eigenvalue problem with iterative solvers (e.g., MATLAB's \texttt{eigs}). 
    \item For each $i \in \{ 1, \ldots d \}$, use $\bm \Phi$ and the $z_m^{(i)}$ samples to compute an $L\times L$ multiplication operator matrix $\bm Z^{(i)}$ via~\eqref{eqAMat} with $ \bm f = (z_0^{(i)}, \ldots, z_{N-1}^{(i)})^\top$. In the stochastic variant of QMCl, we also compute the eigenvalues $a^{(i)}_0, \ldots, a^{(i)}_{L-1} \in \mathbb R$ and corresponding eigenvectors $ \bm u_0^{(i)}, \ldots, \bm u_{L-1}^{(i)}$ for each $\bm Z^{(i)}$. 
    \item Compute the $L\times L$ Koopman operator matrix $\bm U$ from \eqref{eqUMat}.
\end{enumerate}
 
\subsubsection{\label{appAlgInitialization}Initialization}
Let $\hat x_0 \in \mathcal X$ be a given initial classical state and $\hat\rho_0 \in Q(H_{L,N})$ an initial quantum state that may depend on $\hat x_0$. In the experiments of \cref{secL63,secL96}, we set $\hat \rho_0$ to the uninformative state $\bar\rho_{L,N}$ from~\eqref{eqInit2}. This is a pure state induced by the unit vector $\bar \xi_{L,N} = \bar \psi_{L,N} / \lVert \bar \psi_{L,N} \rVert_{\hat H_N} $, where $ \bar \psi_{L,N} = \sum_{l=0}^{L-1} c_l \phi_{l,N}$ and $c_l = \langle 1_\Omega, \phi_{l,N} \rangle_N $. In column vector notation, the expansion coefficients $c_l$ are given by $c_l = \bm \phi_l^\top \bm 1_N / N$, where $\bm 1_N = (1,\ldots, 1)^\top $ is the column vector in $\mathbb C^N$ whose all elements are equal to 1. Correspondingly, in the $\{ \phi_{l,N} \}$ basis of $H_{L,N}$ the state vector $\bar \xi_{L,N}$ is represented by the column vector $\bar{\bm \xi} = \beta_{L,N} \bar\xi_{L,N}$ given by
\begin{displaymath}
    \bar{\bm \xi} = \frac{\bm c}{ \lVert \bm c \rVert_2}, \quad \bm c = (c_0, \ldots, c_{L-1} )^\top,
\end{displaymath}
and the density operator $\bar\rho_{L,N}$ is represented by the rank-1 density matrix $\bm{\bar \rho}_{L,N} = \bar{\bm \xi} \bar{\bm \xi}^\top$. Note that if the kernel $\kappa_{\mathcal W}$ used to build the basis is normalized to a Markov kernel (see, e.g., the QMCl experiments in \cref{secL63Small} and the QMDA experiments in refs.~\cite{Giannakis19b,FreemanEtAl22}), then the leading basis element $\bm \phi_0$ can be chosen as $\bm 1_N$, and $\bar{\bm \xi}$ simplifies to the vector $(1,0,\ldots, 0)$. 

Besides the choice $\hat \rho_0 = \bar\rho_{L,N}$, an alternative approach to quantum state initialization (which we do not employ in the experiments presented in this paper) is to use the feature map $\mathcal F_{L,N}$ to set $\hat \rho_0 = \mathcal F_{L,N}(\hat x_0) / \tr(\mathcal F_{L,N}(\hat x_0))$. As noted in \cref{secInit}, $\hat \rho_0$ obtained with this approach is not, in general, a pure state.   

\subsubsection{Simulation/prediction\label{appAlgClosure}}
Given the initial data $(\hat x_0, \hat{\bm \rho}_0)$ obtained by any of the two methods described in \cref{appAlgInitialization}, the QMCl parameterized system evolves by alternating between classical and quantum evolution, as described in \cref{secStateUpdate,secDataDrivenQMCl}. \cref{algQMCl} implements this evolution over $r$ timesteps given that the state of the system at the $n$-th timestep is $(\hat x_n, \hat{\bm \rho}_n)$. We recall that $r \in \mathbb N$ is the number of timesteps between each update of the quantum state by the feature map (see \cref{secStateUpdate}). Thus, the output of \cref{algQMCl} is a sequence of classical states $\hat x_{n+1}, \ldots, \hat x_{n+r} \in \mathcal X$ at timesteps $n+1, \ldots, n + r$, respectively, and a posterior quantum state $\hat{\bm \rho}_{n+r} $ at timestep $n+r$. To continue to march the system forward, \cref{algQMCl} is executed with initial conditions $(\hat x_{n+r}, \hat{\bm \rho}_{n+r})$.

\begin{algorithm}
    \caption{\label{algQMCl}QMCl forecast--analysis cycle. The algorithm assumes that training has been performed as described in \cref{appAlgTraining}.}
    \Inputs
    \begin{itemize}
        \item Resolved variables $\hat x_n \in \mathcal X$ and density matrix $\hat{\bm \rho}_n \in \mathbb M_L$ at the $n$-th timestep.
        \item Number of timesteps $r \in \mathbb N$ per quantum Bayesian update.
    \end{itemize}
    \Outputs
    \begin{itemize}
        \item Resolved variables $\hat x_{n+1}, \ldots, \hat x_{n+r} \in \mathcal X$ at timesteps $n+1, \ldots, n+r$.
        \item Posterior density matrix $\hat{\bm \rho}_{n+r} \in \mathbb M_L$ at timestep $n+r$.
    \end{itemize}
    \Steps
    \begin{enumerate}
        \item Set $ \tilde{\bm \rho}_n = \hat{\bm \rho}_n$.
        \item For $j \in \{ 1, \ldots, r \}$:
            \begin{enumerate}
                \item Compute the fluxes $ z_{n+j-1} = (z_{n+j-1}^{(1)}, \ldots, z_{n+j-1}^{(d)}) \in \mathbb R^d$ with $ z_{n+j-1} = \tr(\hat{\bm \rho}_{n+j-1} \bm Z^{(i)}) $. 
                \item Update the resolved variables: $\hat x_{n+j} = \tilde \phi(\hat x_{n+j-1},z_{n+j-1})$.
                \item Update the density matrix with the transfer operator: $\tilde{\bm \rho}_{n+j} = \tilde{\mathsf P}(\tilde{\bm \rho}_{n+j-1})$.
            \end{enumerate}
        \item Evaluate the effect-valued feature map: $\bm e_{n+r} = \tilde{\bm F}(\hat x_{n+r})$.
        \item Compute the posterior density matrix $\hat{\bm \rho}_{n+r} = \tilde{\bm \rho}_{n+r} |_{\bm e_{n+r}}$ via~\eqref{eqQMPosteriorMat2}.
        \item \Return $\hat x_{n+1}, \ldots, \hat x_{n+r}, \hat{\bm \rho}_{n+r} $.
    \end{enumerate}
\end{algorithm}

If the initial state $\hat{\bm \rho}_0 = \hat{\bm \xi}_0 \hat{\bm \xi}_0^\dag$ is pure, then under iteration of \cref{algQMCl} all subsequent states $\hat{\bm \rho}_r, \hat{\bm \rho}_{2r}, \ldots $ are pure. \cref{algQMClPure} specializes \cref{algQMCl} to that setting, which results in a reduction of computational cost from $O(L^3)$ to $O(L^2)$ (see \cref{secCost}).

\begin{algorithm}
    \caption{\label{algQMClPure}QMCl forecast--analysis cycle for pure states. The algorithm assumes that training has been performed as described in \cref{appAlgTraining}.}
    \Inputs
    \begin{itemize}
        \item Resolved variables $\hat x_n \in \mathcal X$ and state vector $\hat{\bm \xi}_n \in \mathbb C^L$ at the $n$-th timestep.
        \item Number of timesteps $r \in \mathbb N$ per quantum Bayesian update.
    \end{itemize}
    \Outputs
    \begin{itemize}
        \item Resolved variables $\hat x_{n+1}, \ldots, \hat x_{n+r} \in \mathcal X$ at timesteps $n+1, \ldots, n+r$.
        \item Posterior state vectors $\hat{\bm \xi}_{n+r} \in \mathbb M_L$ at timestep $n+r$.
    \end{itemize}
    \Steps
    \begin{enumerate}
        \item Set $ \tilde{\bm \xi}_n = \hat{\bm \xi}_n$.
        \item For $j \in \{ 1, \ldots, r \}$:
            \begin{enumerate}
                \item Compute the fluxes $ z_{n+j-1} = (z_{n+j-1}^{(1)}, \ldots, z_{n+j-1}^{(d)}) \in \mathbb R^d$ with $ z_{n+j-1} =\tilde{\bm \xi}_{n+j-1}^\dag \bm Z^{(i)} \tilde{\bm \xi}_{n+j-1}$. 
                \item Update the resolved variables: $\hat x_{n+j} = \tilde \phi(\hat x_{n+j-1},z_{n+j-1})$.
                \item Update the state vector with transfer operator: $\tilde{\bm \xi}_{n+j} = \tilde{\bm P}(\tilde{\bm \xi}_{n+j-1})$.
            \end{enumerate}
        \item Evaluate the effect-valued feature map: $\bm e = \tilde{\bm F}(\hat x_{n+r})$.
        \item Compute the state vector $\hat{\bm \xi}_{n+r} = \tilde{\bm \xi}_{n+r} |_{\bm e_{n+r}}$ via~\eqref{eqQMPosteriorPure2}.
        \item \Return $\hat x_{n+1}, \ldots, \hat x_{n+r}, \hat{\bm \xi}_{n+r} $.
    \end{enumerate}
\end{algorithm}

Finally, \cref{algQMClStoch} implements the stochastic variant of QMCl for general (mixed) states. The specialization of this algorithm to pure states is entirely analogous to \cref{algQMClPure} so we do not include it here in the interest of brevity.

\begin{algorithm}
    \caption{\label{algQMClStoch}Stochastic QMCl forecast--analysis cycle. The algorithm assumes that training has been performed as described in \cref{appAlgTraining}.}
    \Inputs
    \begin{itemize}
        \item Resolved variables $\hat x_n \in \mathcal X$ and density matrix $\hat{\bm \rho}_n \in \mathbb M_L$ at the $n$-th timestep.
        \item Number of timesteps $r \in \mathbb N$ per quantum Bayesian update.
    \end{itemize}
    \Outputs
    \begin{itemize}
        \item Resolved variables $\hat x_{n+1}, \ldots, \hat x_{n+r} \in \mathcal X$ at timesteps $n+1, \ldots, n+r$.
        \item Posterior density matrix $\hat{\bm \rho}_{n+r} \in \mathbb M_L$ at timestep $n+r$.
    \end{itemize}
    \Steps
    \begin{enumerate}
        \item Set $ \tilde{\bm \rho}_n = \hat{\bm \rho}_n$.
        \item For $j \in \{ 1, \ldots, r \}$:
            \begin{enumerate}
                \item For $i \in \{1, \ldots, d \}$:
                    \begin{enumerate}
                        \item Compute the probability vector $\bm p_{n+j-1}^{(i)} \in \mathbb R^L$ associated with $\tilde{\bm\rho}_{n+j-1} $ and $\bm Z^{(i)}$ using~\eqref{eqPVector}.   
                        \item Draw a sample $z_{n+j-1}^{(i)}$ from the spectrum $\{ a^{(i)}_0, \ldots, a^{(i)}_{L-1} \}$ of $\bm Z^{(i)}$ with distribution $\bm p_{n+j-1}^{(i)}$.  
                    \end{enumerate}
                \item Set the flux term $z_{n+j-1} = (z_{n+j-1}^{(1)}, \ldots, z_{n+j-1}^{(d)}) \in \mathbb R^d$, and update the resolved variables: $\hat x_{n+j} = \tilde \phi(\hat x_{n+j-1},z_{n+j-1})$.
                \item Update the density matrix with the transfer operator: $\tilde{\bm \rho}_{n+j} = \tilde{\mathsf P}(\tilde{\bm \rho}_{n+j-1})$.
            \end{enumerate}
        \item Evaluate the effect-valued feature map: $\bm e_{n+r} = \tilde{\bm F}(\hat x_{n+r})$.
        \item Compute the posterior density matrix $\hat{\bm \rho}_{n+r} = \tilde{\bm \rho}_{n+r} |_{\bm e_{n+r}}$ via~\eqref{eqQMPosteriorMat2}.
        \item \Return $\hat x_{n+1}, \ldots, \hat x_{n+r}, \hat{\bm \rho}_{n+r} $.
    \end{enumerate}
\end{algorithm}

\subsection{\label{appHistograms}Histograms and autocorrelation functions}

For time-ordered data $f_0, \ldots, f_{N-1} \in \mathbb R$ sampled at a fixed interval $\Delta t$, the value of the time-autocorrelation function $C_{f,N}(\tau_j)$ at timeshift value $\tau_j = j\, \Delta t$, $j \in \mathbb N$, is given by 
\begin{displaymath}
    C_{f,N}(\tau_j) = \frac{1}{N} \sum_{n=0}^{N-j} f_n f_{n+j}.  
\end{displaymath}
In the main text, we show plots of the normalized autocorrelation function $\bar C_{f,N}(\tau_j) := C_{f,N}(\tau_j) / C_{f,N}(0)$, where $\bar C_{f,N}(0) = 1$ by construction.
If the $f_n$ are samples of an observable $ f \in L^2(\mu)$ taken on an orbit $\omega_0, \ldots, \omega_{N-1}$ of the dynamics, i.e., $f_n = f(\omega_n)$ and $ \omega_n = \Phi^n(\omega_0)$, then by the pointwise ergodic theorem, as $N \to \infty$, $C_{f,N}(\tau_j)$ converges for $\mu$-a.e.\ $\omega_0 \in \Omega$ to $C_f(\tau_j) = \langle f, U^j f \rangle$, where $U : H \to H$ is the Koopman operator induced by $\Phi$.   

To generate histograms based on the data $f_0, \ldots, f_{N-1}$, we split the interval $ [\min \{ f_n \}_{n=0}^{N-1}, \max \{ f_n \}_{n=0}^{N-1} ]$ into $B$ uniformly-sized bins $S_0, \ldots, S_{B-1}$, and compute the normalized counts $(N_0/N, \ldots, N_{S-1}/ N) $ where $N_i$ is the number of datapoints $f_n$ lying in $S_i$. For our graphs, the value $B=45$ was chosen.

\bibliographystyle{siamplain}

\end{document}